\documentclass[12pt]{article}
\usepackage{latexsym}
\usepackage{amssymb}
\usepackage{wasysym}
\usepackage{mathrsfs}
\DeclareMathAlphabet{\zap}{OT1}{pzc}{m}{it}
\DeclareMathAlphabet{\mathzap}{OT1}{pzc}{m}{it}
\usepackage[OT2,OT1]{fontenc}
\usepackage{graphicx}
\DeclareGraphicsExtensions{.eps, .jpg}
\def\cyr{%
\renewcommand\rmdefault{wncyr}%
\renewcommand\sfdefault{wncyss}%
\renewcommand\encodingdefault{OT2}%
\normalfont
\selectfont}
\DeclareTextFontCommand{\textcyr}{\cyr}

\def\half{{\textstyle \frac{1}{2}}}
\def\hook{~\lrcorner ~}
\newenvironment{romenumerate}{\begin{list}{$(\roman{enumi})$}
{\usecounter{enumi}\setlength{\rightmargin}{\leftmargin}}}{\end{list}}
\newenvironment{alphenumerate}{\begin{list}{$(\alph{enumi})$}
{\usecounter{enumi}\setlength{\rightmargin}{\leftmargin}}}{\end{list}}
\newcommand{\wot}{{\mho}}

\newcommand{\Dye}{\mbox{\cyr   D}}
\newcommand{\Bye}{\mbox{\cyr   B}}
\newcommand{\bye}{\mbox{\cyr   b}}
\newcommand{\fye}{\mbox{\cyr   f}}
\newcommand{\Pye}{\mbox{\cyr  \em  P}}

\newcommand{\CC}{\mathbb{C}}

\newcommand{\CP}{\mathbb C\mathbb P}
\newcommand{\RP}{\mathbb R\mathbb P}
\newcommand{\PP}{\mathbb P}
\newcommand{\RR}{\mathbb R}

\newcommand{\ZZ}{\mathbb{Z}}

\newcommand{\D}{\mathrm{D}}

\renewcommand{\O}{\mathcal{O}}
\renewcommand{\P}{\mathbb{P}}

\newcommand{\zz}{{\mathfrak z}}

\newtheorem{main}{Theorem}

\newtheorem{thm}{Theorem}[section]
\newtheorem{defn}[thm]{Definition}
\newtheorem{prop}[thm]{Proposition}
\newtheorem{cor}[thm]{Corollary}
\newtheorem{lem}[thm]{Lemma}

\newenvironment{remark}{\medskip \noindent {\bf Remark.}}{\hfill $\diamondsuit$\\}
\newenvironment{ack}{\medskip \noindent {\bf Acknowledgment.}}{\hfill ~\\}
\newenvironment{proof}{\medskip \noindent {\bf Proof.}}{\hfill $\blacksquare$\\}
\begin{document}
\title{Nonlinear Gravitons,  Null Geodesics, \\
 and  Holomorphic Disks}
  
\author{Claude LeBrun\thanks{Supported 
in part by  NSF grant DMS-0305865.}     ~and  
L.J. Mason}

 \date{April 14, 2005}

\maketitle

\abstract{
We develop a  global twistor correspondence 
for  pseudo-Riemannian conformal structures 
of signature $({+}{+}{-}{-})$ 
 with self-dual
  Weyl  curvature. 
Near the conformal class of the standard indefinite product metric
 on $S^2 \times S^2$, there is
an infinite-dimensional moduli space of such conformal structures, 
and each of these has the surprising  global property that its null
geodesics are all periodic.  
Each such 
 conformal structure arises  from   a
family of holomorphic disks in  $\CP_3$
with boundary on some totally real embedding of  $\RP^3$ into $\CP_3$.  
An  interesting   sub-class  of these conformal structures are 
represented by  scalar-flat  indefinite K\"ahler metrics, and 
our methods give particularly sharp results in this 
more restrictive setting.}

\section{Introduction}

Twistor correspondences, as pioneered by Roger Penrose \cite{pnlg}, 
provide a way of understanding certain differential geometries as 
fundamentally arising from moduli  spaces of compact complex curves in  a complex 
manifold. It  has  emerged  only recently, however,   that an 
analogous pattern of phenomena can also be expected to arise from moduli spaces of 
compact complex curves-with-boundary in a complex manifold,
where the boundaries of the curves are constrained to lie in a maximal totally real
submanifold. Our previous work in this direction  \cite{lmzoll}  focused on spaces of  holomorphic
disks in $\CP_2$, with boundaries on a totally real embedding of $\RP^2$. 
In the present article, we will  see that a similarly rich geometric  story  arises 
from the moduli space of holomorphic disks in $\CP_3$ with boundaries
on a totally real $\RP^3$.

Penrose's original twistor correspondence,
which he called the {\sl nonlinear graviton}, 
hinged on the idea  that 
 {self-dual} conformal 
  metrics on $4$-manifolds tend to 
   arise from suitable holomorphic families of  $\CP_1$'s  in 
  complex $3$-manifolds. 
Penrose's formulation of these ideas  
involved    local analytic continuations of real-analytic geometries into the 
  complex domain, thereby making   the metric signature    essentially 
  irrelevant. Nonetheless,  it was  specifically the positive-definite  realm 
  of  Riemannian
geometry that  
witnessed 
 the most intensive  subsequent  
cultivation  of these ideas,  a development largely attributable to   
  the elegant and  definitive global
Riemannian 
reformulation of  the 
Penrose correspondence discovered by     Atiyah, Hitchin, and Singer \cite{AHS}.
   By contrast, however, the present article will     focus entirely  on 
 $4$-manifolds with 
  {\em split-signature} metrics, meaning pseudo-Riemannian metrics 
  of signature $(++--)$; these  have elsewhere been called
  {\em neutral metrics} \cite{law}, and are characterized by the 
  fact that they have components 
$$ \left[\begin{array}{cccc}+1 &  &  &  \\ & +1 &  &  \\ &  & -1 &  \\ &  &  & -1\end{array}\right]$$
 in a suitably chosen basis for any given tangent space. 
 What we will develop here is  a global twistor correspondence 
 for self-dual split-signature $4$-manifolds in  which 
 {\em every null geodesic is a  simple closed curve}. 
Such metrics will turn out to naturally arise  as moduli spaces for holomorphic
disks in $\CP_3$  with boundary on a fixed  totally real submanifold. 
 
As in the Riemannian case, a  split-signature metric $g$ 
on an oriented $4$-manifold $M$ is said to be 
{\em self-dual} if its Weyl (or conformal) curvature tensor, considered as a bundle-valued
$2$-form, is its own Hodge star; cf.  \S \ref{selfdual} 
below. This  is a conformally invariant condition, and should therefore
primarily 
be thought of as a constraint on the conformal class
$[g]=\{ {\zap f} g~|~{\zap f}\neq 0\}$ of the metric. Notice that
any locally conformally flat  split-signature metric on an oriented
$4$-manifold is automatically self-dual. 

For us, the protypical  example  is 
 the indefinite product metric $$g_0=\pi_1^*h-\pi_2^*h$$
on $S^2\times S^2$,
where $\pi_1 , \pi_2 : S^2\times S^2\to S^2$ are the two
factor projections and  $h$ is the standard homogeneous  metric on 
$S^2$.  This metric is actually 
conformally flat, since, thinking of  $S^2\times S^2$ as the locus
$$x_1^2+x_2^2+x_3^2=1, ~~ y_1^2+y_2^2+y_3^2=1,$$
in $\RR^{3,3}= \RR^3\times \RR^3$,  and introducing  `stereographic' coordinates  by \begin{equation}
\label{stereo}
\begin{array}{ccccccc}{\mathfrak x}_1 & = &\displaystyle 
 \frac{x_1}{2(x_3-y_3)} & ~~~~ & {\mathfrak x}_2 & = & \displaystyle \frac{x_2}{2(x_3-y_3)}
 \\ &&&&&&\\
 {\mathfrak y}_1 & = &\displaystyle  \frac{y_1}{2(x_3-y_3)}  & ~~~~ & {\mathfrak y}_2 & = &\displaystyle 
  \frac{y_2}{2(x_3-y_3)} ~,\end{array}
\end{equation}
 the  metric can be re-expressed in the form 
$$g_0 = \frac{d{\mathfrak x}_1^2+ d{\mathfrak x}_2^2-d{\mathfrak y}_1^2-d{\mathfrak y}_2^2}{{\mathfrak x}_1^2+
{\mathfrak x}_2^2+ [{\mathfrak y}_1^2+{\mathfrak y}_2^2-{\mathfrak x}_1^2-{\mathfrak x}_2^2+ \frac{1}{4}]^2} ~.$$

However, this example has a second fundamental property that will
play a crucial r\^ole  in this paper. 
Indeed,  the null geodesics of 
$(S^2\times S^2, g_0)$  are all {\em embedded circles}, 
since each is obtained   
  by simultaneously traversing 
 a great circle in  each  $S^2$  with  equal speed. 
Following Guillemin \cite{mogul}, 
we  will use  the word {\em Zollfrei}  to describe 
pseudo-Riemannian metrics with this property;  
for a detailed discussion, see \S \ref{zollfrei} below.
The Zollfrei condition is also 
conformally invariant, so that  we may consider it as yet another 
 property of the conformal class $[g]$.

 Among all split-signature metrics on a given manifold, the Zollfrei condition 
 is highly non-generic. It may therefore seem surprising that it becomes 
 an  {\em open}
 condition when restricted to the subspace of self-dual metrics: 

\begin{main}\label{auto}
Let $(M,g)$ be a self-dual Zollfrei $4$-manifold. Then,
with respect to the $C^2$ topology, there is 
  an open  neighborhood of $g$ 
  in the space of 
pseudo-Riemannian
metrics on $M$
  such that every  self-dual metric contained in this neighborhood is 
also Zollfrei. 
\end{main}
For the purpose of studying the moduli  of self-dual conformal structures,
it thus seems reasonable    to focus for the present 
on  understanding those   
self-dual metrics which are also Zollfrei. 

But this point of view  immediately prompts us to ask, 
``Which $4$-manifolds  admit self-dual Zollfrei metrics?''
 We have just seen that $S^2\times S^2$ is one such manifold. 
Another example is given by the 
 projective quadric
$${\mathbb M}^{2,2}= \Big\{ [x_1: x_2 : x_3 : y_1: y_2: y_3 ] 
\in \RP^5~\Big|~ |\vec{x}|^2 - |\vec{y}|^2 = 0\Big\}~,$$ 
which may be viewed as
 the  quotient
of $(S^2\times S^2, g_0)$
by the isometric $\ZZ_2$-action generated by the double antipodal map 
$$(\vec{x},\vec{y}) \mapsto (-\vec{x},-\vec{y}).$$
However,  we will show in \S \ref{prelim} that these are the {only} topological 
possibilities: 

\begin{main}\label{topprop}
Let $(M,g)$ be a connected oriented split-signature $4$-manifold which is 
both Zollfrei and self-dual. Then $M$ is homeomorphic to  
either $S^2\times S^2$ or  ${\mathbb M}^{2,2}$. 
\end{main}

This  topological rigidity, however,  is by no means  symptomatic of any 
kind of underlying {\sl geometric} rigidity. 
 To the contrary,  our central purpose here is  
to prove the following {\em flexibility} result: 

\begin{main} 
\label{zfsd}
There is a natural one-to-one correspondence between 
\begin{itemize}
\item equivalence classes of 
smooth self-dual  split-signature conformal
structures on $S^2\times S^2$;  and
\item  equivalence classes of  totally real embeddings $\RP^3\hookrightarrow \CP_3$,
\end{itemize} at least in a
neighborhood of the standard conformal metric $[g_0]$ and 
the standard  embedding of $\RP^3$. 
\end{main}

 Here, two conformal structures are   considered to be equivalent  iff one
is the pull-back of the other via some orientation-preserving 
self-diffeomorphism of $S^2\times S^2$;  two embeddings 
 $\RP^3\hookrightarrow \CP_3$ are considered to be equivalent iff 
 they are interrelated by a reparameterization of $\RP^3$ and/or 
 the action of  $PSL(4,\CC)$ on $\CP_3$. 
 In particular, the moduli space of self-dual Zollfrei conformal
structures on $S^2\times S^2$ is infinite-dimensional; and, roughly speaking, the 
general such conformal structure depends on  $3$ free  functions of $3$  variables. 
The correspondence between the two kinds of structures depends 
on the existence of an $(S^2\times S^2)$-family of holomorphic 
disks with boundary on a given totally real $\RP^3\subset \CP_3$.

By contrast, the same arguments also show that $({\mathbb M}^{2,2}, [g_0])$
has no non-trivial  self-dual deformations. Indeed, by analogy 
with the Blaschke conjecture \cite{beszoll,lmzoll}, we are tempted to speculate   
 that, up to conformal isometry,  $({\mathbb M}^{2,2}, [g_0])$ might well  be the 
 only non-simply-connected self-dual Zollfrei   $4$-manifold.

Finally, it is interesting to observe that $g_0$ may be viewed as 
an indefinite scalar-flat K\"ahler metric on $\CP_1\times \CP_1$, 
and that, conversely, any 
indefinite scalar-flat K\"ahler metric on a complex surfaces is 
automatically self-dual.
In this regard, our techniques lead to the following: 

\begin{main}\label{D}
The only complex surface $(M,J)$ admitting 
Zollfrei scalar-flat indefinite  K\"ahler metrics is $\CP_1\times \CP_1$.
Every such metric arises from a family of analytic disks 
in $\CP_3$ with boundary on a totally real $\RP^3$. 
Near the standard metric $g_0$, moreover, indefinite scalar-flat K\"ahler  metrics  
of fixed total volume 
are in one-to-one correspondence with  those totally real embeddings 
$\RP^3\hookrightarrow \CP^3-Q$ on which the pull-back  of 
the $3$-form 
$$\phi = \Im m \frac{z_1 dz_2 \wedge dz_3\wedge dz_4 - \cdots - z_4 dz_1 \wedge dz_2 \wedge  
dz_3
}{(z_1^2+z_2^2+z_3^2+z_4^2)^2}$$
vanishes. Here $Q$ denotes the quadric surface  $z_1^2+z_2^2+z_3^2+z_4^2=0$. 
\end{main}
In particular, the moduli space of such metrics is once again 
infinite-dimensional.

In the special setting of metrics with  circular symmetry,
 indefinite scalar-flat K\"ahler metrics on $\CP_1\times \CP_1$
were previously investigated  by Tod \cite{tod} and, independently, by
 Kamada \cite{kamada}, both of whom discovered that 
 infinite-dimensional families of such metrics can  be written down in closed 
form by means of the Lorentzian analogue of the first author's hyperbolic ansatz \cite{mcp2}. 
We thus believe that the chief interest of the present
article must be found, not in the mere infinite-dimensionality of the relevant moduli space, 
 but, rather, in the manner in which our    holomorphic disk picture allows one to
explore this interesting, geometric,  
 non-linear ultra-hyperbolic second-order 
 equation in terms of a first-order elliptic boundary-value problem.

\section{Zollfrei Metrics}
\label{zollfrei}

If $(M,g)$ is an indefinite pseudo-Riemannian manifold, a geodesic
$\gamma \subset M$ is said to be a {\em null geodesic} if 
$g(v,v)=0$ for any vector $v$ tangent to $\gamma$. 
We will primarily consider these null geodesics as 
{\em unparameterized} curves, even though $g$
endows them with a preferred class of so-called affine parameters. 
The reason behind this point of view is that 
the null geodesics of a  pseudo-Riemannian  manifold  $(M,g)$
are {\em conformally
invariant} as  unparameterized curves; that is, 
${\zap f} g$ has the same null geodesics as $g$, for any
non-zero function $\zap f$ on $M$. Indeed, let ${\mathscr H}\subset T^*M$
be the  hypersurface
 of non-zero null co-vectors, and notice that ${\mathscr H}$ is foliated by 
 a unique system of curves tangent to $\ker (\omega|_{\mathscr H})$,
 where $\omega$ is the usual symplectic form
on $T^*M$. 
The Hamiltonian formalism then tells us that that the projections into $M$ of these integral 
 curves are precisely the 
 null geodesics of $g$. The  conformal invariance of null geodesics 
 is thus an immediate consequence of the conformal invariance of  ${\mathscr H}$.
 
 Manifolds for which every null geodesic is a simple closed curve will play a central r\^ole in this
 paper, and, 
 following Guillemin \cite{mogul}, we will therefore  introduce some convenient 
 terminology to describe such spaces: 
 
 \begin{defn}\label{simple}
 An indefinite  pseudo-Riemannian manifold $(M,g)$ will be called {\em Zollfrei} if
 the image of  each of its
maximally extended null geodesics is an embedded circle  $S^1\subset M$. 
 \end{defn} 
 
 Notice that this condition is conformally invariant, so it makes perfectly good sense 
 to say that $(M,[g])$ is Zollfrei, where 
$$[g] = \{ {\mathzap f} g ~|~ {\mathzap  f} : M \to \RR^\times\}$$
 denotes the conformal class determined by the metric $g$. 
 
 Guillemin's definition \cite{mogul} is actually a good deal more stringent than Definition
 \ref{simple}. 
Let ${\mathzap Q}\subset \PP T^*M$ denote
the quotient  ${\mathscr H}/\RR^\times$,
where $\RR^\times = \RR - \{ 0\}$ acts by scalar multiplication on $T^*M$.  
The lifts of  null geodesics then define a foliation of ${\mathzap Q}$ by curves. 
Guillemin's definition then amounts to the following: 

\begin{defn}\label{guil}
Let  $(M,g)$ be a  Zollfrei manifold. We will say that $M$ is {\em strongly Zollfrei} 
if the  foliation of ${\mathzap Q}$ by lifted null geodesics is 
a (locally trivial) circle fibration. 
\end{defn}
Since this condition  is obviously also conformally invariant,  the strongly Zollfrei 
condition will  also  be  considered as primarily pertaining to the  conformal class
$[g]$ rather than to the particular metric $g$  representing it. 

If $(M,[g])$ is a strongly Zollfrei $n$-manifold, we may 
defined its {\em space of null geodesics}
$N$ 
to be the leaf-space of the null-geodesic foliation of ${\mathzap Q}$. 
Because the foliation  is assumed to be a locally trivial circle fibration, 
$N$ is then 
automatically a smooth manifold of dimension $2n-3$. The symplectic
description of the foliation endows $N$ with a {\em contact structure},
meaning a maximally non-integrable codimension-$1$ sub-bundle 
$C\subset TN$ of the tangent bundle. Concretely, the tangent space
$T_\gamma N$ of $N$ at a null geodesic $\gamma$ is locally represented
on $M$ as equivalence classes of solutions $w$ of Jacobi's equation
$$\nabla_v\nabla_v w = R_{vw}(v)$$
subject to the constraint
$$g(v,  w) = \mbox{constant}$$
and the equivalence relation 
$$w \sim w + (a + bt) v,$$
where $R$ is the curvature tensor of $g$, 
$t$ is a local affine parameter for $\gamma$,
$v= d/dt$, and $a$ and $b$ are constants; in these terms, the contact 
sub-bundle $C\subset TN$ then corresponds to 
those Jacobi fields $w$ which satisfy the constraint 
$$g(v,w)=0.$$
If ${\mathzap L}\to N$ is the line-subundle of $T^*N$ consisting of
the $1$-forms which annihilate $C$, then  ${\mathzap L}^\times := 
{\mathzap L}-0_N$ is a
symplectic submanifold of $T^*N$,  called \cite{arnold} the 
{\em symplectification} of the contact manifold $(N,C)$. 
However, the pull-back of ${\mathzap L}^\times$ to ${\mathzap Q}$ can 
be canonically identified with ${\mathscr H}$, 
and  the Marsden-Weinstein reduction
of ${\mathscr H}\subset T^*M$ is therefore globally well defined. This shows that
the null geodesic foliation must actually be 
periodic up  on ${\mathscr H}$, and not just down on ${\mathzap Q}$. 
Thus, 
no matter which metric $g$ we choose in the 
conformal class $[g]$, the  the null geodesics of a strongly 
Zollfrei manifold are all automatically {\em periodic} 
with respect to their affine parameters. (This
conclusion should be contrasted with the 
closed but non-periodic null geodesics   \cite{hawkell} of the Taub-NUT metric and
related examples.) For this reason, Definition \ref{guil} is logically equivalent
to the definition  used by  Guillemin  in \cite{mogul}.

\section{Self-Duality}
\label{selfdual}

Suppose that $M$ is an oriented $4$-manifold, and that 
$g$ is a split-signature pseudo-Riemannian metric. Then, as in the Riemannian case, 
the Hodge star operator $\star : \Lambda^2 \to \Lambda^2$
satisfies  $\star^2=+{\mathbf 1}$, so  there is an invariant splitting
$$\Lambda^2 = \Lambda^+\oplus \Lambda^-$$
of the $2$-forms into the $(\pm 1)$-eigenspaces of $\star$. 
The inner product induced by $g$ is of Lorentz signature on both 
$\Lambda^\pm$, reflecting the fact that $SO_+(2,2)$ is a 
double cover of $SO_+(1,2)\times SO_+(1,2)$. 
Sections of $\Lambda^+$ (respectively, $\Lambda^-$) are called
{\em self-dual} (respectively, {\em anti-self-dual}) forms. 
Thinking of the curvature tensor ${\mathcal R}$ of $g$
as a linear map ${\mathcal R}: \Lambda^2 \to \Lambda^2$, 
we thus obtain a decomposition 
$$
{\mathcal R}=
\left(
\mbox{
\begin{tabular}{c|c}
&\\
$W_++\frac{s}{12}$&$\mathring{r}$\\ &\\
\cline{1-2}&\\
$\mathring{r}$ & $W_-+\frac{s}{12}$\\&\\
\end{tabular}
} \right) . 
$$
of the Riemann tensor into simpler pieces. 
Here $W_+$ and $W_-$ are the trace-free pieces of the appropriate blocks,
and  are  called the
{\em self-dual} and {\em anti-self-dual Weyl curvatures}, respectively. 
The scalar curvature  $s$ is understood to act by scalar multiplication,
whereas $\mathring{r}$ is a disguised form of the   
     trace-free part of the Ricci curvature tensor. 

\begin{defn} An oriented split-signature pseudo-Riemannian $4$-manifold
$(M,g)$ is called {\em self-dual} if it satisifes $W_-\equiv 0$. 
\end{defn}

This condition is conformally invariant, in the sense that if $g$ is
self-dual, so is the metric ${\mathzap f} g$, where ${\mathzap f}: M\to \RR^\times$ is any 
non-zero function. 
Thus  the self-duality 
condition should fundamentally be understood    as pertaining to a   conformal class
$$[g] = \{ {\mathzap f} g ~|~ {\mathzap f} : M \to \RR^\times\}$$
rather than to a particular metric $g$  representing it. 

If $(M,g)$ is a pseudo-Riemannian manifold, we will say that a real linear subspace 
$\Pi $ of a tangent space $T_xM$ is {\em isotropic} if it consists entirely of null
vectors. Notice that this is a conformally invariant condition. If 
$(M,g)$ is an oriented split-signature  $4$-manifold,
then the space of isotropic $2$-planes in $TM$ has 
two connected components, each of which is a 
circle bundle over $M$. Indeed, if $\Pi \subset T_xM$ is an
isotropic $2$-plane, then $\wedge^2\Pi $ corresponds, by index-lowering,
to a null $1$-dimensional subspace of either $\Lambda^+$ or
$\Lambda^-$. In the first case, one says that $\Pi $ is an {\em $\alpha$-plane}, 
whereas in the second case one says that $\Pi $ is a {\em $\beta$-plane}. 
We will henceforth use ${\zap  p} : F\to M$ to denote the circle-bundle of
$\beta$-planes over an oriented split-signature $4$-manifold $(M,g)$. 

\begin{defn}
An  immersed  connected surface in $S\looparrowright M$ will be called a
{\em proto-$\beta$-surface} if its tangent space $T_{x}S$ is a 
$\beta$-plane for all $x\in S$.  If, in addition, the proto-$\beta$-surface 
$S$ is maximal, in the
sense that it is not a proper subset of a larger  proto-$\beta$-surface, we will say
that $S$ is a 
{\em $\beta$-surface}.
 \end{defn}

\begin{lem}\label{noii} 
Let $(M,[g])$ be an oriented $4$-manifold with split-signature conformal metric, and 
let 
$S\looparrowright M$ be any proto-$\beta$-surface. Then 
the second fundamental form of $S$ vanishes. Consequently,  
$S$ is totally geodesic. 
\end{lem}
\begin{proof} The tangent bundle of 
any proto-$\beta$-surface $S$ is locally spanned by 
vector fields $v$ and $w$ with $[v,w]=0$ and 
$$g(v,v)=g(w,w)=g(v,w)=0.$$
Now notice $\nabla_vw=\nabla_wv$,
and hence
$$g(v,\nabla_vw)=g(v,\nabla_wv)= \half wg(v,v)=0,$$
whereas we also have
$$g(w, \nabla_vw) = \half vg(w,w)=0.$$
However, $TS= TS^\perp$  with respect to $g$, 
so we conclude that 
$$\nabla_vw\in TS.$$
Similarly, 
$$g(w, \nabla_ww)= \half w g(w,w)=0,$$
and
$$g(v,\nabla_ww) = wg(v,w)- g(\nabla_vw, w) = 0,$$
so we must have  
$$\nabla_ww\in TS,$$
too. Thus
$$\nabla : \Gamma (TS) \times \Gamma (TS) \to \Gamma (TS).$$
In other words,  the second fundamental form
\begin{eqnarray*}
\gemini : TS\times TS & \to  & TM/TS \\
(v,w) & \mapsto & \nabla_vw \bmod TS 
\end{eqnarray*}
 vanishes. Equivalently, 
$S$ is totally geodesic, in the sense that any geodesic
tangent to $S$ at some point necessarily  remains within $S$.  
\end{proof}

This observation then allows one to prove the following:

 \begin{lem}\label{proteus} 
 Let $(M,[g])$ be an oriented $4$-manifold with 
 split-signature conformal structure. Then the
 following are equivalent: 
 \begin{romenumerate}
 \item 
 $[g]$ is 
 self-dual;
  \item  every $\beta$-plane $\Pi \subset TM$
 is tangent to some proto-$\beta$-surface;  
 \item if   $\Pi  \subset TM$ is any $\beta$-plane, and if   
 $v,w\in \Pi $, then ${\mathcal R}_{vw}v\in \Pi $, too. 
 \end{romenumerate}
 \end{lem}
 \begin{proof}
 Suppose that $\Pi \subset TM$ is a $\beta$-plane, and let $v$ and $w$ span $\Pi $. Then,
 using $g$ to freely identify vectors and $1$-forms via index lowering, 
  $v\wedge w\in \Lambda^-$ and $\langle v\wedge w , v\wedge w\rangle =0$.
  Hence 
 \begin{eqnarray*}
g (w, {\mathcal R}_{vw}v) & = & \left\langle v\wedge w , {\mathcal R}(v \wedge w) \right\rangle \\
 & = &  \left\langle v\wedge w , \left( W_-+\frac{s}{12}\right)(v \wedge w) \right\rangle\\
 & = &  \left\langle v\wedge w , W_-(v \wedge w) \right\rangle .
\end{eqnarray*}
On the other hand, 
 $$g (v, {\mathcal R}_{vw}v)=0$$
 by the  Bianchi identities. Since $\Pi =\Pi ^\perp$ with respect to $g$, and because $W_-$
 is a trace-free quadratic form on $\Lambda^-$, it therefore follows
 that  $(iii)$ is equivalent to  requiring that 
  $W_-\equiv 0$. Hence  $(i) \Longleftrightarrow (iii)$.
 
 Now if $S$ is any proto-$\beta$-surface in $M$, 
 and if $v$ and $w$ are any vector fields 
  on $S$, then $\nabla_vw\in TS$ by  Lemma
 \ref{noii}.
 The Riemann curvature tensor of ${\mathcal R}$ of $g$
 thus satisfies  
 $${\mathcal R}_{vw}v = \nabla_v\nabla_wv- \nabla_w\nabla_vv - \nabla_{[v,w]}v\in TS$$
whenever  $S$ is a proto-$\beta$-surface  and  $v,w\in TS$. When every
 $\beta$-plane $\Pi $ can be expressed as $T_xS$ for some proto-$\beta$-surface,
 it thus follows that 
 condition $(iii)$ holds. Hence $(ii) \Longrightarrow (iii)$.  
 
 Conversely, suppose that $(iii)$  holds. Let $\Pi \subset T_xM$ be a $\beta$-plane, and 
 let $\gamma$ be any null geodesic through $x$  tangent to $\Pi $.
 Let ${\Pye }\to \gamma$ be the 
 rank-$2$  sub-bundle of $TM|_\gamma$
 obtained from $\Pi $ by parallel transport of along $\gamma$, 
 and notice that each fiber of $\Pye $ is the unique $\beta$-plane
 containing $v=\gamma^\prime$. Hypothesis $(iii)$  therefore guarantees
 that $$w\in {\Pye } \Longrightarrow 
 {\mathcal R}_{vw}v\in{\Pye },$$ and a solution of  Jacobi's equation
 \begin{equation}
\label{jac}
\nabla_v\nabla_v w =   {\mathcal R}_{vw}v
\end{equation}
 along $\gamma$ is therefore a section of ${\Pye }$ iff  
 $w|_x , (\nabla_vw)|_x\in \Pi $.   Now let  $U\subset \Pi $ 
 be an 
 open disk about $0\in \Pi $
which   is  sufficiently small so as to be mapped diffeomorphically  
to a surface  $S=\exp (U)$ with $T_xS=\Pi $ by the exponential map of $g$. Then $S$ is 
a union of null geodesics $\gamma$ through $x$, and along each such 
 $\gamma$ we  have 
$$T_{\tilde{x}}S
= \left\{ w|_{\tilde{x}}~\Big|~ w \mbox{ solves (\ref{jac}) along } \gamma,  w|_x=0,
(\nabla_v w)|_x\in \Pi \right\}$$  
for each $\tilde{x}\neq x$.
The above argument  thus
 shows that  the tangent spaces of $S$ are precisely the $\beta$-planes
 obtained from $\Pi $ by parallel transport along  radial null geodesics. 
 In particular,  $S$ is a proto-$\beta$-surface tangent to the given $\beta$-plane $\Pi $. Thus 
$(iii)\Longrightarrow (ii)$, and our proof is complete. 
  \end{proof}
  
  Now given an oriented split-signature $4$-manifold $(M,g)$, 
  let us consider the bundle
  ${\zap  p} : F\to M$   of $\beta$-planes. We may define 
a $2$-dimensional real distribution $E\subset TF$ \label{dodger} by declaring that its 
value at $\Pi $ is the horizontal lift of $\Pi \subset TM$ to $T_\Pi F$. 
Then every proto-$\beta$-surface in $M$ has  a canonical lift
as an integral surface of $E$, and conversely every integral surface 
of $E$ projects to a proto-$\beta$-surface in $M$. In this way,
we see that $E$ is integrable iff $(M,g)$ is self-dual. 
In particular, when $(M,g)$ is a self-dual, there is a foliation 
$\mathscr F$ of $F$ tangent to $E$, and we  can then obtain (maximal)  
$\beta$-surfaces in  $(M,g)$ by projecting the
leaves of  $\mathscr F$  into $M$ via ${\zap  p} : F\to M$. 
 Lemma \ref{proteus} thus  implies  

\begin{prop}\label{roger}
The following 
assertions regarding an  oriented split-signature \linebreak $4$-manifold $(M,g)$ are 
logically equivalent: 
\begin{itemize}
\item $g$   is self-dual;
\item   each $\beta$-plane  $\Pi\subset TM$ is tangent to a unique $\beta$-surface  
$S\looparrowright M$;
\item the distribution of $2$-planes $E\to F$ is Frobenius integrable. 
\end{itemize}
\end{prop}

Since Lemma \ref{noii} tells us that each $\beta$-surface $S$ is totally geodesic,  
the Levi-Civita connection $\nabla$ of the ambient  metric $g$ 
induces a torsion-free connection $\triangledown$ on $S$ whose  
geodesics of $\triangledown$ are precisely 
those null geodesics of $(M,g)$ which are contained in $S$.
But, as we saw in  \S \ref{zollfrei}, null geodesics are conformally invariant as 
unparameterized curves, so   the {\em projective class} $[\triangledown ]$ of
this induced  connection \cite{lmzoll,schouten}
therefore depends only  on the ambient conformal class $[g]$.

\begin{prop}\label{projflat}
Let $(M,g)$ be a self-dual split-signature $4$-manifold,  let
$S\looparrowright M$ be any $\beta$-surface, and let 
$\triangledown$ be the  connection induced on  $S$
  by restriction of the 
Levi-Civita connection $\nabla$ to $S$. 
Then $\triangledown$ is {\em projectively flat}.
Indeed, there is a local diffeomorphism  
$\fye : \tilde{S}\to \RP^2$, where $\tilde{S}$ is the universal cover of 
$S$, which  maps each  geodesic  to  a  portion of some  projective line. 
\end{prop}
\begin{proof}
Locally, we have a $3$-parameter family of $\beta$-surfaces in $M$, and, by  
taking the derivatives at $S$, these define 
 a $3$-dimensional  space of sections of the 
normal bundle $TM/TS$ of $S$. These are the covariantly constant local sections for 
 the natural flat connection on the  
normal bundle $TF/TS$ of 
lift of $S$
to $F$ induced by the integrable distribution $E$; and we  obtain a natural $3$-dimensional
space of local sections of $TM/TS$ by pushing forward parallel local sections of
$TF/TS$ via the derivative of $\zap p$.   Moreover, these sections can be taken to be  
{\em global} sections  on the universal cover $\tilde{S}$
of $S$, since  the pull-back of $TF/TS$ to $\tilde{S}$ is not only flat,  but actually
has {trivial} holonomy.

We can describe these local sections more concretely by exploiting our fixed metric
$g$ in the self-dual conformal class $[g]$. Indeed, since $TS$ is maximally isotropic
with respect to $g$, our metric induces a non-degenerate pairing 
$$TS\times (TM/TS)\to \RR,$$
thus giving us an isomorphism between the cotangent bundle $T^*S$ and 
the normal bundle $TM/TS$ of $S$. Thus a section of the normal bundle precisely 
corresponds to a $1$-form $\varphi$ on $S$. We claim that a $1$-form arises 
from a $1$-parameter family of $\beta$-surfaces iff it satisfies the 
{\em generalized Killing equation}
\begin{equation}
\label{gke}
\triangledown \varphi = \half d\varphi~,
\end{equation}
which, according to ones taste, can be rewritten either as  
$$\triangledown_j \varphi_k + \triangledown_k \varphi_j=0$$
or as 
$$
(\triangledown \varphi ) (v,v) = 0~~~\forall v.
$$

Let us first demonstrate the `only if' direction of this assertion. Suppose
we have a $1$-parameter family of proto-$\beta$-surfaces obtained by 
moving some open subset  $U\subset S$. Then we can foliate these surfaces
by null geodesics in a smooth manner, say with tangent vector field $v$.
The vector field $u$ representing the variation then satisfies
$[u,v]=0$, and projects to 
 the section of $TM/TS$ 
along $U$ which represents the first variation of  the family. 
The $1$-form $\varphi$ on $U$ 
 representing the first variation is then given by
$$\varphi (w) = g (u, w) ~~~\forall w\in TS.$$
But now
\begin{eqnarray*}
(\triangledown \varphi ) (v,v)  & = & g (\nabla_v u , v) \\
 & = & g (\nabla_u v , v) \\
 &=& {\textstyle \frac{1}{2}} u g (v,v) \\
 &=& 0.
\end{eqnarray*}
Since $v$ can be chosen to point in any direction at any point of
$U$, it follows that $\triangledown \varphi$ must be skew-symmetric, 
and $\varphi$ is therefore a solution of (\ref{gke}). 

Now (\ref{gke}) is an over-determined equation, and a solution $\varphi$ is 
completely determined by its $1$-jet at a point of $S$. To see this, observe that 
 we certainly have   
 $$\mbox{Alt }( \triangledown\triangledown \varphi ) =0, $$
 since $S$  does not carry any non-zero $3$-forms. 
 Using (\ref{gke}), however, this six-term identity can be rewritten as the three-term identity  
 $$
 \triangledown_j\triangledown_k \varphi_\ell = 
  \triangledown_\ell\triangledown_k \varphi_j
  -\triangledown_k\triangledown_\ell \varphi_j ,
 $$
 and we may then notice  that the right-hand side is just a curvature term. 
 Along a null geodesic $\gamma \subset S$ with tangent field $v$, $\varphi$  therefore
 satisfies the  ordinary  differential equation 
 \begin{equation}
\label{gje}
\triangledown_v \triangledown_v \varphi = \varphi( {\mathcal K}_{v\bullet}v) ,
\end{equation}
 where ${\mathcal K}$ is the curvature tensor of $\triangledown$, and where the 
 right-hand-side denotes the $1$-form
 $$w\mapsto \varphi ( {\mathcal K}_{vw}v) .$$
 Since a solution of (\ref{gje}) is completely determined by 
 the value of $\varphi$ and $\triangledown_v \varphi$ at one point, 
 it follows that solutions of (\ref{gke}) are completely
 determined by the value of $\varphi$ and $\triangledown \varphi$ at
 one point of a convex subset $U\subset S$. But (\ref{gke}) then tells us
  that a solution is consequently determined by the value of the 
 $1$-form $\varphi$ and the $2$-form $d\varphi$ at one point. 
 This shows that the space of solutions is at most $3$-dimensional. 
 But since the codimension of the leaves in $F$ is exactly $3$, 
  we conclude that 
 the space of solutions of (\ref{gke}) must be exactly $3$-dimensional
 up on the universal cover $\tilde{S}$ of $S$. 
  
  Thus, let $\mathbb V\cong \RR^3$ be the  space of solutions of (\ref{gke})
  on $\tilde{S}$, and let $\P ({\mathbb V})\cong \RP^2$ be the corresponding
  real projective space. For each $x\in \tilde{S}$, set 
  $${\mathbb L}_x=  \{ \mbox{solutions $\varphi$ of (\ref{gke}) on $\tilde{S}$ for which } \varphi|_x=0\},$$
  and notice that this is a $1$-dimensional subspace of $\mathbb V$, since 
  the freedom in choosing such a solution amounts to specifying the value
  of the $2$-form $d\varphi$ at $x$. 
  We may thus define a map 
  \begin{eqnarray*}
\fye : \tilde{S} & \longrightarrow & \P ({\mathbb V}) \\
x & \mapsto & {\mathbb L}_x.
\end{eqnarray*}
Let us then first notice that any geodesic $\gamma$ is sent to a projective line
by this map, because equation (\ref{gke}) implies that $\varphi (v) = \mbox{constant}$,
where $v$ is an autoparallel tangent field for $\gamma$; thus $\fye (\gamma ) \subset
\P ({\mathbb V}_0 )$, where ${\mathbb V}_0\subset {\mathbb V}$ is the plane
defined by $(\varphi|_x)(v)=0$ for some arbitrary $x\in \gamma$. 
Now let $t$ be an affine parameter along $\gamma$, with $v=d/dt$,
and let $w$ be a parallel vector field along $\gamma$ which is linearly
independent from $v$. Then the restriction of an element of ${\mathbb V}_0$
to $\gamma$ satisfies $\varphi (v) \equiv 0$ and $\varphi (w) = f(t)$,
where $f$ is a solution of the second order linear ordinary differential  equation
\begin{equation}
\label{ode}
\frac{d^2f}{dt^2}+ \kappa f=0,
\end{equation}
where $\kappa (t)= {{\mathcal K}^2}_{121}$ with respect to the frame $e_1=v$, $e_2=w$. 
If $\{f_1, f_2\}$ is a basis for the solution space of (\ref{ode}), then $\fye$ sends 
$\gamma$ to $\P ({\mathbb V}_0)\cong \RP^1$ by 
$t\mapsto [f_2 (t) : - f_1(t)]$. However, equation  (\ref{ode}) implies that the
Wronskian $W= f_1f_2^\prime - f_1^\prime f_2$ is a non-zero constant along $\gamma$.
Thus at least one of the expressions  
$$\frac{d}{dt}\left( \frac{-f_1}{f_2}\right) = \frac{W}{f_2^2}~~~\mbox{ and }~~~
 \frac{d}{dt}\left( \frac{f_2}{-f_1}\right) = -\frac{W}{f_1^2}$$
is defined and non-zero at each point of $\gamma$, and $\fye$ thus sends $\gamma$ to the projective
line $\P ({\mathbb V}_0) \subset \P ({\mathbb V})$ via a smooth immersion.
Since the geodesic $\gamma$ is arbitrary, it follows that $\fye : \tilde{S}\to \P ({\mathbb V})$
is an equidimensional  smooth immersion  sending each geodesic to 
a portion of a projective line. In particular, the connection $\triangledown$ 
induced on the $\beta$-surface $S$ 
is projectively flat. 
\end{proof}

The above proof  is loosely based on a  spinor  argument   given  in \cite{thick}.
The careful reader may notice that, in its present form, the proof is not manifestly 
conformally invariant. However, it is not terribly difficult to embellish the argument
so as to achieve this end. The main point is that the $g$-induced map $TM/TS\to T^*S$
actually carries a conformal weight, so that the $1$-form fields $\varphi$
under discussion may better be described as $1$-forms with values in a
line bundle.

It is also worth pointing out  that the above result   depends quite strongly  on the assumption
that $M$ is self-dual.
Indeed, it is not difficult to construct  non-self-dual $4$-manifolds
with isolated $\beta$-surfaces on which  the induced connection is {\em not} projectively 
flat. For example, let $(\Sigma , h)$ be an oriented  Riemannian $2$-manifold of
{\em non-constant} Gauss curvature. 
Since \cite{schouten} a torsion-free connection $\triangledown$ on a surface
is projectively flat   iff its Ricci curvature $\rho$ satisfies 
$$
2\triangledown_j \rho_{k\ell }  - 2 \triangledown_k \rho_{j\ell} + 
\triangledown_j \rho_{\ell k}  - \triangledown_l \rho_{\ell j} =0,
$$
it follows that the Riemannian connection of such a generic metric $h$
is  not projectively flat. 
Now give $\Sigma\times \Sigma$ the 
 indefinite product metric $\pi_1^*h-\pi_2^*h$, and observe that 
 the diagonal 
$\Sigma \hookrightarrow \Sigma \times 
{\Sigma}$ becomes a $\beta$-surface if we endow $\Sigma\times \Sigma$
with the non-product orientation. However, the induced connection
on this $\beta$-surface is just the Riemannian  connection of $h$,
so this $\beta$-surface is {\em not} projectively flat. Of course, this
example 
in no way contradicts Proposition \ref{projflat}, since the $4$-manifold
in question is non-self-dual.

\section{Projectively Flat Surfaces}

 Proposition \ref{projflat} reveals that  surfaces with 
 flat projective connections  play an important r\^ole in the 
 theory of split-signature self-dual manifolds. The systematic study
 of surfaces
 with flat projective structures, also known as $\RP^2$-structures, 
 was begun by Kuiper \cite{nico}, who in particular observed
 that if $(S,[\triangledown ])$ is any projectively flat surface, 
 the locally trivial nature of the geometry always
 gives rise to  a developing map $\fye : \tilde{S}\to \RP^2$, defined 
on    the universal cover  $\tilde{S}$ of $S$, as well as a representation
of $\phi : \pi_1 (S) \to PGL (3, \RR )$, both of which are unique up to 
an overall $PGL(3, \RR)$ transformation.
Crucial explorations of this idea by 
Sullivan and Thurston \cite{sulthurs} eventually allowed 
 Choi and Goldman \cite{chogo} to develop a 
 substantially complete theory of flat projective structures on {\em compact} surfaces.

 In this article, we will be 
specifically 
 interested in the case when the relevant projective structure is 
 {\em Zoll}, meaning \cite{lmzoll} that every geodesic is a simple closed curve. 
It seems quite plausible that 
 a connected surface which admits a 
Zoll projective connection must necessarily be compact, 
but, to our knowledge, this  still seems to be an open problem.
Fortunately, however, the  
{\em projectively  flat}  case of the problem is a bit more  manageable.

\begin{lem}\label{closed} 
Let $(S, [\triangledown])$ be a connected surface with
flat projective structure, and suppose that every maximal geodesic
of $[\triangledown ]$ is a simple closed curve in $S$. 
Then $S$ is compact. 
\end{lem}
\begin{proof}
It suffices to consider the case when $S$ is 
 {\em orientable}, since otherwise we may pass to an oriented
 double cover without sacrificing the assumption that every
 geodesic is a simple closed curve.

 Since $S$ is assumed to be Zoll, 
every (maximal) geodesic $\gamma\subset S$ is an embedded circle; and  
because we may now  assume that $S$ is orientable, any 
such $\gamma$ has an open  neighborhood
$U\subset S$ 
diffeomorphic to an annulus $S^1 \times (-\epsilon , \epsilon )$. 
Now develop the universal cover $\tilde{U}$ of $U$ onto the $2$-sphere
in such a manner that $\gamma$ is sent to some portion of the equator $x_3=0$,
sending a chosen base-point to $(1,0,0)$.We orient the equator in the usual
west-to-east manner, and give $\gamma$ the corresponding orientation. 
Let $I\subset \tilde{U}$ be an arc (that is, an embedded closed interval)
such that $\gamma \subset U\subset S$ is obtained from $I$
by identifying its two endpoints via the covering map $\tilde{U}\to U$,
and such that the initial end-point is a pre-image of the chosen base-point for $S$.
Then the restriction of the developing map 
to some open  neighborhood of $I\subset \tilde{U}$
 lifts to the universal cover $V$ of $S^2-\{ (\pm 1, 0, 0)\}$, 
 where $V$ may  may be explicitly identified with  
$\RR\times (-\pi/2, \pi/2)$
through the use of spherical coordinates
$$ (x_1, x_2, x_3) = (\cos \theta \cos \varphi , \sin \theta \cos \varphi , \sin \varphi ),
~~~~ (\theta , \varphi )\in \RR\times (-\pi/2, \pi/2).$$
 This lift of the  development then takes $I$ diffeomorphically
onto a closed interval, say,  $\tilde{I}= [0,L]\times \{0\}$ in $V$. Thus, a perhaps 
smaller  neighborhood $U^\prime$ of $\gamma\in S$ can be
obtained from  a neighborhood $V^\prime$ 
of $\tilde{I}\subset V$ by identifying some neighborhood $V_1$
of $(0,0)$ with a neighborhood $V_2$ of $(L,0)$; moreover, 
this identification is carried out via a lift of the action of some $A\in SL(3, \RR)$ 
of $S^2= (\RR^3-0)/\RR^+$. Notice that this transformation $A$ must  send
 the equator to itself,  in an orientation-preserving manner. Hence  
 $(0,0,1)$
must be an eigenvector of $A^t$, with  eigenvalue $\lambda > 0$. 

Let us now examine the action of $A^t$ on the space $\RP^{2*} = \P (\RR^{3*})$
of great circles in $S^2$. We have just observed that $[0,0,1]$ must be a fixed point
of this action. But our hypotheses also preclude the existence of a point $p\in \RP^{2*}$,
$p\neq [0,0,1]$,  such that 
$\lim_{n\to \infty} (A^t)^n(p)= [0,0,1]$. Indeed, if there were such $p$, the great circles
corresponding to the iterates $(A^t)^n(p)$ would, for $n\gg 0$, link up end-to-end
via $A$ to form part of a geodesic $\gamma^\prime\neq \gamma$ in the annulus $U \subset S$
which spiraled into $\gamma$; 
 every point of   $\gamma$ would then be an accumulation point of $\gamma^\prime$, and 
as the Zoll hypothesis implies that $\gamma^\prime$ must be a closed
subset of $S$,  this would imply that $\gamma\subset \gamma^\prime$, 
contradicting the fact that $\gamma$ is a maximal 
geodesic. We can also run this argument backwards in parameter time by
replacing $A$ with $A^{-1}$, and thereby  deduce that 
  there cannot be any point $p\in \RP^{2*}$, $p\neq [0,0,1]$,
   such that $\lim_{n\to \infty} (A^t)^{-n}(p)= [0,0,1]$.
The complex eigenvalues of 
$A$ must therefore all have the same modulus. Moreover, there cannot be a vector 
$v\in \RR^{3*}$ such that $A^t(v) = \lambda v + (0,0,1)$. 
Hence $A\in SL(3,\RR)$ can be put in  one of the normal 
forms
$$(a) ~~
\left[\begin{array}{ccc}\cos \psi & -\sin \psi & 0 \\ \sin \psi & \hphantom{-}\cos \psi & 0 \\0 & 0 & 1\end{array}\right]
~~~\mbox{or}~~~(b)~~
\left[\begin{array}{ccc}\pm 1 & \hphantom{\pm}1 & \hphantom{\pm}0 \\\hphantom{\pm}0 & \pm1 & \hphantom{\pm}0 \\\hphantom{\pm}0 & \hphantom{\pm}0 & \hphantom{\pm}1\end{array}\right]
$$
by an appropriate change of basis of the $x_1x_2$-plane. 

  Now suppose that $A$ has normal form $(a)$. Then any geodesic with
  initial point and tangent direction close to that of $\gamma$ will remain
  in our annular neighborhood $U$; indeed,  every such  geodesic is explicitly represented
  in our $(\theta, \varphi)$ coordinates as the union of the graphs 
  $$  \varphi  = \tan^{-1} (t\sin (\theta -\theta_0 + k\psi )), ~~\theta \in [0,\psi ], ~~k\in \ZZ,$$
  for $t$ and $\theta_0$ given constants, with $t$ is sufficiently small,
  and where the ostensible $\bmod$-$2\pi$ ambiguity of $\psi$ has been remedied by 
  setting $\psi=L$.  But the Zoll condition stipulates that
  every geodesic is a simple closed curve, and a simple closed curve in an
  annulus necessarily has winding number   one. Thus the Zoll assumption
  guarantees that $\psi$ is a multiple of $2\pi$, and the developing map 
  will therefore  be well defined on a neighborhood of any such geodesic
  $\gamma\subset S$, even though general principles had led us to expect that
  it would merely be defined up on the universal cover $\tilde{S}$.   
  Moreover, the subset of $\P(TS)$ consisting of directions tangent to 
  geodesics $\gamma$ with this normal form is {\em open}. 
  
  On the other hand, if $A$ has normal form $(b)$, then 
  the $\pm$ sign must be $+$; if not, an annular neighborhood of the given geodesic
  $\gamma$ would contain closed geodesics with self-crossings, obtained by 
  gluing together great circles near the equator with their reflections through the
  $x_3$-axis. Thus, the transformation $A$ must take the normal form 
  $$(a) ~~
\left[\begin{array}{ccc}1 & 0 & 0 \\ 0& 1 & 0 \\0 & 0 & 1\end{array}\right]
~~~\mbox{or}~~~(b)~~
\left[\begin{array}{ccc} 1 & 1 & 0 \\0 & 1 & 0 \\ 0 & 0 & 1\end{array}\right]~,
$$
and we will henceforth say that a given geodesic $\gamma$ is type $(a)$ or type
$(b)$ depending on which one of these normal forms occurs.

   Now recall that    
the developing-map construction  gives us an immersion 
$\tilde{\fye} : \tilde{S}\to S^2$ and a group homomorphism 
$\phi : \pi_1 (S) \to SL(3,\RR )$ such that the deck-transformation action of $\pi_1 (S)$
on the universal cover $\tilde{S}$ is compatible with the 
action of $SL(3,\RR )$ on $S^2$. In particular, we may define  a
natural intermediate cover $\hat{S}$ of $S$ by setting $\hat{S}= \tilde{S}/\ker \phi$,
so that $S$ is then obtained from $\hat{S}$ by dividing by the action of a
matrix group $G= \phi [\pi_1 (S)]\subset SL(3,\RR )$, and such that
we still have a developing map $\hat{\fye}: \hat{S}\to S^2$ which 
correctly intertwines the effective actions of $G$ on $\tilde{S}$ and $S^2$. 
Let $\varpi : \hat{S}\to S$ be the covering map. 
If $\gamma\subset S$ is a geodesic of type $(a)$, then
 $\varpi^{-1}(\gamma)=\coprod_j \hat{\gamma}_j$, where each $\gamma_j\subset
 \hat{S}$ is a  closed geodesic  of type $(a)$, and where $\varpi|_{\hat{\gamma}_j}
 : \hat{\gamma}_j\to \gamma$ is a diffeomorphism for each $j$;
 this follows immediately from the fact that every non-trivial deck transformation 
 of $\hat{S}$ must act non-trivially on $S^2$, whereas  non-trivial
 coverings of a closed geodesic of type $(a)$ are invisible to the 
developing map. 
%Notice that every element of $G-\{ 1\}$ 
% must therefore act on $\varpi^{-1}(\gamma)$ by sending each copy
% $\hat{\gamma}_j$ of $\gamma$ diffeomorphically to a different one. 
 
On the other hand,
if $\gamma^\prime\subset S$ is a geodesic of type $(b)$, then  
$\varpi^{-1}(\gamma^\prime)=\coprod_j \hat{\gamma}_j^\prime$,
where each $\hat{\gamma}_j^\prime\approx \RR$ is a non-closed geodesic
in $\hat{S}$; moreover, the conjugates of the image of $[\gamma]\in \pi_1(S)$ in $G$ 
give us  deck transformations of $\hat{S}$ which 
 roll
up the various  $\hat{\gamma}_j^\prime$ into  copies of $\gamma$, while
simultaneously acting on $S^2$ via  linear transformations of normal form $(b)$. 
Such a matrix acts on $S^2$ in a manner fixing a great circle, and on this 
great circle there is a preferred antipodal pair of points, given by  
$(\pm 1 , 0, 0)$ for the standard model, which are the accumulation points
of the non-closed orbits, 
 and which we will refer to as the {\em goals} of $\hat{\gamma}_j^\prime$.
Since $S$ is paracompact, $G= \pi_1(S)/\ker \phi$ is countable, so
it follows that only countably many points of $S^2$ occur as goals of geodesics of 
type $(b)$. 

Now let $\hat{x}\in \hat{S}$ be any point that is not sent to a goal, and 
let $x=\varpi (\hat{x})$ be its projection to $S$. Then only countably many
geodesics through $x$ can be of type $(b)$, since any such geodesic would 
develop onto a great circle joining $\hat{\fye} (x)$ to a goal. Moreover, 
the set of directions in $\PP (T_xS)\approx S^1$ which are tangent to
geodesics of type $(a)$ is open. Thus the set  $B\subset \PP (T_xS)$
of directions tangent to geodesics of type $(b)$ is a countable closed
subset of the circle. We claim that $B=\emptyset$. If not, 
choose a base-point for $\PP (T_xS)$ which is not
in $B$ and use the counter-clock-wise angle from this direction to define 
 a homeomorphism $\PP (T_xS)\approx [0, \pi]/\{ 0,\pi\}$
which sends the base-point to the equivalence class $\{ 0, \pi\}$. 
Then $B$ then becomes a non-empty countable closed subset 
$\Bye \subset (0,\pi )$.  Let $\bye = (\sup\Bye) \in\Bye$, and let 
$b\in \PP (T_xS)$ be the corresponding direction. Let 
${\mathscr I}\subset \PP (T_xS)$ be the open subset corresponding to 
the open interval
$(\bye ,\pi )$. Every direction in $\mathscr I$ is tangent to a 
geodesic of type $(a)$, and since every such 
geodesic $\gamma$ has an annular neighborhood
which   looks like a finite  covering of a band around the equator,
all the geodesics $\gamma_t$ tangent to elements of $\mathscr I$ form a 
smooth family of maps of the circle, and in particular are all 
homotopic to one another.  We can then uniquely lift this family of 
geodesics of type $(a)$ as a smooth family $\hat{\gamma}_t$, $t\in (\bye ,\pi )$, 
of closed geodesics though $\hat{x}\in \hat{S}$. Let $Y\subset \hat{S}$ 
be the union of these curves $\hat{\gamma}_t$, and notice that
$\varpi|_Y$ is  injective, since all the $\gamma_t$  are homotopic. 
If  $a\in G-\{ 1\}$, it thus follows that 
 $a(Y)\cap Y = \emptyset$. Now  let $\gamma^\prime\subset S$ 
 be  the geodesic of type $(b)$ through $x$ with tangent $b$,  and
 let $\hat{\gamma}^\prime$ be its lift through $\hat{x}$.
 By composing $\hat{\fye}$ with an element of $SL(3,\RR)$ 
 if necessary, we can henceforth assume that $\hat{\fye}[\hat{\gamma}^\prime]$
 is a  subset of the equator $z=0$, that  
 $\hat{\fye} (\hat{x})= (0,1,0)$, and that there is an element $a\in G$ 
 which sends $\hat{\gamma}^\prime$ to itself, while acting on 
 $S^2$ by $(x,y,z)\mapsto (x+y,y, z)/\|(x+y,y,z)\|$. 
 Let $\sigma\subset  \hat{\gamma}^\prime$ consist of 
 those points of $\hat{\gamma}^\prime$ for which every neighborhood
 meets every $\hat{\gamma}_t$ for $t\in (\bye ,\bye + \epsilon )$, where 
 $\epsilon> 0$ is allowed to depend on the neighborhood.
Now 
  the developing map $\hat{\fye}$ is a local diffeomorphism, and carries $Y$ 
  onto $\{ (0 , \pm 1 , 0)\} \cup \{ y < z \cot \bye, z > 0\} \cup \{ y >  z \cot \bye,  z <  0\}$
by a finite covering map. It follows that 
the non-empty subset $\sigma \subset \hat{\gamma}^\prime$ is therefore
both open and closed. Hence 
 $\sigma = \hat{\gamma}^\prime$. In particular, 
   any  point of $\hat{\gamma}^\prime$ which is sent to
the semi-circle $\{ z=0, y > 0\}$ is contained in an open  disk which intersects
$Y$ in an open half-disk consisting of points south  of $\hat{\gamma}^\prime$. 
Hence each of the iterates $a^n (x)$, $n > 0$, 
has an open neighborhood $U_n$ such that $a(U_n\cap Y) \cap (U_{n+1}\cap Y)
\neq \emptyset$. But this means that $a(Y)\cap Y \neq \emptyset$, which is
a contradiction. Hence $\Bye = \emptyset$, and every geodesics through 
$x$ is of type $(a)$. 

The set of all geodesics though $x$ therefore forms a smooth family
of embedded circles.
If $\tilde{X}\subset \PP (TS)$ denotes the union of all the lifts of geodesics 
through $x$, then $\tilde{X}$ is a smooth compact surface---in fact, a Klein bottle. 
Moreover, $\PP (T_xS)$ is a subset of $\tilde{X}$, and this circle has non-orientable
normal bundle. Let $X$ be the smooth compact 
surface---actually, a projective plane---obtained from $\tilde{X}$
by blowing this circle down to a point $x_0\in X$.  The projection $\tilde{X}\to S$
then induces a smooth proper map $f: X\to S$ such that $f_* : T_{x_0}X\to T_xS$
is an isomorphism. But, by assumption,  any geodesic passes though 
$x$ only once. Thus $f^{-1}(x)= \{ x_0\}$, and the mod-$2$  degree of 
$f$ is therefore $1\in \ZZ_2$. If $f$ were not onto, this would be a contradiction, 
since any regular value must have an odd number of points in its pre-image. 
Hence $f$ is onto, and  $S=f(X)$ is compact, as claimed. 
\end{proof}

We therefore obtain  the following useful result: 
 
\begin{thm} \label{standard} 
Let $(S, [\triangledown])$ be a connected surface with
flat projective structure, and suppose that every maximal geodesic
of $[\triangledown ]$ is a simple closed curve in $S$. 
Then, up to diffeomorphism,  $(S,[\triangledown])$ is  either $S^2$ or $\RP^2$,
equipped with the standard  projective connection. 
\end{thm}
\begin{proof}
By Lemma \ref{closed}, $S$ is a compact Zoll manifold. Hence
 \cite[Lemma 2.8]{lmzoll} tells us that $S$ is diffeomorphic to 
 either $S^2$ or $\RP^2$. In particular, the universal cover $\tilde{S}$ 
 of $S$ is compact, so the developing map $\tilde{\fye}: \tilde{S}\to 
 S^2$ must be a covering map. Hence $\tilde{\fye}$ a diffeomorphism. 
 
 If $S\approx S^2$, $\tilde{\fye}$ is now a diffeomorphism $S\to S^2$ which 
 sends the given flat projective structure to the standard one, and 
 we are done. 
  
 If $S\approx \RP^2$, the non-trivial deck
 transformaton of $\tilde{S}\approx S^2$ defines a linear involution
 for which $+1$ is not an eigenvalue. But the only such involution is
 $-1$. Thus we actually obtain a developing map induces $S\to \RP^2$, and this
  gives us  the promised diffeomorphism 
 sending  the given  projective structure to the standard one.  \end{proof}
 
 In the next section, we will see that this has some interesting ramifications for the theory of 
 Zollfrei self-dual $4$-manifolds.

\section{Topological Implications}
\label{prelim}

In this section,  we will
show that,  
 up to homeomorphism,  the only oriented  $4$-manifolds
which admit  self-dual Zollfrei metrics are  $S^2 \times S^2$ and
the real projective quadric 
${\mathbb M}^{2,2} = [S^2 \times S^2]/\ZZ_2$. 
We begin our proof with the following observation:

\begin{lem} \label{beta}
Let $(M,[g])$ be a Zollfrei self-dual $4$-manifold. Then every 
$\beta$-surface $S\looparrowright M$ is  an embedded $S^2$ or $\RP^2$
in $M$. Moreover, every two points of such a
$\beta$-surface $S$ are joined by a null geodesic $\gamma$. 
\end{lem}
\begin{proof}
Let $(M,[g])$ be a Zollfrei self-dual $4$-manifold, and let $S\looparrowright M$
be a  $\beta$-surface. By Lemma \ref{noii}, $S$ is totally geodesic, so 
the immersion $S\looparrowright M$ is injective outside a discrete subset,
where the various tangent spaces of $S$ must be transverse to each other. 
Moreover, 
 $\nabla$  induces a connection $\triangledown$ on $S$.
 Proposition \ref{projflat} asserts that the associated projective 
structure $[\triangledown ]$ is {\em flat}. But the geodesics of $(S, [\triangledown ])$
are all null geodesics of $[g]$, so the assumption that $(M,[g])$ is 
Zollfrei implies that $(S, [\triangledown ])$ is a projectively flat 
surface in which every geodesic is a simple closed curve. Theorem \ref{standard} therefore 
 tells us that 
$S$ is diffeomorphic to either $S^2$ or $\RP^2$, in such a manner that 
$[\triangledown ]$ becomes the standard projective structure. 
In particular, every pair of points of $S$ can be joined by a geodesic
of $[\triangledown ]$. Since the restriction of $S\looparrowright M$
to any geodesic yields an immersion which is one-to-one outside a finite
number of double-points with distinct tangents, the assumption  that every null
geodesic of $[g]$ is a simple closed curve therefore implies that 
$S\looparrowright M$ is actually an embedding. 
\end{proof}

When we say that $(M,[g])$ is self-dual, it is already implicit that
$M$ is oriented. However, $O(2,2)$ has  {\em four} components;
indeed, the inclusion $O(2) \times O(2) \hookrightarrow O(2,2)$
is a homotopy equivalence.
We will say that an oriented  
 split-signature pseudo-Riemannian $4$-manifold  is {\em space-time-orientable} if its
structure group can be reduced to the identity component
$SO_+(2,2)$  of $O(2,2)$. Obviously this is automatically the case if
$H^1(M,\ZZ_2)=0$, and so in particular holds whenever $M$
 is simply connected. If $M$ is {\em not} space-time-orientable, there is always a 
 double cover $\tilde{M}\to M$ which {\em is} space-time orientable 
 with respect to the pull-back of the metric. Moreover, 
 $\tilde{M}$ will then be  Zollfrei if $M$ is, since all the null geodesics of 
 $\tilde{M}$ are at worst double covers of those in $M$.

Now suppose that $(M,g)$ is a space-time-orientable split-signature self-dual 
$4$-manifold. Then we may express $TM$ as a direct sum $T_+\oplus T_-$,
where $T_+$ and $T_-$ are mutually orthogonal with respect to $g$, and where
the restriction of $g$ to $T_+$ (respectively, to $T_-$) is positive (respectively, 
negative) definite; for example, this may be done by  choosing some background
Riemannian metric $h$ on $M$, and then diagonalizing $g$ with respect to $h$
at each point. A space-time orientation for $M$ then amounts to a choice
of orientations for the bundles $T_\pm$. Notice that  this then allows us to 
express $TM$ as the sum of two complex line bundles; indeed, a reduction
of the structure group of $(M,g)$ to $SO(2) \times SO(2) = U(1)\times U(1)$
is equivalent \cite{matsushita} to the choice $( {\mathfrak J} , \tilde{{\mathfrak J}})$ 
of a pair of $g$-compatible almost-complex structures, 
where ${\mathfrak J}$ is compatible with the given orientation 
of $M$, and where $\tilde{{\mathfrak J}}$ is compatible with the opposite orientation. 
An isotropic $2$-plane $\Pi \subset T_xM$ then becomes the graph of an isometry
from $(T_{-x}, -g)$ to $(T_{+x},g)$, and such an isotropic  subspace $\Pi $ is then 
a $\beta$-plane iff this isometry is {\em orientation-reversing}. 
In particular, the orientation of $T_-$ induces an orientation on every $\beta$-plane,
and hence on any $\beta$-surface;  what is more,  any $\beta$-surface is a 
 $\tilde{{\mathfrak J}}$-holomorphic curve in $M$. We thus obtain the following:

\begin{lem}\label{sphere}
Let $(M,[g])$ be a  space-time-orientable Zollfrei self-dual $4$-manifold. 
Then every $\beta$-surface $S$ in $M$ is an embedded $2$-sphere. 
\end{lem}
\begin{proof}
A space-time orientation induces an orientation of each 
$\beta$-surface. Lemma \ref{beta} therefore tells us that 
each $\beta$-surface must be an embedded $2$-sphere. 
\end{proof}

The following observation is therefore pertinent:   

\begin{lem}\label{unisphere} 
Suppose that  $(M,[g])$ is a  split-signature  self-dual  
$4$-manifold in which every $\beta$-surface is an
immersed $2$-sphere. Let ${\zap p}: F\to M$
be the bundle of $\beta$-planes over $M$. 
Then the canonical foliation $\mathscr F$  of $F$ by lifted $\beta$-surfaces
is locally trivial, in the sense that every leaf has
a neighborhood which is diffeomorphic to $S^2 \times \RR^3$ 
in such a manner that each first-factor sphere $S^2 \times \{ *\}$ is 
a leaf. Moreover, this diffeomorphism can be chosen in such a way that each
great circle in each first-factor sphere projects to a null geodesic in $M$. 
\end{lem}
\begin{proof}
Since every leaf of $\mathscr F$ is compact and simply connected, 
the holonomy of $\mathscr F$  around any leaf is trivial, and 
$\mathscr F$ is therefore a fibration \cite{thurstfol}. In particular, we can choose 
a transversal $U$ through a given leaf which meets every nearby leaf exactly once. 
Assume, without loss of generality, that 
$U\approx \RR^3$, and let  $V\approx U \times S^2$ be the corresponding neighborhood of the leaf. 
Since $V$ is simply connected, the line bundle $\ker {\zap p}$ becomes
trivial when restricted to $V$, and we can therefore choose a non-zero vector  
field $u$ on $V$ which spans $\ker {\zap p}$. The $U$-component of this
vector field  then defines a function $V\to (\RR^3 - \{ 0\})$, and we thus
get a map  $V\to S^2$ by composing with the radial projection $(\RR^3 - \{ 0\})\to S^2$.
Modulo the action of the 
$SL(3, \RR)$, however, the restriction of this map to any 
 leaf $S$ is really just the   developing map  $\fye : S\to \RP^2$
 constructed in Proposition \ref{projflat},  lifted 
to the universal cover $S^2$ of $\RP^2$. Taking the Cartesian product  with
leaf projection $V\to U\approx \RR^3$, we thus obtain a local trivialization 
$V\to S^2\times \RR^3$  of $\mathscr F$ in which every lifted null 
geodesic becomes a great circle in a first-factor $S^2$. 
\end{proof}

This gives us a more transparent understanding  of the 
Zollfrei condition: 

\begin{thm}\label{tasty}
Let $(M,[g])$ be a space-time-orientable self-dual  
$4$-manifold. Then the following conditions are equivalent:
\begin{romenumerate}
\item $(M,[g])$ is Zollfrei;  
\item $(M,[g])$ is strongly Zollfrei;  
\item  every $\beta$-surface 
is an embedded $2$-sphere  in $M$. 
\end{romenumerate} 
\end{thm}
\begin{proof}
Definition \ref{guil} tells us that 
 $(ii) \Longrightarrow (i)$, while Lemma \ref{sphere}
 asserts  that $(i) \Longrightarrow (iii)$. It  therefore suffices 
 to show that $(iii) \Longrightarrow (ii)$.

Thus, suppose that every $\beta$-surface of  $(M,[g])$ is an
embedded $2$-sphere in $M$. 
Let ${\zap Q}$ be the bundle of null directions of $(M,[g])$, and notice that 
the  bundle projection ${\zap Q}\to M$   factors through  
an $S^1$-fibration  ${\zap Q}\to F$, since every non-zero null vector is an
element of exactly one $\beta$-plane. But Lemma \ref{unisphere}
tells us that the foliation of ${\zap Q}$
by lifted null geodesics  simplifies when restricted to the null
geodesics in a $\beta$-surface, where it just becomes  the
standard fibration $\P (TS^2) \to \RP^2$;   moreover, 
this picture applies uniformly
in a neighborhood of each leaf of the foliation $\mathscr F$
of $F$. Hence the foliation of $\zap Q$ by lifted null geodesics
is a locally trivial circle fibration. Since each null
geodesic lifts to a great circle in a leaf of $\mathscr F$, and
each leaf embeds into $M$ via ${\zap p}: F\to M$, each
null geodesic is also an embedded circle. It therefore follows that  $(M,[g])$
is strongly Zollfrei, and we are done. 
\end{proof}

We also obtain the following crucial fact: 

\begin{lem}
Suppose that  $(M,[g])$ is a  space-time-orientable self-dual  Zollfrei 
$4$-manifold, and let ${\zap p}: F\to M$
be the bundle of $\beta$-planes over $M$. 
Then there is a smooth  $3$-manifold $P$ and a smooth proper 
submersion ${\zap q}: F\to P$ whose fibers are exactly the 
leaves of the foliation $\mathscr F$. 
\end{lem}
\begin{proof}
By  Lemma \ref{sphere}  and Lemma 
\ref{unisphere},    $\mathscr F$ must be  a locally
trivial fibration by $2$-spheres, and the leaf space $P$ is therefore a manifold. 
\end{proof}

The situation is thus encapsulated by  the diagram 
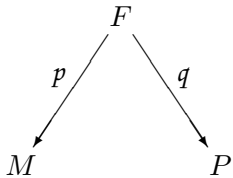
\setlength{\unitlength}{1ex}
\begin{center}\begin{picture}(20,17)(0,3)
\put(10,17){\makebox(0,0){$F$}}
\put(2,5){\makebox(0,0){$M$}}\put(18,5){\makebox(0,0){$P$}}
\put(15,12){\makebox(0,0){${\zap  q}$}}
\put(5,12){\makebox(0,0){${\zap  p}$}}
\put(11,15.5){\vector(2,-3){6}}
\put(9,15.5){\vector(-2,-3){6}}
\end{picture}\end{center}
which we shall refer to as the (real) {\em  double fibration} 
of $(M,[g])$.

Now since $F$ is connected, so is $P={\zap q}(F)$, and we may therefore 
join any two distinct points of $P$ by  a smoothly embedded arc. Trivializing the restriction of 
${\zap q}$ to this arc then results in a free homotopy of the corresponding
leaves of $\mathscr F$. Finally,  pushing this homotopy down via ${\zap p}$
 produces a free homotopy of any two given $\beta$-surfaces
in $M$. In particular, any two $\beta$-surfaces are homologous:

\begin{lem}\label{homotopic}
Let $(M,[g])$ be a space-time oriented Zollfrei self-dual $4$-manifold.
Then any two   $\beta$-surface $S, S^\prime \subset M$  
are freely homotopic, and so, in particular, represent the same 
homology class in  $H_2 (M, \ZZ)$. \end{lem} 

Now since $M^4$ is oriented, there is a well defined intersection 
form $$H_2(M,\ZZ)\times H_2(M,\ZZ)\to \ZZ $$ even if $M$ is non-compact;
for example, 
this reflects the fact that we always have  a Poincar\'e-duality isomorphism 
$H_2(M)\cong H^2_c(M)$ as well as 
a natural homomorphism $H^2_c(M)\to H^2(M)$. If $S_1$ and $S_2$ are compact 
embedded oriented surfaces in general position, one assigns a local intersection
index of $\pm 1$ to each intersection point $x\in S_1\cap S_2$
so as to indicate whether the given orientations of 
 $T_xS_1\oplus T_xS_2$ and $T_xM$ agree or disagree,
 and the  homological intersection
number $[S_1]\cdot [S_2]$ is then precisely the sum of these intersection 
indices. When  $S_1$ and $S_2$ happen to be $\beta$-surfaces, 
we thus obtain the following: 

\begin{lem}\label{insect}
If $S_1$ and $S_2$ are  distinct compact embedded $\beta$-surfaces
in  a  space-time-orientable  self-dual $4$-manifold  $(M,[g])$, then 
 their  homological 
intersection number $[S_1]\cdot [S_2]$ equals $-\# (S_1\cap S_2 )$. 
\end{lem} 
\begin{proof}
Each $\beta$-surface is totally geodesic, so two distinct 
$\beta$-surfaces can never share the same tangent space. 
Since distinct $\beta$-planes in any tangent space are transverse,
this shows that $S_1$ and $S_2$ are necessarily in general position. 
Now since any $\beta$-plane may be viewed as the graph of 
an orientation-reversing isometry $T_-\to T_+$, 
 the intersection index assigned to each point of intersection is 
 $-1$. Summing over intersection points thus yields 
 $[S_1]\cdot [S_2]=-\# (S_1\cap S_2 )$. 
\end{proof}

When $(M,[g])$ is Zollfrei, we thus obtain the following: 

\begin{lem}\label{any} 
Let $(M,[g])$ be a  space-time-orientable  Zollfrei self-dual manifold.
Then any two $\beta$-surfaces in $M$ have non-empty intersection. 
Moreover, any two distinct $\beta$-surfaces meet in exactly  $m$ points, where
the homological self-intersection of  any $\beta$-surface $S$ is 
given by 
$ [S] \cdot [S] = -m <0$. 
\end{lem}
\begin{proof}
Let $S$ be a reference $\beta$-surface, and suppose that we wish to 
understand the intersection of two given 
$\beta$-surfaces $S_1$ and $S_2$. If 
 $S_1=S_2$, they certainly intersect, and 
there is nothing to prove. Otherwise, Lemma \ref{homotopic}
tells us that $[S_1]=[S_2]=[S]$, and Lemma \ref{insect}
then yields $\# (S_1\cap S_2 )= - [S_1]\cdot [S_2]= -[S]\cdot [S]$.
In particular, the number $m$ of points of interesection is independent of {\em which} 
pair of distinct $\beta$-surfaces we choose to consider. But
since every $\beta$-plane is tangent to a $\beta$-surface, 
and since we have a circle's worth of different $\beta$-planes in 
each tangent space $T_xM$, 
we can certainly find   pairs $(S_1, S_2)$ with $S_1\neq S_2$ and $S_1\cap S_2\neq \emptyset$.
Thus $m > 0$, and we are done.   
\end{proof}
 
 \begin{lem} \label{compact}
 If  $(M,[g])$ is a Zollfrei self-dual $4$-manifold, then $M$ is compact. 
 \end{lem}
 \begin{proof}
 By passing to a double cover if necessary, we may assume that 
 $(M,g)$ is space-time-oriented. 
 
 Fix a reference $\beta$-surface $S\subset M$. Then for any point 
  $x\in M$,   there
 is a $\beta$-surface through $x$ that meets $S$; indeed, there certainly
 {\em are} $\beta$-surfaces through $x$, and Lemma \ref{any} tells us that
 {\em any} of these must meet $S$. But this statement can be rewritten as 
 the assertion that 
 $$M= {\zap p}\left[{\zap q}^{-1}\left({\zap q}\left[{\zap p}^{-1}(S)\right]\right)\right].$$
 Since ${\zap p}$ and ${\zap q}$ are both proper  maps, and since $S$ is compact,
 it therefore follows that $M$ is compact, too. 
 \end{proof}

If $(M,g)$ is a space-time-oriented Zollfrei $4$-manifold,  each fiber ${\zap  p}^{-1}(x)$ 
of $F\to M$ is an oriented circle, and  its image
${\zap  q}[{\zap  p}^{-1}(x)]$ 
 in $P$ may be thought of as
 a map $\gamma_x : S^1 \to P$,  which  we will  call a {\em standard loop}. 

\begin{prop} \label{structure} Let $M$ be  space-time-orientable self-dual
Zollfrei, and let $P$ be its space of 
$\beta$-surfaces. Then $P$ is diffeomorphic to $\RP^3$, and 
  $\pi_1 (P) \cong \ZZ_2$ is generated by any  
standard loop $\gamma_x$, $x\in M$. 
Moreover, any two distinct $\beta$-surfaces in $M$ meet in exactly two points. 
 \end{prop}
\begin{proof}
Let $y\in P$ be any base point, and let $S\subset M$ be the
corresponding $\beta$-surface. Every other $\beta$-surface in $M$
meets $S$ in $m$ distinct  points, where $[S]\cdot [S] = -m$. Moreover, through every
point of $S$, there passes a circle's worth of $\beta$-surfaces, only
one of which is $S$. Now recall that the structure group of ${\zap p}$ is $O(1,2)=PSL(2,\RR)$. 
 Thus, by removing the one point representing $TS$ from each fiber of 
${\zap  p}^{-1} (S) \to S$, we obtain  an affine 
$\RR$-bundle $L$ over $S\approx S^2$ which, via ${\zap  q}$,   maps locally diffeomorphically
onto 
$P-\{ y \}$ in an $m$-to-$1$ fashion.  Since any affine $\RR$-bundle over $S^2$ is trivial,
it follows that 
universal cover of $P-\{ y \}$ is $L\approx S^2\times \RR= S^3 -\{ 2 \mbox{ points}\}$; 
and since the order of this covering is $m=-[S]\cdot [S]$, we also see that 
$|\pi_1(P-\{ y\})|=m$. 
But $\pi_1(P)= \pi_1(P-\{ y\})$, since removing a point from 
a $3$-manifold doesn't change its fundamental group.
The universal cover of $P-\{ y \}$ is therefore gotten from the
universal cover $\tilde{P}$ of $P$ by removing $|\pi_1(P)|=m$
points. Since the universal cover of $P-\{ y \}$ has 
$m$ ends, whereas $S^2\times \RR$ has just two ends, it follows 
that $m=2$. Thus
$P= S^3 /\ZZ_2$ for some free $\ZZ_2$-action, and a theorem of Livesay 
\cite{live}  tells us  that $P\approx \RP^3$. 

Finally, notice that  the fiber of $L\to S$
defines a lift $\tilde{\gamma}_x$ 
of $\gamma_x$  to $\tilde{P}\approx S^3$. Since this lift is  
 not a loop, but rather is  a curve 
joining the two pre-images of $y$, it follows  that $\gamma_x$ is non-trivial in $\pi_1 (P)$.
Thus $[\gamma_x]$ generates $\pi_1 (P)\cong\pi_1 (\RP^3)\cong 
\ZZ_2$, and we are done. 
\end{proof}

Imitating the proof of Lemma \ref{unisphere} now gives us the following: 

\begin{lem}\label{folia}
 If  $(M,[g])$ is space-time-orientable and self-dual
Zollfrei, then   $F$ is diffeomorphic to $\RP^3\times S^2$ in  such a manner
 that ${\zap  q}$ becomes  the first-factor projection
  $\RP^3\times S^2\to \RP^3$. 
\end{lem}
\begin{proof}
Let $\varpi: S(TP) \to P$ denote the sphere bundle  defined by
$$S(TP)= (TP-0_P)/\RR^+,$$ where $0_P\subset TP$ denotes the zero section, and
where the positive  reals  $\RR^+$ act on $TP$ 
by scalar multiplication. That is,  $S(TP)$ may be thought of as the 
unit tangent bundle of $P$ for any choice of Riemannian metric on $P$. 

Let $u$ be a non-zero vector field on $F$ which spans $\ker {\zap p}_*$ 
at each point; this is possible because the choice of a
 metric-compatible decomposition $TM= T_-\oplus T_+$
 allows one to 
realize ${\zap p}: F\to M$ as the principle $SO(2)$-bundle of orientation-reversing 
isometries $T_-\to T_+$. 
 Since the fibers
of ${\zap p}$ and ${\zap q}$ are nowhere tangent,
we can therefore  define a map 
\begin{eqnarray*}
\Phi: F & \to & S (TP) \\
z & \mapsto  &  [{\zap q}_*u]
\end{eqnarray*}
which makes the diagram 
\setlength{\unitlength}{1ex}
\begin{center}\begin{picture}(40,17)(0,3)
\put(10,17){\makebox(0,0){$F$}}
\put(18,19){\makebox(0,0){$\Phi$}}
\put(18,5){\makebox(0,0){$P$}}
\put(12,11){\makebox(0,0){$\zap q$}}
\put(24,11){\makebox(0,0){$\varpi$}}
\put(26,17){\makebox(0,0){$S (TP)$}}
\put(11,15.5){\vector(2,-3){6}}
\put(25,15.5){\vector(-2,-3){6}}
\put(12,17){\vector(1,0){10}}
\end{picture}\end{center}
commute. Over each point of $P$, the map $\Phi$
is just the lift  of the developing map $\fye : S \to \RP^2$ constructed in Proposition \ref{projflat}
to the universal cover $S^2$ of $\RP^2$, and  so is a  diffeomorphism. 
Hence $\Phi$ is a bijection. Moreover, 
since $\zap q$ and $\varpi$ are
both submersions, it follows that $\Phi_*$ has maximal rank everywhere,
and $\Phi$ is therefore a diffeomorphism. 
However, $P\approx \RP^3$  is parallelizable, 
so $F\approx S (TP) \approx  \RP^3 \times S^2$, as claimed.
\end{proof}

\begin{thm}\label{s2xs2} If $(M,[g])$ is a space-time-orientable self-dual
Zollfrei $4$-manifold, 
then $M$ is homeomorphic to $S^2 \times S^2$.  
\end{thm}
\begin{proof}
Since the standard loop 
$\gamma_x={\zap  q}[{\zap  p}^{-1}(x)]$ generates
$\pi_1(P)$, the pull-back map 
$${\zap  q}^*: H^1(P, \ZZ_2)\to
H^1 (F, \ZZ_2)$$ sends the generator  of 
$H^1(P, \ZZ_2)\cong \ZZ_2$ to an element 
of $H^1(F, \ZZ_2)$ which 
is non-trivial on the fiber class of ${\zap  p}$. 
This shows that  there is a double cover  $\tilde{F}\to F$ which restricts to a 
double cover  $S^1\to S^1$  of each fiber of ${\zap  p}$. 

Now choose any $g$-adapted, orientation compatible  almost-complex structure
${\mathfrak J}$ on ${M}$. The $S^1$-bundle ${\zap  p}: F\to M$
can then be identified with the unit circle bundle  of the 
canonical line bundle $K$ of $({M}, {\mathfrak J})$. 
The double cover $\tilde{F}\to F$ then becomes 
the unit circle bundle of a square-root $K^{1/2}$
of $K$. Hence $c_1({M}, {\mathfrak J})$ is divisible
by $2$ in $H^2(M,\ZZ)$. Because  
$w_2(M)$ is the mod-$2$ reduction \cite{milnorstaf} of $c_1({M}, {\mathfrak J})$,
and because the sequence 
$$\cdots \to 
H^2(M,\ZZ ) \stackrel{2}{\longrightarrow}
H^2(M,\ZZ )\to H^2(M,\ZZ_2 )\to \cdots $$
is exact, it follows that 
$w_2(M)=0$. Thus $M$ is a spin manifold. 

Since $F\approx \RP^3 \times S^2$, its universal cover must
be $\tilde{F}\approx S^3\times S^2$. Hence   the long exact 
homotopy sequence   \cite{span} 
$$\cdots \to  \pi_2 (S^1) \to \pi_2(\tilde{F})\to \pi_2 (M) \to 
\pi_1 (S^1)\to \pi_1 (\tilde{F})\to \pi_1 (M)\to 0$$
of the fibration $S^1\to \tilde{F}\to M$ 
now tells us that $\pi_1(M)=0$ and $\pi_2(M) = \ZZ\oplus \ZZ$. 
Thus  $M$ is a simply connected compact  $4$-manifold with $b_2=2$
and even intersection form. 
The Freedman classification of simply connected $4$-manifolds \cite{freedman}
therefore tells us that $M$ is homeomorphic to $S^2 \times S^2$. 
\end{proof}

In fact, it seems reasonable to  conjecture that any 
space-time-orientable self-dual
Zollfrei $4$-manifold must actually be {\em diffeomorphic} to 
$S^2\times S^2$. While we have not managed 
to prove this stronger statement in  general, 
we will eventually see,  in Theorem \ref{bihol} below, 
that 
it at least turns out to be true if $[g]$ is represented by an
indefinite   K\"ahler metric.
\bigskip 

We now turn  to the {non}-space-time-orientable case. 

\begin{prop}\label{rp2}
Let $(M,[g])$ be a  Zollfrei  self-dual
$4$-manifold which is {\em not} space-time-orientable. 
Then every $\beta$-surface in $M$ is an embedded $\RP^2$,
and every pair of distinct $\beta$-surfaces intersects in exactly
one point. 
\end{prop}
\begin{proof}
Notice that our definition of self-duality requires that $M$ be
orientable. Thus the set $\tilde{M}$ of orientation-compatible local 
space-time orientations 
of $(M,[g])$ is a  double cover of $M$. Notice that 
$\tilde{M}$ is space-time-orientable and self-dual Zollfrei
with respect to the pulled back metric. 
Let $a: \tilde{M}\to  \tilde{M}$ be the non-trivial 
deck transformation. 

If $S\subset \tilde{M}$ is any $\beta$-surface, then we claim that 
$a [S] = S$. Indeed, suppose not. Then $a [S] = S^\prime$
would be a different $\beta$-surface, and  hence 
$S\cap S^\prime$ would consist of exactly two points
by Proposition \ref{structure}; and since 
$a[S\cap S^\prime ] = a [S] \cap a[a[S]] = S^\prime \cap S$,
these two points would necessarily be interchanged by 
the fixed-point-free involution $a$. On the other hand, all
the other points of $S$ would be moved to the complement of 
$S$ by $a$. Hence the image of $S$ in $M= \tilde{M}/\langle a \rangle$
would be an imersed sphere with a single self-intersection. 
But this contradicts Lemma \ref{beta}. Thus every $\beta$-surface
in $\tilde{M}$ must be sent to itself by $a$.

It follows that every $\beta$-surface in $\tilde{M}$ is the double cover
of a $\beta$-surface in $M$. Since all  the $\beta$-surfaces in $\tilde{M}$
are $2$-spheres by Lemma \ref{sphere}, 
and since every $\beta$-surface in $M$ must be the image of a 
$\beta$-surface in $\tilde{M}$, it follows that every $\beta$-surface
in $M$ must be an 
 $\RP^2$. Moreover, since $\tilde{M}\to M$ is a double cover, and since
 $\beta$-surfaces  in $\tilde{M}$ intersect in pairs of 
 points, pairs of distinct $\beta$-planes in $M$ must always 
 intersect in a unique point. 
 \end{proof}

 \begin{thm}\label{tastier}
 Let $(M,[g])$ be a 
   self-dual split-signature $4$-manifold.
   Then the following are equivalent:
   \begin{romenumerate}
   \item  $(M,[g])$ is Zollfrei;
     \item  $(M,[g])$ is strongly Zollfrei;
\item exactly  one of the following holds: 
  \begin{alphenumerate}
  \item every $\beta$-surface is an embedded $S^2\subset M$; or 
  \item every $\beta$-surface is an embedded $\RP^2\subset M$.
  \end{alphenumerate}
  \end{romenumerate}
  \end{thm}
  \begin{proof}
 Notice that $(iii)\Longrightarrow (i)$ by Proposition \ref{projflat} and the uniqueness
  \cite{nico} 
 of the flat projective structures on $\RP^2$ and $S^2$. 
  Thus Lemma \ref{sphere} and Proposition \ref{rp2} tell us that 
  $(iii)(a)$ can only occur if $(M,[g])$ is space-time orientable,
  whereas $(iii)(b)$ can only occur  if $(M,[g])$ is {\em not} space-time orientable. 
  
  If $(M,[g])$ is space-time orientable, 
   the desired equivalence is therefore given by Theorem \ref{tasty}.

  If, on the other hand, $(M,[g])$ is not space-time orientable, then 
  $(ii)\Longrightarrow (i)$ by Definition \ref{guil},
  and $(i)\Longrightarrow (iii)$ by Proposition \ref{rp2}. 
  On the other hand, $(iii)(b)\Longrightarrow (ii)$, too. Indeed, 
   the space-time-orientable 
  double cover $\tilde{M}$  of $M$ is Zollfrei, and hence strongly Zollfrei  by  Theorem \ref{tasty}.
  The non-trivial deck transformation $a$ of $\tilde{M}\to M$ must therefore send each 
  null geodesic to itself by the uniqueness \cite{nico} 
  of the flat projective structure on $\RP^2$.   \end{proof}
 
 Proposition \ref{rp2} also allows us to deduce the following: 
 
 \begin{lem}\label{nonspin}
Let $(M,[g])$ be a  Zollfrei  self-dual
$4$-manifold which is {not} space-time-orientable. 
Then $M$ is non-spin. 
 \end{lem}
 \begin{proof}
 Let 
  ${\mathbf b}\in H^2(M,\ZZ_2 )$
  denote the Poincar\'e dual of the $\ZZ_2$-homology class of any $\beta$-surface 
  $S\subset M$. Since any two distinct $\beta$-surfaces are freely homotopic and
  intersect transversely in exactly one  
   point, we have 
  ${\mathbf b } \cdot {\mathbf b}= 1\in \ZZ_2$, 
  where 
  $$\cdot :   H^2(M,\ZZ_2 ) \times H^2(M,\ZZ_2 ) \to \ZZ_2$$
  is the intersection form of $M$ with $\ZZ_2$ coefficients. 
  But since $M$ is orientable, Wu's formula  \cite{hirzwu} asserts that $w_2(M)$ satisfies 
  $$w_2  \cdot {\mathbf x}= {\mathbf x} \cdot {\mathbf x}$$
  for any ${\mathbf x}\in H^2 (M, \ZZ_2),$
  so we have
  $$w_2 \cdot {\mathbf b} = {\mathbf b } \cdot {\mathbf b}= 1.$$
  Thus  $w_2(M)\neq 0$, and $M$ is non-spin, as claimed. 
 \end{proof}
 
 \begin{thm} \label{quadric} 
Let $(M,[g])$ be a  Zollfrei  self-dual
$4$-manifold which is {not} space-time-orientable. 
Then $M$ is homeomorphic to the real projective
quadric ${\mathbb M}^{2,2}$. 
 \end{thm}
 \begin{proof}
 Freedman's topological classification of  simply connected $4$-manifolds
 has been extended to compact oriented $4$-manifolds with finite cyclic 
 fundamental group by Hambleton and Kreck \cite[Theorem C]{hmbkrck}. 
 They show that such manifolds are classified up to homeomorphism by 
 their fundamental groups,  their
 intersection forms on $H^2(\bullet  , \ZZ)/\mbox{torsion}$,
 their $w_2$-types, 
 and 
 their Kirby-Siebenmann  invariants.
The Kirby-Siebenmann invariant   vanishes if a manifold
admits a smooth structure.  The $w_2$-type of a 
$4$-manifold says whether the manifold and its
universal cover are spin; 
namely,  an oriented manifold
$M$ with universal cover $\tilde{M}$ is said to be of type (I) if 
$w_2(\tilde{M})\neq 0$, type (II) if 
$w_2(M)=0$,  and type (III) if $w_2(\tilde{M})=0$, but 
 $w_2(M)\neq 0$. 

Now assume that   $(M,[g])$ is a non-space-time-orientable 
self-dual Zollfrei manifold. Then 
 $M$ is smooth,  and so has vanishing Kirby-Siebenmann invariant.
 Also, $M$  is  oriented, as  is required by our definition of 
 self-duality. Now recall that the double cover $\tilde{M}$ of $M$ by its local space-time
 orientations is a space-time orientable  Zollfrei self-dual $4$-manifold,
 and so is homeomorphic to  $S^2\times S^2$ by Theorem \ref{s2xs2}.
 Since $S^2\times S^2$ is simply connected, $\tilde{M}$ is actually the
 universal cover, and we therefore have $\pi_1(M)= \ZZ_2$. 
 Moreover, the Euler characteristic of $M$ must be $\chi (M)= \chi(S^2\times S^2)/2=2$,
 so 
  $H^2(M,\ZZ)/\mbox{torsion}=0$, and the intersection form of
 $M$ must therefore be trivial. Finally, $w_2(M)\neq 0$ by 
 Proposition \ref{nonspin}, whereas $w_2(\tilde{M})=w_2(S^2\times S^2)=0$,
 so $M$ is of type (III). Hambleton and Kreck therefore tell us that there
 is only one possible homeomorphism type for such an $M$. 
 The quadric ${\mathbb M}^{2,2}\subset \RP^5$ 
 therefore represents the only topological possibility.  
 \end{proof}

 Combining   Theorems \ref{s2xs2} and \ref{quadric}, we have thus proved  
 {\bf Theorem 
 \ref{topprop}}. 
 
\section{Stability of the Zollfrei Condition}
 
We now turn to the important assertion that the Zollfrei condition is open 
among self-dual metrics. This phenomenon is actually a manifestation 
of aspects of the theory of foliations arising from  Thurston stability for   
compact leaves of foliations \cite{thurstfol}. 
The  result we will need  is originally due to Langevin and Rosenberg \cite{stability},
although the formulation  given here is actually that of Epstein and 
Rosenberg \cite{epros}. 

\begin{thm}[Langevin-Rosenberg] \label{lep}
Let $\mbox{\cyr  p} : X\to Y$ be a $C^1$ fiber bundle with compact fibers and compact base, 
where the fibers of $\mbox{\cyr  p}$ have $b_1=0$ over $\RR$. 
Let $\mathfrak F$ be the foliation of $X$ by the fibers
of $\mbox{\cyr  p}$. Then the foliation $\mathfrak F$ has a neighborhood $\mathscr V$ in the 
$C^1$ Epstein topology on the space of foliations of $X$ such that 
every foliation ${\mathfrak F}^\prime\in {\mathscr V}$  is of the 
form $\phi^* {\mathfrak F}$ for some $C^1$-diffeomorphism 
$\phi : X\to X$. 
\end{thm}

Here two $C^1$ foliations of $X$ are close in the  $C^1$ {\em Epstein topology}  \cite{eptop}  if 
there are finite atlases of trivializing charts for the two foliations which are close
in the usual $C^1$ topology on the space of maps. The only thing that need
concern us here is that two $C^1$ integrable distributions of $k$-planes which are $C^1$ close
as sections of the  Grassmann bundle $Gr_k(TX)\to X$ define foliations
which are close in Epstein's sense. 

Combining Theorem \ref{lep}  with our results from  \S \ref{prelim},
we thus obtain 

\setcounter{main}{0}
\begin{main}
Let $(M,g)$ be a self-dual Zollfrei $4$-manifold. 
Then any other self-dual metric $g'$ on $M$ that is sufficiently  close to $g$ 
in the $C^2$ topology is also Zollfrei. 
\end{main} 
\begin{proof}
There is a $C^0$ neighborhood of $g$ in the pseudo-Riemannian
metrics in which every metric $g^\prime$ can be written  as
$g^\prime= A^*g$ for a unique $g$-self-adjoint endomorphism $A: TM\to TM$
which is $C^0$  close to the identity. This endomorphism of 
$TM$ allows one to identify the pseudo-orthonormal frame 
bundles of $g$ and $g^\prime$. Moreover, 
if $g^\prime$ is $C^2$ close to $g$, the corresponding principle connections 
are then $C^1$-close after this correspondence has been made. 
Using $A$ to identify the bundle of $\beta$-planes for $g^\prime$
with the bundle ${\zap p}: F\to M$ of $\beta$-planes for $g$,
we then obtain two distributions $E$ and $E^\prime$ on $F$ which represent the 
horizontal lifts of the $\beta$-planes of $g$ and $g^\prime$, respectively;
and these two distributions will be $C^1$ close if we again assume that 
$g$ and $g^\prime$ are $C^2$ close. 

Now if $g$ and $g^\prime$ are both self-dual, the distributions $E$ and 
$E^\prime$ will both be integrable, and will be tangent to foliations
${\mathscr F}$ and ${\mathscr F}^\prime$ that represent the
canonical lifts of the $\beta$-surfaces of the two metrics. Moreover, 
 ${\mathscr F}^\prime$ will be $C^1$ close to $\mathscr F$ if we
assume that $g^\prime$ is $C^2$ close to $g$. But if, in addition, 
$g$ is Zollfrei, the leaves of the foliation $\mathscr F$ will 
exactly be the fibers of a fiber bundle ${\zap q}: F\to P$.
Now $F$ is necessarily compact by Lemma \ref{compact}, 
while Theorem \ref{tastier} tells us that  the fibers of ${\zap q}$  are spheres or
projective planes. Since these are compact surfaces
with $b_1=0$, we may 
 therefore apply Theorem \ref{lep}
to conclude that there is a $C^1$ diffeomorphism 
$\phi : F\to F$ which sends $\mathscr F$ to ${\mathscr F}^\prime$.
Thus, if $g$ is self-dual 
and Zollfrei, and if $g^\prime$ is self-dual and $C^2$ close to $g$,
then the $\beta$-surface of $g^\prime$ are either all spheres
or all projective planes, and Theorem \ref{tastier} therefore tells us that 
$g^\prime$ is Zollfrei, too, as claimed. 
\end{proof}

\section{Constructing the Twistor Space}

At this point, 
we have already achieved a certain level of intimacy with the  bundle ${\zap p} : F\to M$  of {real} 
$\beta$-planes over an oriented  split-signature
 conformal $4$-manifold
 $(M,[g])$.
It is now time to introduce  the bundle $\wp : {\mathcal Z}\to M$ of 
{\em complex} $\beta$-planes. Just as in the real case, a $2$-dimensional complex subspace 
$\Pi$ of a complexified tangent space $T_\CC M|_x= \CC\otimes T_xM$  of $M$
is called {\em isotropic} if the complex-bilinear extension of $g$ 
vanishes when restricted to $\Pi$. Such isotropic planes come in 
two flavors. The complex $\alpha$-planes are precisely those complex
$2$-planes $\Pi$ such that $\wedge^2 \Pi$ corresponds by index lowering to a complex
null line in $\Lambda^+_\CC$; the complex $\beta$-planes instead
correspond to null $1$-dimensional subspaces of $\Lambda^-_\CC$. 
Thus, the bundle of complex $\beta$-planes on $M$ is exactly given by 
$${\mathcal Z}= \{ [\varphi ] \in \PP (\Lambda^+_\CC)~|~ 
\langle \varphi , \varphi \rangle =0\},$$
where $\langle \varphi , \psi \rangle = \frac{1}{2}\varphi_{ab}\psi_{cd}g^{ac}g^{bd}$
is the complex-bilinear extension to $\Lambda^2_\CC$ of the  inner-product on 
$2$-forms induced by $g$.   
Since $\PP (\Lambda^+_\CC)$ is a $\CP_2$-bundle over $M$,
each fiber of $\mathcal Z$ is a non-degenerate conic in $\CP_2$, 
and so is intrinsically a $\CP_1$. Indeed, $\mathcal Z$ is precisely the 
$\CP_1$-bundle obtained from $F\to M$ by remembering
that $F$ has structure group $PSL(2, \RR)$, and
that one can therefore construct an associated  $\CP_1$-bundle 
over $M$ by including $PSL(2,\RR)$ in $PSL(2,\CC)$ and considering the 
standard action of $PSL(2,\CC)$ on $\CP_1$. In particular, each fiber
of $\wp : {\mathcal Z}\to M$ is a holomorphic curve. Let
${\mathzap V}^{0,1}\subset T_\CC {\mathcal Z}$ be the $(0,1)$-tangent bundle
of the fibers. 

Fix a metric $g$ in the conformal class, and notice that $g$ determines
a connection on $\mathcal Z$, in the sense that $g$ determines
a notion of parallel transport of elements of $\mathcal Z$ along
smooth curves in $M$. Let ${\mathzap H} \subset T{\mathcal Z}$ be the 
horizontal subspace of this connection, so that  the derivative
of the projection gives us a canonical isomorphism 
$\wp_*: {\mathzap H}\to \wp^* TM$. Let ${\mathzap H}_\CC= {\mathzap H}\otimes\CC$.
Then $\mathcal Z$ carries a unique distribution 
${\mathzap E}\subset {\mathzap H}_\CC\subset T_\CC {\mathcal Z}$ of
horizontal complex $2$-planes  such that 
$$\wp_* ({\mathzap E}|_\Pi)= \Pi \subset T_\CC M.$$
Set 
$${\mathzap D}= {\mathzap E}+ {\mathzap V}^{0,1}.$$
Since ${\mathzap E}$ is horizontal and ${\mathzap V}^{0,1}$ is vertical, this
sum is in fact a direct sum, and $\mathzap D$ is therefore a 
distribution of complex $3$-planes on $\mathcal Z$. 

Let us make this discussion more concrete by temporarily 
restricting our attention to an open subset ${\mathscr U}\subset M$
on which we can find an oriented  pseudo-orthonormal frame 
$e_1, \ldots , e_4$ with 
$$g(e_j, e_k) =
\left\{ \begin{array}{rl}
0& \mbox{ if } j\neq k,\\
    1 &\mbox{ if }   j=k\in\{1, 2\}, \\
      -1&\mbox{ if }   j=k \in \{3, 4\}.
      \end{array}
\right.
$$
We remark that if $g$ is of differentiability 
class $C^{k}$, then such frames $e_1, \ldots , e_4$
of class $C^{k}$ can locally be constructed by means of the Gramm-Schmidt procedure. 
This in turn determines a pseudo-orthonormal basis for $\Lambda^{-}|_{\mathscr U}$ by setting
\begin{eqnarray*}
\varphi_1 & = & \frac{1}{\sqrt{2}}(e^1\wedge e^2- e^3\wedge e^4) \\
\varphi_2& = & \frac{1}{\sqrt{2}}(e^1\wedge e^3  -e^2\wedge e^4) \\
\varphi_3& = & \frac{1}{\sqrt{2}}(e^1\wedge e^4 + e^2\wedge e^3) 
\end{eqnarray*}
so that 
$$
\langle \varphi_{\mathzap j} , \varphi_{\mathzap k} \rangle = 
\left\{ \begin{array}{rl}
0& \mbox{ if } {\mathzap j}\neq {\mathzap k},\\
    1 &\mbox{ if }   {\mathzap j}={\mathzap k}=1 \\
      -1&\mbox{ if }   {\mathzap j}={\mathzap k} \in \{2, 3\}.
      \end{array}
\right.
$$
We  can then identify $\CP_1\times {\mathscr U}$ with 
$\wp^{-1}({\mathscr U})\subset {\mathcal Z}$ by 
$$([\zeta_1: \zeta_2], x) \longmapsto \left. 
\left[(\zeta_1^2+ \zeta_2^2)~\varphi_1 +  (\zeta_1^2- \zeta_2^2)~
\varphi_2 - 2\zeta_1\zeta_2 ~\varphi_3 \right] \right|_x ,$$
and it is worth noting that in the process we have identified 
$\RP^1\times {\mathscr U}$ with ${\zap p}^{-1}({\mathscr U})\subset F \subset {\mathcal Z}$.
In particular,  an  open dense subset of $\wp^{-1}({\mathscr U})$ may be parameterized by 
 $\CC  \times  {\mathscr U}$, via the map
 $$(\zeta , x ) \longmapsto 
 [(1+\zeta^2)~\varphi_1 +  (1-\zeta^2)~
\varphi_2 - 2\zeta ~\varphi_3  ]\Big|_x,$$
and in the process we sweep out an open dense subset of 
${\zap p}^{-1}({\mathscr U})$ with $\RR \times {\mathscr U.}$ 
Notice that for each $(\zeta, x)$ with $\zeta\neq \pm i$,  the corresponding $\beta$-plane is 
exactly
$$\Pi = \mbox{span}\left. 
\left\{ (\zeta^2+1)e_1 -2\zeta e_3 + (\zeta^2-1)e_4 ~,~ 
(\zeta^2+1)e_2+ (\zeta^2-1)e_3 + 2\zeta  e_4
\right\} \right|_x.
$$

Now  observe that we  have 
$$\nabla \varphi_{{\mathzap j}} = \theta_{{\mathzap j}}^{{\mathzap k}} \otimes
 \varphi_{{\mathzap k}},$$
for an ${\mathfrak s \mathfrak o}(1,2)$-valued $1$-form 
$[\theta_{{\mathzap j}}^{\mathzap k}]$:
$$
\theta^1_2= \theta^2_1, \hspace{.5cm} \theta^1_3= \theta^3_1, \hspace{.5cm} \theta^2_3= -\theta^3_2,
\hspace{.5cm}\theta^1_1=\theta^2_2=\theta^3_3= 0.
$$
When we then expand these $1$-forms 
 as $\theta_{{\mathzap j}}^{\mathzap k}=\theta_{{\mathzap j}\ell}^{\mathzap k} e^\ell$
 the resulting 
 functions $\theta_{{\mathzap j}\ell}^{\mathzap k} $ are just linear combinations 
of the components of the usual connection symbols of the frame, and so are of class $C^{k-1}$
if our  frame is of class $C^{k}$. 
The distribution $\mathzap D$ now  becomes
$$
{\mathzap D} = \mbox{span} \left\{ 
 {\mathfrak w}_1, {\mathfrak w}_2, \frac{\partial}{\partial\overline{\zeta}}
\right\} $$
on $(\CC -\{ \pm i\}) \times  {\mathscr U}$, where the vector fields 
\begin{eqnarray*}
{\mathfrak w}_1 & = & (\zeta^2+1)e_1 -2\zeta e_3 + (\zeta^2-1)e_4 + Q_1(x,\zeta) \frac{\partial}{\partial \zeta}\\
{\mathfrak w}_2  & = &  (\zeta^2+1)e_2+ (\zeta^2-1)e_3 + 2\zeta  e_4+ Q_2(x,\zeta) \frac{\partial}{\partial \zeta}
\end{eqnarray*}
are defined in terms of the  functions 
\begin{eqnarray*}
Q_1(x,\zeta) & = &
\frac{1-\zeta^2}{2}\Big[  (\zeta^2+1)\theta^3_{11} - 2\zeta\theta^3_{13} + (\zeta^2-1)\theta^3_{14}
\Big]\\ &&\hspace{1cm}+
\zeta \Big[  (\zeta^2+1)\theta^2_{11}- 2\zeta\theta^2_{13}+ (\zeta^2-1)\theta^2_{14} 
\Big]\\ &&\hspace{2cm}
-\frac{1+\zeta^2}{2}\Big[(\zeta^2+1)\theta^2_{31} - 2\zeta\theta^2_{33} + (\zeta^2-1)\theta^2_{34}
\Big] 
\\
Q_2(x,\zeta)  & = & 
\frac{1-\zeta^2}{2}\Big[ (\zeta^2+1) \theta^3_{12} +  (\zeta^2-1)\theta^3_{13} + 2\zeta  \theta^3_{14}
\Big]\\ &&\hspace{1cm}+
\zeta \Big[ (\zeta^2+1) \theta^2_{12}+  (\zeta^2-1)\theta^2_{13} + 2\zeta  \theta^2_{14}
\Big]\\ &&\hspace{2cm}
-\frac{1+\zeta^2}{2}\Big[(\zeta^2+1)\theta^2_{32} + (\zeta^2-1) \theta^2_{33} + 2\zeta  \theta^2_{34}
\Big] 
\end{eqnarray*}

The minuti{\ae}   of these expressions are of little importance, but three
facts are  worthy of emphasis. First of all, 
the components of ${\mathfrak w}_1$ and ${\mathfrak w}_2$ 
in the basis $e_1, \ldots , e_4, \partial/\partial \zeta$ are polynomial
in $\zeta$ for any  fixed 
 $x\in {\mathscr U}$, and so, in particular, 
$$
\left[ \frac{\partial}{\partial \overline{\zeta}}, {\mathfrak w}_1\right] =
\left[ \frac{\partial}{\partial \overline{\zeta}}, {\mathfrak w}_2\right] = 0.
$$
Secondly,   we have chosen the vector fields 
${\mathfrak w}_1$ and ${\mathfrak w}_2$
to be real and {\em horizontal}
along the locus $F$ where $\zeta$ is real\footnote{The  ${\mathfrak w}_j$
 could only be forced to be horizontal {\em everywhere}
 at the price of   adding  muliples of $\partial/\partial \overline{\zeta}$ to them. We have  
  avoided 
  doing so here  because    the  relevant coefficients 
would generally {\em not} be holomorphic 
in $\zeta$, and the Lie brackets of the  
${\mathfrak w}_j$ with $\partial/\partial \overline{\zeta}$
would therefore no longer vanish.} .
Finally, notice  that $\mathzap D$ is spanned by $C^{k-1}$ vector
fields if $g$ is of class $C^{k}$.

\begin{prop}\label{critter}
Let $(M,g)$ be an  oriented   split-signature $C^2$
pseudo-Riemannian  $4$-manifold. 
Let  $\wp : {\mathcal Z}\to M$ be the bundle of 
complex $\beta$-planes in $T_\CC M$, and let 
 ${\mathzap D}\subset T_\CC{\mathcal Z}$ be the $C^1$ distribution of
 complex $3$-planes defined
above. Then ${\mathzap D}$ is involutive, in the sense
that 
$$[C^1({\mathzap D}), C^1({\mathzap D})]\subset C^0({\mathzap D}),$$
iff $(M,g)$ is self-dual. 
\end{prop}
\begin{proof}
Let us begin by noticing that 
$${\mathzap D}\cap T_\CC F= {\mathzap E}|_F= E\otimes \CC ,$$
where the real distribution  of $2$-planes $E$ on $F$  is  defined on page 
\pageref{dodger}. Also recall that   Proposition
\ref{roger} tells us that $E$ is Frobenius integrable iff $g$ is self-dual. 

Now, supose that ${\mathzap D}$ is involutive. Then both $T_\CC F$ and 
$\mathzap D$ are closed under Lie brackets. Hence 
${\mathzap D}\cap T_\CC F=  E\otimes \CC$
is closed under Lie brackets, too. Thus $E$ is Frobenius integrable, and 
 Proposition
\ref{roger} therefore tells us that $g$ is self-dual. 

Conversely, suppose that $g$ is self-dual. Then  Proposition
\ref{roger} tells us that $E\to F$ is involutive. 
Let  ${\mathscr U}\subset M$ be any open set on which there exists 
 a pseudo-orthononormal frame $e_1, \ldots, e_4$, and 
consider the vector fields ${\mathfrak w}_1$ and ${\mathfrak w}_2$ 
constructed on an open dense subset of $\wp^{-1}({\mathscr U})$ above. 
Along $F$, the vector fields ${\mathfrak w}_1$ and ${\mathfrak w}_2$ are linearly independent
sections 
of the involutive rank-$2$ bundle $E\subset TF$, so
$$\left[ {\mathfrak w}_1, {\mathfrak w}_2\right] \wedge {\mathfrak w}_1\wedge {\mathfrak w}_2
= 0 ~~~\mbox{  when   }  \zeta = \overline{\zeta}.$$
However, relative to the frame 
$e_1, \ldots e_4, \partial/\partial \zeta$, 
the components of ${\mathfrak w}_1$ and ${\mathfrak w}_2$
are polynomial in $\zeta$, so it follows that   the components
 of the tensor field $\left[ {\mathfrak w}_1, {\mathfrak w}_2\right]\wedge {\mathfrak w}_1\wedge
 {\mathfrak w}_2$ are polynomial in $\zeta$, too. But we have already seen
 that $\left[ {\mathfrak w}_1, {\mathfrak w}_2\right]\wedge {\mathfrak w}_1\wedge
 {\mathfrak w}_2$ vanishes when $\zeta$ is real. Hence $\left[ {\mathfrak w}_1, {\mathfrak w}_2\right]\wedge {\mathfrak w}_1\wedge
 {\mathfrak w}_2$ vanishes identically, and 
 we therefore have  
 $$
 \Big[ \frac{\partial}{\partial \overline{\zeta}} , {\mathfrak w}_1\Big], \Big[ \frac{\partial}{\partial \overline{\zeta}} , {\mathfrak w}_2\Big] ,  \Big[ {\mathfrak w}_1, {\mathfrak w}_2\Big] 
 \in \mbox{span}\left\{ \frac{\partial}{\partial \overline{\zeta}} , {\mathfrak w}_1, {\mathfrak w}_2\right\} .
 $$
  Thus $\mathzap D$ is  involutive
 on the region of $\wp^{-1}({\mathscr U})$ parameterized by $(\CC -\{ \pm i\}) \times  {\mathscr U}$,
 and   the O'Neill
 tensor
\begin{eqnarray*}
A_{\mathzap D}: {\mathzap D}\times  {\mathzap D}& \longleftarrow  & T_\CC {\mathcal Z}/{\mathzap D} \\
(u,v) &\mapsto  & [u,v] \bmod {\mathzap D} 
\end{eqnarray*}
therefore 
vanishes on an open dense subset of $\wp^{-1} ({\mathscr U})$. But $A_{\mathzap D}$ is 
continuous, so it therefore vanishes on all 
of $\wp^{-1} ({\mathscr U})$. Since such subsets ${\mathscr U}$ cover all of $M$,
it therefore follows that $\mathzap D$ is involutive on all of $\mathcal Z$. 
\end{proof}

Similar reasoning also shows the following:

\begin{prop} \label{conformal}
Let $(M,[g])$ be an  oriented   split-signature self-dual  
$4$-manifold. 
Then the involutive distribution ${\mathzap D}$  
on ${\mathcal Z}$ is conformally invariant ---
that is, it depends only on the conformal class $[g]$, rather than on the metric 
$g\in [g]$.   \end{prop}
\begin{proof}
Since multiplying $g$ by $-1$ does not change the metric connection, and
therefore does not change ${\mathzap D}= {\mathzap E}\oplus {\mathzap V}^{0,1}$,
it suffices to henceforth consider only conformally related 
pairs of metrics $g$ and $\hat{g}={\zap f}g$ for which  the factor  ${\zap f}$
is positive.

Now the distribution $E$ on $F$  only depends on $[g]$, since 
 it is tangent to the foliation $\mathscr F$ of $F$ by lifted $\beta$-surfaces. 
 Now consider two metrics $g$ and $\hat{g}={\zap f}g$ in $[g]$, where ${\zap f}> 0$. If 
 $e_1 , \ldots , e_4$ is a pseudo-orthonormal frame for $g$ on 
 an open subset ${\mathscr U}\subset M$,
 then ${\zap f}^{-1/2}e_1 , \ldots , {\zap f}^{-1/2}e_4$ is a 
 pseudo-orthonormal frame for $\hat{g}$.
 Let ${\mathfrak w}_j$ and $\hat{\mathfrak w}_j$
 be the vector fields on $(\CC -\{ \pm i\}) \times  {\mathscr U}$
 constructed from these two frames and metrics. Then 
 ${\mathfrak w}_j$ and ${\zap f}^{1/2}\hat{\mathfrak w}_j$ 
 coincide
 along $F$, since  they are sections of $E$ with the same projections. 
 But the components of ${\mathfrak w}_j$ and ${\zap f}^{1/2}\hat{\mathfrak w}_j$
 (expressed, say, as linear combinations of the $e_j$ and $\partial/\partial \zeta$) 
 are polynomial in $\zeta$. Since they coincide when $\zeta$ is real, 
 we must therefore have ${\mathfrak w}_j\equiv {\zap f}^{1/2}\hat{\mathfrak w}_j$.
 Hence the distribution $\mathzap D$ determined by $g$ coincides with 
 the distribution $\hat{\mathzap D}$ determined by $\hat{g}$ on an open dense
 subset of $\wp^{-1}(U)$, and we therefore have $\mathzap D\equiv  \hat{\mathzap D}$
 on $\wp^{-1}(U)$ by continuity. Since $M$ can be covered with such open sets
 ${\mathscr U}$, it therefore follows that $\mathzap D=  \hat{\mathzap D}$
 on all of $\mathcal Z$, as claimed. 
\end{proof}

Actually,  the conformal invariance of $\mathzap D$ holds even  in the absence of 
the self-duality hypothesis, but we will never need this fact. It is also worth remarking that 
Proposition \ref{critter}  could instead, for example, have  been proved by imitating the arguments 
of Atiyah-Hitchin-Singer \cite{AHS}. The route we have chosen is not arbitrary,
however, but rather is specifically intended to prepare the reader for the proof of 
Proposition \ref{machine} below. 

%\begin{cor} Let $(M,[g])$ be a self-dual split-signature $4$-manifold, 
%where $[g]$ contains some $C^2$ metric $g$. Then 
%\end{cor}

What is the `real' geometrical meaning of a point of the bundle $\wp : {\mathcal Z}\to M$?
Obviously, the points of $F\subset {\mathcal Z}$ are real totally null $2$-planes, and there is not much 
more to be said. By contrast, a point  of ${\mathcal Z}-F$ is a subspace $\Pi\subset T_x M\otimes \CC$
with the property that $\Pi \cap \overline{\Pi}=0$. Thus $\Pi\oplus \overline{\Pi}= T_x M\otimes \CC$,
and we can therefore define a unique almost-complex structure $\jmath : T_xM\to T_xM$
at $x$
by declaring that $\Pi$ is its $(+i)$-eigenspace. The requirement that $\Pi$ be isotropic
is then equivalent to the condition that $\jmath$ be an orthogonal transformation---
i.e. that $\jmath^*g= g$. Finally, the requirement that $\Pi$ be a $\beta$-plane,
rather than an $\alpha$-plane, is exactly that $\jmath$ determine the {\em given}
orientation of $M$, rather than the opposite one. This last requirement concretely
amounts to asking that there be an oriented pseudo-orthonormal basis 
$e_1, \ldots , e_4$ with $\jmath e_1= e_2$ and $\jmath e_3= e_4$.
Notice that this formulation implicitly is associated with a decomposition 
$T_xM= T_+\oplus T_-$, where $T_+= \mbox{span }\{e_1 , e_2\}$ 
and $T_-=\mbox{span} \{ e_3,e_4\}$, and that $\jmath$ gives us 
a specific orientation of the maximally positive and negative subspaces 
$T_+$ and $T_-$. 

Now suppose that  $(M,[g])$ is space-time orientable. 
It then follows that ${\mathcal Z}-F$ 
has two connected components, depending on whether the associated 
orientation on $T_-$ is the given one, or its reverse. Let $U\subset ({\mathcal Z}-F)$
be the open subset  corresponding to $\jmath$ for which the 
induced orientation on $T_-$ agrees with the previously chosen one. 
Then $\wp|_U:U\to M$ is an open disk bundle over $M$, and corresponds
to the region $\Im m ~\zeta > 0$ in our explicit local description of $\mathcal Z$. 
Let ${\mathcal Z}_+= U \cup F$ be the  closure of $U$ in $\mathcal Z$.
Thus ${\mathcal Z}_+$ is a compact $6$-manifold-with-boundary, 
and $\wp|_{{\mathcal Z}_+}: {\mathcal Z}_+\to M$ is a bundle of closed oriented $2$-disks.

Now  $F$ carries a foliation $\mathscr F$ by lifted $\beta$-surfaces. 
If we  assume that our space-time-oriented
self-dual $4$-manifold 
$(M,[g])$ is also {\em Zollfrei}, then $\mathscr F$ becomes the system of
 fibers of the fibration ${\zap q} : F\to P$,
and  Lemma \ref{folia} tells us, moreover, that  ${\zap q} : F\to P$
 is a trivial $2$-sphere bundle over $P\approx \RP^3$. We can thus 
 give the disjoint union 
 $$Z= U ~{\textstyle \coprod}~ P$$
 the structure of  
 a compact topological $6$-manifold by endowing it with the quotient topology
  induced by the map 
 $$\Psi : {\mathcal Z}_+\to Z,$$
 where 
 the restriction of $\Psi$ to 
 $\mbox{Int } {\mathcal Z}_+=U$ is the identity map $U\to U$,
 and where the restriction of $\Psi$ to 
 $\partial {\mathcal Z}_+=F$ is the fibration ${\zap q}: F\to P$.
 Indeed,  we may do this by  using the `polar coordinate'
 map 
 \begin{eqnarray*}
P\times S^2 \times [0,\infty ) & \longrightarrow & P\times \RR^3 \\
(p, \vec{x},t) & \longmapsto & (p,t\vec{x}) 
\end{eqnarray*}
as our  model for $\Psi$ near $\partial {\mathcal Z}_+= F$.  
 Now   if $g$ is of class $C^{k}$, then ${\zap q}: F\to P$ is of class
 $C^{k-1}$, and the diffeomorphism 
 $\Phi: F\to P \times S^2$ of Lemma \ref{folia} is
 of class $C^{k-2}$, so this picture actually endows  $Z$ with the structure
 of a $C^{k-2}$ manifold, in such a way that $\Psi$ becomes a $C^{k-2}$ map. 
 
 This said, we are now ready for one of the key constructions of this article:
  
\begin{thm} \label{buzz}
Let $(M,[g])$ be a space-time-oriented self-dual Zollfrei manifold, where $[g]$
can be represented by a $C^4$ split-signature metric $g$. 
Let $Z$ be the differentiable $6$-manifold obtained from ${\mathcal Z}_+$ by collapsing
$\partial {\mathcal Z}_+= F$ down to $P$ along the foliation $\mathscr F$. Then 
$Z$ can be made into a compact complex $3$-manifold in a unique 
way such that the quotient map $\Psi : {\mathcal Z}_+\to Z$
satisfies 
$$\Psi_* {\mathzap D}\subset  T^{0,1}Z.$$
Moreover, $\Psi$ is $C^\infty$  with respect to the associated complex-analytic atlas of
 $Z$ if  $g$ is itself assumed to be $C^\infty$.
\end{thm}
\begin{proof}
By construction, $\Psi$ is  a diffeomorphism between ${\mathcal Z}_+-F$
and $Z-P$. Since 
${\mathzap D}\oplus \overline{\mathzap D}=T_\CC{\mathcal Z}_+$
on ${\mathcal Z}_+-F$, it follows that there a unique complex structure $J$
on $Z-P$ with $T^{0,1}= \Psi_* {\mathzap D}$. Moreover, the assumption that
$g$ is $C^4$ guarantees that   $\mathzap D$ is  $C^3$.
Since $\mathzap D$ is involutive by Proposition \ref{critter}, 
 the Malgrange  version \cite{malgrange} of the Newlander-Nirenberg theorem \cite{newnir} 
implies that this almost-complex structure is integrable, in the sense that $Z$ admits 
local complex coordinates  in which $J$ becomes the standard
complex structure on $\CC^3$. Thus the crux of the 
theorem resides in understanding the behavior of $\Psi_* {\mathzap D}$ in
the vicinity of $P$. 

Now let us recall that the proof of Lemma \ref{folia} hinges on the introduction of 
a non-zero vector field $u$ on $F$ which spans $\ker {\zap p}_*$ at each point. 
By rescaling $u$ by an appropriate function, we may now assume henceforth  that 
${\zap q}_* u$ is always a unit vector with respect to, say,  the standard metric 
on $P\approx \RP^3$. With this convention, $S(TP)$ may be identified with 
the concrete $S^2$-bundle of unit vectors on $\RP^3$, and the $C^2$ diffeomorphism 
$\Phi : F\to S(TP)$ is just given by ${\zap q}_*u$. 
Now this  vector field $u$ is tangent to the boundary circles of the disk fibers of 
$\wp : {\mathcal Z}_+\to M$, and the fiber-wise complex structure $\jmath$ of these
disks then sends $u$ to some vector field $v=\jmath u$ along 
$\partial {\mathcal Z}_+=F$ which points inward at every boundary point
of ${\mathcal Z}_+$. Extend this $v$ to a $C^2$ vector field on a collar neighborhood of 
${\mathcal Z}_+$ so that we have $v\in \ker \wp_*$ at every point of the collar,
and then use the flow of $v$ to identify a slightly smaller  collar with $F\times [0, \epsilon )$.
Using $\Phi$ and $\Psi$, we may thus construct a $C^2$  diffeomorphism between
 a tubular 
neighborhood of 
$P$ and the $\epsilon$-tube around the zero section of $TP$, in such a manner that 
the restriction of $\Psi$ to our collar $F\times [0,\epsilon )\approx S(TP)\times [0,\epsilon )$
becomes the map 
\begin{eqnarray*}
S(TP)\times [0,\epsilon )  & \to & TP \\
(\vec{v} ,t)  & \mapsto & t\vec{v} 
\end{eqnarray*}
and so that our vector field $v$ becomes the radial field $\vec{v}/\|\vec{v}\|$. 
In particular, this picture gives us a specific isomorphism 
$$TZ|_P \cong TP\oplus TP,$$ where the first factor is tangent to
$P$, and where the second factor is transverse to it.
Moreover, this isomorphism has been constructed precisely so 
that  $\Psi_*(\jmath u)= J\Psi_* (u)$ 
at each point of $F=\partial {\mathcal Z}_+$, provided that we take  
$J: TP\oplus TP\to TP\oplus TP$ to be the almost complex structure given by 
$$J=
\left[
\begin{array}{cc}
0& -I\\
I& 0
\end{array}\right] ,
$$
where $I: TP\to TP$ denotes the identity map. 
Since the rank of $\Psi_*{\mathzap D}$ is just $1$ along $F$, this
choice of $J$ therefore gives us $\Psi_*{\mathzap D}= T^{0,1}(Z,J)$
along $P$, as desired; 
moreover, this is the only choice of $J$ with this property, since 
every unit element of $TP\subset TZ$ is of the form 
$\Psi_{*z}u$ for some $z\in \partial {\mathcal Z}_+$.
Thus, in conjunction with our previous discussion of $Z-P$, we
see that  there is a unique almost-complex structure $J$
on all of $Z$ such that $\Psi_*{\zap D}\subset T^{0,1}(Z,J)$.
However, it is not yet  clear that this $J$ is even continuous, 
much less integrable! 

We will remedy this  by next showing that $J$ is
actually {\em Lipschitz} continuous, relative to the $C^2$ structure
with which we have provisionally endowed $Z$. 
Of course, this is is only an issue near  $P$, since
$J$ has been constructed so as to be better than $C^1$ on $Z-P$.
It therefore suffices to show that $J$ is Lipschitz along each radial 
line segments $t\mapsto t\vec{v}$, t$\in [0,\epsilon)$ in   our tubular
neighborhood of $P$ modeled on the $\epsilon$ tube in $TP$,
provided we can also show in the process that the Lipschitz constants 
are uniformly bounded.

To this end,  let us therefore recall that  we have written down an 
explicit local basis $({\mathfrak w}_1 , {\mathfrak w}_2, \partial/\partial \overline{\zeta})$ for  
${\mathzap D}$ such that $[{\mathfrak w}_j, \partial/\partial \overline{\zeta}]=0$. Moreover,  the  
${\mathfrak w}_j$ are real along
$F=\partial {\mathcal Z}_+$, where they span the distribution of $2$-planes $E$
tangent to the  foliation $\mathscr F$ of $F$. Now, through a given point of 
${\zap q}^{-1}(y)\subset F$, 
there is a unique curve in the leaf ${\zap q}^{-1}(y)$ with parameter $t$ such that 
$d/dt={\mathfrak w}_1$. For any $C^2$ function $f$ on ${Z}$,
 we  then have
$$
\frac{d}{dt}\left[ \Psi_* (\frac{\partial}{\partial \overline{\zeta}}) f\right]= 
\frac{d}{dt}\frac{\partial}{\partial \overline{\zeta}} \Psi^* f=
{\mathfrak w}_1\frac{\partial}{\partial \overline{\zeta}}\Psi^* f =  \frac{\partial}{\partial \overline{\zeta}}{\mathfrak w}_1\Psi^* f =
  \frac{\partial}{\partial \overline{\zeta}}\left[ \Psi_*({\mathfrak w}_1)f\right] .
$$
Thus, setting $\zeta = \xi + i\eta$,  
$$\frac{d}{dt}\left[ \Psi_* (\frac{\partial}{\partial \overline{\zeta}})\right] 
=  \frac{\partial}{\partial \overline{\zeta}}\left[ \Psi_*({\mathfrak w}_1)\right] =
  \frac{i}{2} \frac{\partial}{\partial \eta}
\left[ \Psi_*({\mathfrak w}_1)\right]
$$
at any $y\in P$, 
since $\Psi_*({\mathfrak w}_1)\equiv 0$ along $F$, where $\eta=0$.
Here the 
 right-hand side should   be  interpreted as the
  invariant  derivative 
{\em at a zero} of a section of a  vector bundle on the disk 
 $D_x := \Psi [{\wp}^{-1} (x)\cap {\mathcal Z}_+]\approx \D^2$. 
On the other hand, 
$$\Psi_* \left(\frac{\partial}{\partial \overline{\zeta}}\right)\in T_y^{0,1}({Z},J)$$
for all $t$, by our previous discussion, so it follows that 
$$\left.
\frac{\partial}{\partial \eta}
\left[ \Psi_*({\mathfrak w}_1)\right]\right|_{\eta =0} \in T_y^{0,1}({Z},J).$$
The same argument, with 
 ${\mathfrak w}_1$ replaced by  ${\mathfrak w}_2$, tells us that  
$$\left.
\frac{\partial}{\partial \eta}
\left[ \Psi_*({\mathfrak w}_2)\right]\right|_{\eta =0} \in T_y^{0,1}({Z},J),$$
too. Along any $D_x$, we therefore have, 
near an arbitrary point $y\in P\cap D_x$,    three 
continuous sections of $T^{1,0}$
given by  
$$
{\mathfrak v}_j = \left\{
\begin{array}{cc}
\left[ \Psi_*({\mathfrak w}_j)\right]/\eta& \eta \neq 0\\
\frac{\partial}{\partial \eta}
\left[ \Psi_*({\mathfrak w}_j)\right]&\eta =0
\end{array}
\right.
$$
for $j=1,2$, and ${\mathfrak v}_3=\Psi_*(\partial/\partial \overline{\zeta})$.
These sections are linearly independent at every point, and so 
span $T^{1,0}_y$, because $\det (\Psi_*)$ only vanishes to second 
order at $Z$. Moreover, since $\Psi$ appears to be  $C^2$ in our coordinates, 
 these sections are all  continuously differentiable along $D_x$,
with coordinate derivatives expressible in terms of partial derivatives of $\Psi$ of 
order $\leq 2$. In particular,  $J$ is  Lipschitz along
$D_x$, with Lipschitz constant controlled by the 
partial derivatives of $\Psi$ of order $\leq 2$.
Since each radial line of our tube is contained in a disk $D_x$, 
and because a finite number of  balls with compact closure within coordinate domains suffice
 to cover the compact manifold 
$P$, it therefore follows that  the tensor field $J$ on ${Z}$ is  Lipschitz near $P$, and
hence  on all of $Z$. 

Since $J$ is $C^{0,1}$ on $Z$, and better than $C^1$ on $Z-P$,
the na\"{\i}ve coordinate  partial derivatives of the  components of $J$  on
$Z-P$  extend to $Z$ as 
locally bounded measurable functions. Integration by parts, however, shows that 
 these $L^{\infty}_{{loc}}$ functions are exactly the 
{\em distributional} partial derivatives of the  components of $J$. 
The Nijenhuis tensor 
$$N^\ell_{jk}= {J_k}^m\partial_m{J_j}^\ell-{J_j}^m\partial_m{J_k}^\ell+ {J_m}^\ell\partial_j{J_k}^m-{J_m}^\ell\partial_k{J_j}^m$$
of our almost-complex structure $J$ is therefore
well-defined in the distributional sense, and has $L^\infty_{{loc}}$ components.
Hence $N$ vanishes in the distributional sense, 
since by construction $N=0$ on a subset $Z-P$
of   full measure.
However,  Hill and Taylor  \cite{hiltay}
 have
 shown that the Newlander-Nirenberg theorem 
holds for Lipschitz almost-complex structures for which  $N =0$
in  the  distributional sense. Thus every point of $Z$
has a neighborhood on which we can find  a triple $(z^1, z^2,z^3)$ of differentiable 
complex-valued functions  with $dz^k\in \Lambda^{1,0}(Z, J)$
and $dz^1\wedge dz^2 \wedge dz^3\neq 0$. Taking these to be the complex
coordinate systems gives $Z$ the structure of a 
compact complex $3$-fold. In particular, this gives $Z$ a
specific preferred $C^\infty$  structure compatible with the $C^1$ structure we built
 by 
hand, so 
$\Psi$
remains a differentiable map even with respect to this brand new  atlas for $Z$. 

Now, if $g$ is actually $C^\infty$, we claim that $\Psi$ is actually a
$C^\infty$ map with respect to the tautological smooth structure on  ${\mathcal Z}_+$ 
and the complex atlas of $Z$. Away from $F\to P$, this is an immediate
consequence of the classical Newlander-Nirenberg theorem \cite{newnir},
 so we need merely
verify this assertion near $P$. To do this, let $(x^1,x^2,x^3)$ be any smooth system 
of local
coordinates on  a region ${\mathcal V}\subset P$, and pull them these functions back to $F$
as three smooth functions ${\zap q}^{*}x^j$ on 
${\zap q}^{-1}({\mathcal V})\subset F=\partial {\mathcal Z}_+$
which  are constant along the leaves of $\mathscr F$. These can then be 
extended  \cite{treves} into ${\mathcal Z}_+$ as smooth complex-valued 
functions ${\mathfrak z}^j$ near $\partial {\mathcal Z}_+$
such that 
$\partial {\mathfrak z}^j/\partial \overline{\zeta}$ vanishes to infinite order along $\eta =0$, and 
the  ${\mathfrak w}_k z^j$ will then also vanish to infinite order along $\eta =0$, too. 
Now the real and imaginary parts of the ${\mathfrak z}^j$ give us a differentiable coordinate
system on $Z$, and in these coordinates we have
$$T^{0,1}Z = \mbox{span}\left\{
\frac{\partial}{\partial \overline{{\mathfrak z}}^j}+ a_j^k ({\mathfrak z})\frac{\partial}{\partial {\mathfrak z}^k}
\right\}$$
where the smooth functions $a_j^k({\mathfrak z}^1,{\mathfrak z}^2,{\mathfrak z}^3)$ 
vanish to infinite order along the locus $P$ given by 
$\Im m~ {\mathfrak z}^j=0$. 
If $(z^1,z^2,z^3)$ is a system of holomorphic local coordinates on ${\mathcal U}\subset Z$,
where ${\mathcal U}\cap {\mathcal V}\neq \emptyset$, 
then   $z^j({\mathfrak z}^1,{\mathfrak z}^2,{\mathfrak z}^3)$ is therefore $C^\infty$ 
by  elliptic regularity. Since,  by construction, 
 each $\Psi^*{\mathfrak z}^k$ is a smooth function on 
${\mathcal Z}_+$,  it thus follows that the $\Psi^*z^j$ are smooth functions, 
too. Hence $\Psi$ is smooth with respect to the complex coordinate
atlas of $Z$,  and we are done. 
\end{proof}

\begin{defn}\label{chubby}
The {\em twistor space} of   a space-time-oriented $C^4$
 Zollfrei  self-dual $4$-manifold  $(M,[g])$   is the 
compact complex $3$-manifold $(Z,J)$ constructed from  $(M,[g])$ via Theorem \ref{buzz}. 
\end{defn}

\begin{defn}\label{checkers}
The twistor space of a non-space-time-orientable $C^4$ Zollfrei  self-dual $4$-manifold
$(M,[g])$ is  defined to be the twistor space  $(Z,J)$ of the space-time-oriented  double cover
$(\tilde{M},[g])$ of $M$.
\end{defn}

\section{Unmasking the Twistor Space}

Our construction of the twistor space of a self-dual
Zollfrei $4$-manifold  may seem rather technical. 
However,  the hidden motivation behind the entire  construction 
is   the  observation that when 
 $(M,[g])$ is one of our prototypical models, 
the associated  twistor space $(Z,J)$
is simply  the familiar complex projective $3$-space $\CP_3$. 
Let us now make this  explicit: 

\begin{lem}
If $(M,[g])$ is either $(S^2\times S^2, [g_0])$ or
$({\mathbb M}^{2,2},[g_0])$, then the twistor space $(Z,J)$ of 
$(M,[g])$, in the sense of Definitions \ref{chubby} and 
\ref{checkers}, is biholomorphic to $\CP_3$ in such a manner
than $P\subset Z$ becomes the standard $\RP^3\subset \CP^3$. 
\end{lem}
\begin{proof}
The relationship between Definitions \ref{chubby} and \ref{checkers}
makes it sufficient to consider the case of ${\mathbb M}^{2,2}$. 
Now this may seem to be a strange choice, because  Definition \ref{checkers}
  ostensibly instructs us to    pass up to the double cover
$S^2\times S^2\to {\mathbb M}^{2,2}$ and then blow down 
$\partial {\mathcal Z}_+
(S^2\times S^2)$ along the foliation $\mathscr F$. However, the quotient 
of ${\mathcal Z}_+(S^2\times S^2)$ by the covering map  action of $\ZZ_2$
 on  $\partial {\mathcal Z}_+
(S^2\times S^2)$ is just ${\mathcal Z}({\mathbb M}^{2,2})$.
Thus, Definition \ref{checkers} can be restated as saying that $Z$ is to be obtained from 
${\mathcal Z}({\mathbb M}^{2,2})$ by blowing down the hypersurface $F\subset {\mathcal Z}$.

In fact, there is a nice way of explicitly realizing of this blowing-down map. 
Let ${\mathbb V}\cong \RR^4$
be a real $4$-dimensional vector space, and let ${\mathbb V}_\CC\cong \CC^4$ be its
complexification. Then ${\mathbb M}^{2,2}$ can be be identified with the real Klein quadric 
$$Q_\RR = \{ [\psi ] \in \PP (\wedge^2 {\mathbb V})~|~ \psi \wedge \psi = 0\}$$
in $ \PP (\wedge^2 {\mathbb V})\cong \RP^5$ by choosing a diagonalizing basis
for the signature $(+++---)$ quadratic form $(\psi, \chi) = \phi\wedge \chi$ 
on $\wedge^2 {\mathbb V}$. For a suitable choice of orientation, 
the  $\beta$-surfaces of $({\mathbb M}^{2,2}, [g_0])$
are  exactly those projective planes $\RP^2\subset Q_\RR \subset \RP^5$
which are of the form 
$$\{ [\psi ] \in Q_\RR ~|~  v \hook  \psi = 0\}$$
for some  $[v ] \in \PP ({\mathbb V}^*)\cong \RP^3$. 
Thus $F({\mathbb M}^{2,2})$ may be concretely realized as
the flag manifold 
$$F_{2,3,4}= \{ ([\psi ] , [v]) \in Q_\RR  \times  \PP ({\mathbb V}^*) 
~|~ v \hook  \psi = 0\}\subset Q_\RR  \times  \PP ({\mathbb V}^*)$$
in such a way that ${\zap p}$ and ${\zap q}$ become the tautological 
projections  $F_{2,3,4}\to Q_\RR=Gr_{2,4}$ and $\PP ({\mathbb V}^*)=Gr_{3,4}$. 
However,  $Q_\RR$ is just a real slice of 
the complex $4$-quadric 
$$Q_\CC = \{ [\psi ] \in \PP (\wedge^2 {\mathbb V_\CC})~|~ \psi \wedge \psi = 0\},$$
so we have a canonical isomorphism $T_\CC Q_\RR =
TQ_\CC|_{Q_\RR}$. Any {\em complex} $\beta$-plane $\Pi \subset T_\CC Q_\RR$
is then tangent to a unique  {\em complex} $\beta$-surface 
$\CP_2\subset Q_\CC\subset \CP_5$ given by 
$$\{ [\psi ] \in Q_\CC ~|~  v \hook  \psi = 0\}$$
for some  $[v]\in \PP ({\mathbb V}_\CC^*)\cong \CP_3$. 
Thus ${\mathcal Z}({\mathbb M}^{2,2})$ may naturally be identified with 
the locus 
$$\{ ([\psi ] , [v]) \in Q_\RR  \times  \PP ({\mathbb V}^*_\CC )~|~  v \hook  \psi = 0\}
\subset F_{2,3,4} (\CC)$$
in such a way that the $\Psi : {\mathcal Z}\to Z$ is just becomes the tautological  projection
to $\PP ({\mathbb V}^*_\CC )\cong \CP_3$.

It remains to show that the constructed complex structure on $Z$ coincides with that
of $\CP_3$. To do this, we first recall that the distribution ${\zap D}$ is conformally invariant
by Proposition \ref{conformal}. Passing to the stereographic coordinates of equation 
(\ref{stereo}), it  thus suffices do our computations for the flat metric 
$d{\mathfrak x}_1^2 + d{\mathfrak x}_2^2 - d{\mathfrak y}_1^2 - d{\mathfrak y}_2^2$
using  the pseudo-orthonormal  frame
$$e_1= \frac{\partial}{\partial {\mathfrak x}_1},~ 
e_2= \frac{\partial}{\partial {\mathfrak x}_2},~ 
e_3= \frac{\partial}{\partial {\mathfrak y}_1},~ 
e_4= \frac{\partial}{\partial {\mathfrak y}_2}.$$
Since the
 connection forms ${\theta}_{\zap j}^{\zap k}$ vanish
 for this frame, 
  the distribution $\zap D$ is thus spanned by 
\begin{eqnarray*}
{\mathfrak w}_1 & = & (\zeta^2+1)\frac{\partial}{\partial {\mathfrak x}_1}
 -2\zeta \frac{\partial}{\partial {\mathfrak y}_1} + (\zeta^2-1)\frac{\partial}{\partial {\mathfrak y}_2} \\
{\mathfrak w}_2  & = &  (\zeta^2+1)\frac{\partial}{\partial {\mathfrak x}_2}+ 
(\zeta^2-1)\frac{\partial}{\partial {\mathfrak y}_1} + 
2\zeta  \frac{\partial}{\partial {\mathfrak y}_2}
\end{eqnarray*}
and $\partial/\partial \overline{\zeta}$. 
 But  the projection $\Psi: {\mathcal Z}\to \CP_3$ coming from the Klein quadric picture
is just given by 
\begin{eqnarray*}
z_1 & = & ({\mathfrak x}_1+ {\mathfrak y}_2)+ ({\mathfrak y}_1-{\mathfrak x}_2)\zeta \\
z_2 & = & ( {\mathfrak y}_1+{\mathfrak x}_2)+ ({\mathfrak x}_1-{\mathfrak y}_2)\zeta\\
z_3  & = & \zeta 
\end{eqnarray*}
in suitable affine coordinates $(z_1,z_2,z_3)$ for $\CP_3$.
Since ${\mathfrak w}_1$, ${\mathfrak w}_2$, and $\partial/\partial \overline{\zeta}$
all annihilate $z_1$, $z_2$ and  $z_3$, it 
 follows that the  complex structure $J$ we have constructed on $Z=\CP_3$ coincides
with the usual one on an open dense set, and hence everywhere. 
Thus, 
for both $({\mathbb M}^{2,2}, [g_0])$
and $(S^2\times S^2, [g_0])$,
 the twistor space is just $\CP_3$, with its standard complex structure.
\end{proof}

Now recall that the complex structure of 
$\CP_3$ is rigid,  in the sense of  Kodaira and Spencer 
\cite{KS}. In other words, because $H^1 (\CP_3 , {\mathcal O} (T^{1,0} \CP_3 ))=0$,
any complex-analytic family of deformations of the complex structure is trivial 
for small values of the perturbation parameter. It might therefore 
seem reasonable to expect that the twistor space of any Zollfrei self-dual
$4$-manifold, in the sense of Definitions \ref{chubby} and \ref{checkers},
will {\em always} turn out simply to be $\CP_3$,  with its usual complex structure. 
Our goal in this section will be to show that this is indeed
the case provided that  suitable extra hypotheses are imposed. 
To this end, we will use a beautiful  circle of   characterizations of  the  standard complex structure on 
$\CP_3$ due to Nakamura \cite{nakamura}. One such result is the following: 

\begin{thm}[Nakamura]
Let $(Z,J)$ be a compact complex $3$-manifold homeomorphic to $\CP_3$. 
If  $H^q(Z, {\mathcal O})= 0$ for all $q > 0$,  and if 
$h^0(Z, {\mathcal O}(K^{-m}))\geq 2$ for some $m> 0$,
then $(Z,J)$ is biholomorphic to $\CP_3$. 
\end{thm}

Nakamura then used this to show that any Moishezon $3$-fold homeomorphic
to $\CP_3$ must be biholomorphic to $\CP_3$ unless it is of general type. 
Recall that a compact complex $n$-fold $Z$ is said to be {\em Moishezon}  
 if there exist $n$  meromorphic 
functions $f_1, \ldots, f_n : Z\dashrightarrow \CC$ 
which give local complex coordinates near some point 
$z\in Z$; this is holds, in particular, if \cite{ueno} there is some holomorphic 
line bundle $L\to Z$ with $h^0 (Z, {\mathcal O} (L^m )) > c m^n$ for some $c> 0$ and 
all $m \gg 0$. Koll\'ar \cite{kollar} eventually  improved Nakamura's result by excluding the 
possibility that $Z$ might be of general type. Thus: 

\begin{thm}[Nakamura/Koll\'ar]  \label{nakol}
A  Moishezon manifold
  is homeomorphic to $\CP_3$ iff it is  
 biholomorphic to $\CP_3$. 
\end{thm}
 
 The following standard piece of folklore  is  a minor variation on one of   Nakamura's
results \cite{nakamura}. We  include a proof here only because one does not seem to appear
  elsewhere
 in the literature.
  
\begin{cor}\label{bigdef}
Let $J_t$ be a  family of smooth, integrable almost-complex structures
on a smooth compact $6$-manifold $Z$,  which, in the $C^\infty$ topology, 
depends continuously on an auxiliary  real variable 
$t\in [0,1]$. If $(Z, J_0)$ is biholomorphic to the standard $\CP_3$, so is 
$(Z, J_1)$.
\end{cor} 
\begin{proof}
Kuranishi  \cite{kuranishi} has shown that whenever two smooth complex structures are 
close enough  in a sufficiently  high Sobolev norm, 
they can be  joined by a complex-analytic family in the
sense of Kodaira-Spencer. Hence there is a finite subset $\{ t_0=0, t_1, \ldots , t_\ell = 1\}$
of $[0,1]$
such that, for each $j=1, \ldots , \ell$,   $(Z, J_{t_{j-1}})$ and $(Z, J_{t_{j}})$ both 
occur as fibers of a single holomorphic 
 family of complex manifolds over over the  unit disk  $\subset \CC$.
 
Now Kodaira-Spencer theory \cite{KS} tells us that 
if $(Z, J_{t_{j-1}})$ is biholomorphic to $\CP_3$,  every nearby fiber is,
too. Hence there is a non-empty open set
in the disk for which every  correspondiing fiber satisfies $h^0 ({\mathcal O}(K^{-m}))
> m^3$ for all $m> 0$. But, by the semi-continuity principle \cite{bast}, 
the set
of parameter values for which $h^0 ({\mathcal O}(K^{-m}))
> m^3$ for a particular $m$ must be closed in the analytic Zariski topology --- i.e. 
either discrete, or the whole disk. Hence every fiber must have 
$h^0 ({\mathcal O}(K^{-m}))> m^3$ for all $m >0$, and this conclusion applies,
in particular, to 
$(Z, J_{t_{j}})$. 
Hence $(Z, J_{t_{j}})$ is   Moishezon.  Theorem 
\ref{nakol} therefore shows  that $$(Z, J_{t_{j-1}})\cong \CP_3 ~~
\Longrightarrow ~~(Z, J_{t_{j}})\cong \CP_3.$$ 

Since  $(Z,J_0)$ is biholomorphic to $\CP_3$ by hypothesis, 
it therefore follows by 
 induction on $j$  that $(Z, J_1)$ is also 
biholomorphic to $\CP_3$, as claimed. 
\end{proof}

Note  that an analogous rigidity assertion  also holds for  any $\CP_n$, 
even if $n$ is large, 
as a consequence of    an entirely different
circle of ideas due to  Siu \cite{siudef}.

Now the proof of Theorem \ref{critter} shows that 
two   self-dual Zollfrei metrics which are close in the $C^\infty$ topology will
give rise to two complex structures on $Z$ which are close
in the  $C^\infty$ topology. If $g_t$ is a continuous curve in the space of  
of $C^\infty$ self-dual Zollfrei metrics, with the $C^\infty$ topology,  
Corollary \ref{bigdef} then immediately implies
that if one of the relevant twistor spaces  is biholomorphic to $\CP_3$, so
are all the others. When this happens, the smooth submanifold $P=\Psi (F)$
 thus becomes a smoothly embedded totally real submanifold of $\CP_3$,
 and every fiber of ${\mathcal Z}_+\to M$ is then sent by $\Psi$ to 
 an embedded  holomorphic disk in $\CP_3$ with boundary on $P$. Thus: 

\begin{thm}\label{zorro} 
Let ${\mathcal C}$ be the space of $C^\infty$ self-dual Zollfrei conformal classes metrics on 
$S^2\times S^2$, endowed with the smooth topology. Let
${\mathcal C}_0\subset {\mathcal C}$ be the path component containing
our protypical example 
$[g_0]$. Then, for each  conformal class $[g]\in {\mathcal C}_0$,
the corresponding  twistor space $(Z,J)$ is biholomorphically
equivalent to $\CP_3$, equipped with its standard complex structure. 
In particular, every conformal class in ${\mathcal C}_0$ gives rise to a
smooth totally real submanifold  $P\approx \RP^3$ of $\CP^3$ and
a $4$-parameter family of embedded  holomorphic disks $(D^2, \partial D^2)\hookrightarrow
(\CP_3, P)$. 
\end{thm}

Unfortunately, however, we cannot {\em a priori}  expect 
an indefinite self-dual metric to be highly differentiable, as 
 the relevant partial differential equation is 
ultra-hyperbolic rather than  elliptic. It thus behooves us to 
see what we can say about solutions with comparatively little 
regularity. However, even trying to understand $C^4$ self-dual metrics 
will lead us to consider families of  twistor spaces with 
 so little regularity that the results of Kodaira-Spencer and 
Kuranishi cannot be invoked with confidence. Fortunately, however, 
Nakamura's results are more than enough to deal with the matter at hand:

\begin{thm}
Let $g_0$ be the standard indefinite product metric on $S^2 \times S^2$.
Then $g$ has a neighborhood $\mathscr U$ in space of $C^4$ pseudo-Riemannian metrics
such that any self-dual metric $g\in {\mathscr U}$ is Zollfrei and has
twistor space $(Z,J)$ biholomorphic to $\CP_3$. 
\end{thm}
\begin{proof}
By Theorem A, there is a $C^2$ neighborhood of $g_0$ in which every 
self-dual $g$ is Zollfrei, and if $g$ is also assumed to be 
$C^4$ close to $g_0$, then  the proof of Theorem \ref{critter}
 shows that there is a diffeomorphism between the twistor spaces
of $g$ and $g_0$ such that the almost-complex structure $J$ associated with 
$g$ is close to the  almost-complex structure $J_0$ associated
with $g_0$ in the 
$C^{0,1}$ topology on tensor fields on $Z$. Choose a biholomorphism, once and for all, between 
$(Z,J_0)$ and $\CP_3$. Then, by shrinking our neighborhood $\mathscr U$ if necessary,
we may identify the $(p,q)$-forms for $J$ with those of $J_0$ via the tautological 
projections, and it therefore makes sense to think of the operators $D$ and $D_0$
given by  $\overline{\partial}+\overline{\partial}^*$
associated to these two complex structures as being defined on the same spaces, 
even after twisting with any power of the canonical line bundle. 
Thus, for example, if we consider $D$ and $D_0$ applied to $(0,1)$-forms, 
then, for every $\varepsilon > 0$ there exists a  $\mathscr U$ such that for
every $g\in {\mathscr U}$ we have 
 $\|(D-D_0)f\|^2 \leq \varepsilon (\|\nabla f\|^2+ \|f\|^2)$
 for each and every smooth  $(0,1)$-form $f$,
where $\|~\|$ denotes the $L^2$ norm on $Z=\CP_3$ with respect to, say, 
the Fubini-Study  metric.
Now assume that  such an   elliptic operator $D_0$ has trivial kernel.
By G{\aa}rding's inequality for  $D_0$ we therefore 
have 
$$\| (D-D_0)f\|^2_{L^2}\leq \varepsilon \|f\|^2_{L^2_1}\leq 
C\varepsilon \| D_0f\|^2_{L^2}$$
so that 
$$\| Df\|_{L^2}\geq (1-\sqrt{C\varepsilon }) \|D_0f\|,$$
and we therefore $D$ has trivial kernel, too, provided that we take $\varepsilon < 1/C$. 
Thus, by shrinking our neighborhood $\mathscr U$ if necessary, we may
 arrange that every associated twistor space has 
$H^1(Z,{\mathcal O})=0$, just like $\CP_3$. Similarly, we may arrange that
$H^q(Z,{\mathcal O})=0$
 and $H^q(Z, {\mathcal O}(K^{-1}))=0$ for 
$q =1,2,3$ by further shrinking $\mathscr U$. Since $Z$ also has the same Chern classes
as $\CP_3$, the index theorem then gives us $h^0(Z, {\mathcal O}(K^{-1}))={7\choose 4}$,
so Nakamura's result certainly  guarantees that  there is a biholomorphism
between $Z$ and $\CP_3$. 
\end{proof}

The holomorphic rigidity of the twistor space implies the following geometric rigidity result: 

\begin{thm}
Let $g_0$ be the standard conformally flat split-signature metric on ${\mathbb M}^{2,2}=
(S^2\times S^2)/\ZZ_2$.
Then, in the  $C^4$ topology on the space of pseudo-Riemannian metrics,
 $g_0$ has a neighborhood $\mathscr U$ 
such that any  other self-dual metric  $g\in {\mathscr U}$ is 
of the form ${\zap f}\phi^*g_0$ for some diffeomorphism $\phi: {\mathbb M}^{2,2}\to {\mathbb M}^{2,2}$ and some function ${\zap f}\neq 0$. 
\end{thm}
\begin{proof}
If  $\mathscr U$ is small enough, every  self-dual $g\in  \mathscr U$ is Zollfrei and 
has a twistor space $(Z,J)$ which is biholomorphic to $\CP_3$ by the
previous result. This twistor space can be obtained by blowing 
$\mathcal Z$ down along $F$. Complex conjugation
in ${\mathcal Z}$ therefore induces an anti-holomorphic involution 
$\varrho : Z\to Z$ with fixed point set $P\approx \RP^3$. 
 By a change of homogeneous coordinates,
any  such  $\varrho$  can be put into  the standard form 
$$[z_0: z_1: z_2 : z_3 ]\mapsto [\overline{z} _0: \overline{z} _1: \overline{z} _2 : \overline{z} _3 ],$$
as may be seen by considering the induced action on the sections of the hyperplane line bundle,
thought of as meromorphic functions with simple poles along an invariant hyperplane.  
Thus $P$ becomes the standard $\RP^3\subset \CP_3$ in these coordinates. 
Let $Q$ denote the quadric given in these coordinates by $z_1^2+z_2^2+z_3^2+z_4^2=0$,
and observe that 
 $[Q]$ now generates $H^2(\CP_3-P, \ZZ)$. However, 
 any fiber disk of ${\mathcal Z}_+\to \tilde{M}$ generates 
$H_2({\mathcal Z}_+, \partial {\mathcal Z}_+; \ZZ)$,
where $\tilde{M}=S^2\times S^2$ is the 
 space-time-oriented double cover  of $M={\mathbb M}^{2,2}$. Since $\Psi$ induces a 
homotopy equivalence between $\CP_3-P$ and ${\mathcal Z}_+$, Poincar\'e duality
now tells us that each of these holomorphic disks must meet $Q$ 
in exactly one point.  Thus $\Psi^{-1}(Q)$ is a section of $(\mbox{Int }{\mathcal Z}_+)\to 
\tilde{M}$. Moreover, the non-trivial deck transformation $\tilde{M}\to \tilde{M}$
acts on $Q$ via the complex conjugation map $\varrho$, so we have
constructed a diffeomorphism $\phi : (Q/\varrho )\to M$, and 
since $Q$ is a complex submanifold of $\mbox{Int }{\mathcal Z}_+$,
our construction of ${\zap D}=T^{0,1}(\mbox{Int }{\mathcal Z}_+)$ also shows that
$\phi$ is of class $C^{k,\alpha}$ if $g$ is of class $C^{k,\alpha}$. But the two holomorphic
disks that make up $C_x=\Psi [\wp^{-1}(x)]\subset \CP_3$ have the same boundary along 
$P= \RP^3$, and their union is therefore a rational curve in $\CP_3$,
for any $x\in M$. Each such curve meets $Q$ in a conjugate pair of points; and since
$Q\subset \CP_3$ has degree $2$, this means that $C_x$ has degree $1$. 
Hence each $C_x$ is a projective line $\CP_1\subset \CP_3$. 
However, $P=\RP^3$ is the space of $\beta$-surfaces of 
$(M,[g])$, and, for any $x\in M$, ${\zap q}[{\zap p}^{-1}(x)]= C_x\cap \RP^3$.
Thus any $\beta$ surface in $M=Q/{\varrho}$ is obtained by choosing some point 
$y\in \RP^3$, looking at all the $\RP^2$-family of all $\varrho$-invariant projective lines 
in $\CP_3$ that pass through 
$y$, and  tracing out the intersections of these lines with $Q$. 
But this same picture also, in particular, describes the 
$\beta$-surfaces of $g_0$. 
We have thus found a diffeomorphism $\phi$ between 
$M$ and $Q/\varrho = (S^2\times S^2)/\ZZ_2= {\mathbb M}^{2,2}$
which sends $\beta$ surfaces to $\beta$-surfaces. 
Since  this last statement means that $\phi$ takes null vectors to null 
vectors, we have $\phi^*[g_0]= [g]$, and hence
$g= {\zap f}\phi^*g_0$, as promised. 
\end{proof}

It will turn out that the situation on $S^2\times S^2$ is far different. Nonetheless, 
we do get some interesting  immediate  geometric pay-off from the present
discussion: 
  
  \begin{thm} Let $g_0$ be be the standard indefinite product 
 metric on 
$S^2\times S^2=\CP_1\times \CP_1$.
Then $g_0$ has a neighborhood $\mathscr U$ in the  space of $C^4$ pseudo-Riemannian metrics
such that any self-dual metric $g\in {\mathscr U}$ is of the form 
 $g=\psi^*h$, where $h$ is  an indefinite
  {\em  Hermitian} metric on $\CP_1\times \CP_1$, and where 
   $\psi$ is  a self-diffeomorphism of $S^2\times S^2$. 
 \end{thm}
 \begin{proof}
 The quadric $Q\subset \CP_3$ given by $z_1^2+z_2^2+z_3^2+z_4^2=0$ 
  does not meet the standard $\RP^3$. 
 For every self-dual metric $g$ close to $g_0$ in the $C^4$ topology, 
 $P$ will be $C^1$ close to the standard $\RP^3$,
 and so will also not meet $Q$ if our neighborhood $\mathscr U$
 is small enough. The inverse image of $Q$ under $\Psi: {\mathcal Z}_+\to M$ 
 is therefore a complex submanifold of $\mbox{Int }{\mathcal Z}_+$. 
 Moreover, the fibers of ${\mathcal Z}_+$ have intersection number $1$
 with $Q$, and as both $Q$ and these disks are complex submanifolds, 
 it follows that each fiber meets $Q$ transversely in one point. Thus
 $Q$ is the image of a smooth section $\mathfrak J$ of $\mbox{Int }{\mathcal Z}_+$.
 But this section is a bihomorphism between $(M, {\mathfrak J})$ and $Q\cong \CP_1\times \CP_1$;
 in particular, $\mathfrak J$ is integrable. On the other hand,  $\mathfrak J$ is,
 by construction, a $g$-compatible almost-complex structure. Thus
 what we have constructed is a diffeomorphism $\psi: \CP_1\times \CP_1 \to M$ such that
 $\psi^*g$ is an indefinite  Hermitian metric.
  \end{proof}

Finally, we observe that the smooth topology of the twistor space is
always standard, even without restrictions on our  Zollfrei self-dual $4$-manifold. 
This will turn out to be  quite useful in \S \ref{kahler} below.

\begin{thm} \label{difftwist}
Let $(M,[g])$ be a self-dual Zollfrei $4$-manifold, and let $Z$ be the twistor space 
of $(M,[g])$, as defined in Definitions \ref{chubby} and  \ref{checkers}. 
Then $Z$ is diffeomorphic to 
$\CP_3$ in such a manner  that  the Chern classes $c_j(Z,J)$ are sent  to the usual 
Chern classes of  $\CP_3$.
\end{thm}
\begin{proof}
By passing to a double cover if necessary, we may assume that 
$M$ is space-time orientable.
Thus $M$ is  homeomorphic 
to $S^2\times S^2$, by Theorem \ref{s2xs2}. 
Let $Y\subset Z$ be the closure of a small tubular neighborhood of $P\approx \RP^3$, and let 
$X=Z-(\mbox{Int } Y)$. Thus   $Y\approx \RP^3 \times D^3$, 
$X\cap Y\approx \RP^3 \times S^2$, and 
$X\approx {\mathcal Z}_+$.

Next,  choose an almost-complex structure ${\mathfrak J}$ on $M$ which is compatible
with $g$ and the space-time orientation. Then ${\mathcal Z}_+$
is diffeomorphic to the unit disk bundle in  the anti-canonical line bundle 
 $\Lambda^{0,2} (M ,{\mathfrak J} )$.
 In particular, ${\mathcal Z}_+$ deform retracts to a copy of  $M$.
 Moreover, $T{\mathcal Z}_+|_M= TM\oplus \nu$, where the normal bundle
 $\nu$ of $M$ is exactly the anti-canonical line bundle. Since
 we therefore have 
 $c_1(\nu)=  c_1(M, {\mathfrak J})$, so $w_2 (T{\mathcal Z}_+)|_M= 2w_2(TM)= 0$.
 It follows that $X$ is spin. 

Now  $X$ is simply connected, and since the inclusion $X\cap Y\hookrightarrow
Y$ induces an isomorphism of fundamental groups, the Seifert-van Kampen 
theorem tells us that $Z$ is simply connected, too. Since the inclusion 
$\partial X \hookrightarrow X$ is homotopic to  an $S^1$-bundle projection 
$\RP^3 \times S^2 \to M$, the 
Mayer-Vietoris sequence of $X\cup Y$ now becomes 
$$\begin{array}{rrr}
   &  \cdots ~\to  & H^1(\RP^3 \times S^2) ~\to  \\
     H^2 (Z) ~\to &H^2 (\RP^3) \oplus H^2 (S^2 \times S^2) ~\to & H^2 (\RP^3\times S^2) ~\to \\
        H^3 (Z) ~\to & H^3 (\RP^3) \oplus H^3 (S^2 \times S^2) ~\to  &H^3 (\RP^3\times S^2) ~\cdots
\end{array}$$
and so tells us  that   $H^2(Z, \ZZ) = \ZZ$ and $H^3 (Z, \ZZ ) =0$. In the same way, 
we also see that the inclusions $X\hookrightarrow Z$ and
$Y\hookrightarrow Z$ induce an injection 
 $$H^2(Z, \ZZ_2) \hookrightarrow H^2 (X, \ZZ_2) \oplus H^2 (Y, \ZZ_2),$$
 so the fact that $X$ and $Y$ are both spin implies that that $Z$ is spin, too. 
 
 Now a theorem of C.T.C. Wall \cite{wall6} asserts the diffeotype of a simply connected 
compact spin $6$-manifold with torsion-free $H^2$ and $H^3$  is completely determined by 
the ranks of these groups, 
the Pontrjagin class $p_1(TX)$, and 
the   trilinear form 
$$\smile : H^2(X, \ZZ ) \times H^2(X, \ZZ ) \times 
H^2(X, \ZZ ) \to \ZZ.$$
To finish the proof, it thus just remains  to check that  $Z$ and $\CP_3$
have the same Pontrjagin class and trilinear form. 

To this end, notice that, since  $M$ is homeomorphic to $S^2 \times S^2$, 
 our almost-complex structure ${\mathfrak J}$ must have 
\begin{eqnarray*}
c_1&\equiv& w_2=0  \bmod 2 ,\\
c_1^2 & = & 2\chi+3\tau = 8 ,  
\end{eqnarray*}
and we must therefore have 
$c_1(M,{\mathfrak J})=(2,2)\in \ZZ\oplus \ZZ= H^2(M,\ZZ)$ after correctly orienting each factor $S^2$
 of $S^2 \times S^2$. Since $c_1(\nu)= c_1(M,{\mathfrak J})$, 
 the Poincar\'e dual of $M\subset Z$ has evaluation
$2$ on a factor $S^2$,  and since the above Mayer-Vietoris sequence
shows that this evaluation map $H^2(Z,\ZZ)\to H^2(S^2,\ZZ )$ is an isomorphism, 
  it follows that $[M]= 2\alpha$
for a generator  $\alpha\in H^2 (Z, \ZZ)\cong \ZZ$.
But since    $c_1(\nu) = c_1 (M,{\mathfrak J})= (2,2)$, it follows that   
$(2\alpha)^3= [M]^3= (2,2)\cdot (2,2)=8$, so that  $\alpha^3=1$.
This shows that  $Z$ has the same trilinear form as $\CP_3$. 

Now  notice that  $p_1(TZ|_M) = p_1 (TM)+ p_1(\nu )$. However, since 
$M$ has an orientation-reversing homeomorphism, it has
vanishing signature, and we therefore have  $p_1 (TM)=0$ by the
Hirzebruch signature theorem \cite{milnorstaf}. 
Thus $ p_1 (TZ)\cdot (2\alpha )= \langle p_1 (TZ), [M]\rangle = [c_1 (\nu) ]^2= 8$,
and hence $p_1(TZ)=4\alpha^2$.  Since this is the same answer  one obtains for $\CP_3$, 
Wall's theorem now allows us to 
conclude that $Z\approx \CP_3$.  
Moreover, this diffeomorphism can be chosen so that the pull-back of the hyper-plane class 
in $H^2(\CP_3, \ZZ)$ is $\alpha \in H^2(Z, \ZZ)$. Since we have also shown 
that $c_1(Z,J)= 4\alpha$,  this diffeomorphism also takes 
 the Chern classes of 
$(Z,J)$  to those of the usual complex structure on $\CP_3$, as promised. 
\end{proof}

\section{Families of Holomorphic Disks}

In this section, we will show that every small perturbation of the standard embedding 
$\RP^3\hookrightarrow \CP_3$ gives rise to a self-dual Zollfrei conformal structure
on $S^2\times S^2$.

First let us recall that there is a standard  $(S^2\times S^2)$-family of holomorphic disks
in $\CP_3$ with boundaries on the standard $\RP^3\subset \CP_3$. Indeed, the 
boundary circles of these disks are exactly the real projective lines $\RP^1\subset
\RP^3$. Each such real projective line is contained in a unique complex
projective line $\CP_1\subset \CP_3$, and divides it into two hemispheres. 
A choice of orientation for such an $\RP^1$ then uniquely determines 
a hemisphere for which it 
 is the oriented boundary. These hemispheres  are the 
promised  holomorphic disks. 
 
 A complex projective line $\CP_1\subset \CP_3$ is the complexification
 of a real projective line $\RP^1\subset \RP^3$ iff it ia 
 $\varrho$-invariant, where 
 $\varrho: \CP_3\to \CP_3$
 denotes the complex-conjugation map
 $$\varrho ([z_1 : z_2 : z_3 : z_4]) = [\bar{z}_1 : \bar{z}_2 : \bar{z}_3 : \bar{z}_4].$$
 Now, for reasons of degree, 
   every $\varrho$-invariant $\CP_1\subset \CP_3$
 must meet  the standard quadric 
$${\mathcal Q} = \Big\{ [z_1 : z_2 : z_3 : z_4] \in \CP_3~\Big|~ z_1^2+z_2^2+z_3^2+z_4^2=0\Big\}$$
  in a conjugate pair of points; and exactly one of these points 
 will lie in each of the hemispheres into which the  $\CP_1$ is divided by  
 the fixed-point set  $\RP^3$ of $\varrho$. 
 Conversely, each point 
 $z\in{\mathcal Q}$ is joined to its conjugate point $\varrho (z)$
 by a unique $\varrho$-invariant $\CP_1$, and so is contained in exactly
 one of such hemisphere. 
Thus,  the parameter space of our  family may conveniently be identified with  
${\mathcal Q}\approx S^2\times S^2$. Moreover,  the standard conformal
structure on $S^2\times S^2$ is completely encoded by this picture, in the 
sense that each $\beta$-surface is precisely 
the family of disks whose boundaries  pass through some given point $y\in \RP^3$. 

Although this entire story   takes place in projective space, 
each of the individual disks in question actually lies in an affine subset.
To see this, we once again let $[z_1 : z_2 : z_3:z_4 ]$ be the 
 standard homogeneous coordinates on $\CP_3$,
so that standard $\RP^3\subset \CP_3$
is  represented by $z_1,\ldots , z_4$ real, 
and consider the affine chart  $(\zz_1 , \zz_2 , \zz_3)$ on $\CP_3$ defined by 
$$
\zz_1  =  \frac{z_1-iz_2}{z_1+iz_2} ~, \hspace{0.5in}
\zz_2  =  \frac{z_3}{z_1+iz_2}~, \hspace{0.5in} \zz_3 =  \frac{z_4}{z_1+iz_2} ~. 
$$
This chart realizes  the complement of the line 
$z_1=z_2=0$ 
in $\RP^3$ as the totally real submanifold 
$B$ of  $\CC^3$ 
given by \begin{equation}
\label{affine}
\zz_1\overline{\zz}_1  =  1  ~, \hspace{0.5in}
 \zz_1\overline{\zz}_2 = \zz_2  ~, \hspace{0.5in}  \zz_1\overline{\zz}_3 = \zz_3 ~.
\end{equation}
For each $a, b\in \CC$,   the disk
$$
|\zz_1|  \leq  1~, \hspace{0.5in} \zz_2  =  a + \bar{a} \zz_1 ~,\hspace{0.5in} 
\zz_3  =  b+\bar{b}\zz_1 
$$
has boundary on $B$, and 
belongs to the family under discussion. 
Notice that, as promised, these unparameterized disks 
depend on $4$ real parameters. 
Of course, each of these may in turn be
realized as a {parameterized} holomorphic disk
in a $3$-parameter family of ways by also setting 
$$
\zz_1=  \frac{c\zeta  + d}{\overline{c}+ \overline{d}\zeta} ~, ~~~~|\zeta|\leq 1~, ~~~|c|^2-|d|^2 =1.
$$
 In this manner, we actually obtain a 
 $7$-parameter family of
{\em parameterized} disks. 
In any case, it will  suffice for our purposes to primarily focus on  the
particular 
parameterized disk
$$\zz_1=\zeta ~,  \hspace{0.5in}
|\zeta |\leq 1~,
 \hspace{0.5in}
\zz_2=\zz_3=0 ~,
$$ 
since all the other disks in the family can be obtained from this one via the action of 
$PSL(4,\RR)$ on $\CP_3$.

We will now appeal to some general results concerning holomorphic 
disks in $\CC^n$ with boundary on a totally real submanifold. 
Suppose that  $X^n\subset \CC^n$ is a maximal totally real differentiable submanifold, 
in the sense that $T\CC^n|_X=TX\oplus J(TX)$. The first result we will need
is a regularity result  \cite{chirka}: 

\begin{lem}[Chirka] \label{regularity}
Suppose that $X^n \subset \CC^n$ is a totally real submanifold
of class $C^{\ell+1}$, $\ell\geq 2$, and that $\digamma: (D, \partial D)\to (\CC^n , X)$
is a $C^1$-map which is holomorphic in the interior of the disk. Then 
$\digamma$ is actually a $C^{\ell}$ map. 
\end{lem}

Now suppose that $X$ is a
 is a maximal  totally real submanifold of $\CC^n$,
 and that $\digamma: (D, \partial D)\to (\CC^n , X)$ is a holomorphic disk with
 boundary on $X$. Then $\digamma$ is said to have {partial indices}  
 $\kappa_1, \ldots , \kappa_n$ if there is a map
 $A: D \to GL (n , \CC )$ which is  holomorphic on the interior 
 of $D$ and continuous up to the boundary  
 such that $TX|_{\digamma(\zeta)}\subset \CC^n$ is the real span of the 
 columns of the matrix 
 $$
 A(\zeta ) \left[\begin{array}{ccc}\zeta^{\kappa_1/2} & 0 & 0 \\0 & \ddots & 0 \\0 & 0 & 
 \zeta^{\kappa_n/2}\end{array}\right]
 $$
 for all $\zeta\in \partial D$. 
 These partial indices turn out to be well defined up to permutation. Their sum
$$\kappa = \kappa_1+ \cdots + \kappa_n$$
is called the {\em Maslov index} of the holomorphic disk $\digamma$.  
An application of the Banach-space implicit function theorem 
to the Hilbert transform on the circle 
leads to the following result \cite{glob,ohrk}: 

\begin{prop}[Globevnik/Oh] \label{globo}
Suppose that $\digamma: D\to \CC^n$ is a holomorphic map of the unit disk
whose boundary is contained in a totally real submanifold $X$ of class 
$C^{2\ell+1}$. Suppose, moreover, that all the partial indices $\kappa_1,  \ldots ,  \kappa_n$
 of $\digamma$ satisfy $\kappa_j\geq -1$.
Then, for any totally real submanifold $X^\prime$ of $\CC^n$ which
is sufficiently close to $X$ in the $C^{2\ell+1}$-topology, 
 there is a $(\kappa+n)$-real-parameter family 
  of holomorphic embeddings 
$(D,\partial D) \hookrightarrow (\CC^n , X^\prime)$,
where $\kappa = \kappa_1+ \cdots + \kappa_n$
is the Maslov index of $\digamma$. 
This family is of class $C^\ell$, depends  in 
a $C^\ell$ manner on the choice of $X^\prime$, 
and 
sweeps out all holomorphic maps of the disk which satisfy the relevant boundary conditions
and which are
$C^\ell$ close to $\digamma$. 
\end{prop}

Let us now apply these ideas  to the case at hand. If we take $X$ to be the 
 submanifold $B=\RP^3-\RP^1$ of $\CC^3$  defined by (\ref{affine}), 
 and consider  the holomorphic disk $\digamma: D\to \CC^3$ given by
  $\zeta \mapsto (\zeta, 0 , 0)$ for  $|\zeta | \leq 1$, then  
 $TB$ is spanned over $\RR$ by the columns of the matrix 
$$\left[\begin{array}{ccc}i\zeta & 0 & 0 \\0 & \zeta^{1/2} & 0 \\0 & 0 & \zeta^{1/2}\end{array}\right]$$
for all $\zeta \in \partial D$. 
The partial indices of this disk are thus  $\kappa_1=2$, $\kappa_2=1$, and 
$\kappa_3=1$, and  its Maslov index is consequently $\kappa = 4$. Proposition
\ref{globo} thus asserts that the $7$-parameter  family of perturbations of $\digamma$ we previously
found by hand is actually stable under
deformations of $B$. That is, for any $B^\prime$ represented by a 
a section of the normal bundle of $B\subset \CC^3$ of small $C^3$ norm 
on a neighborhood of $f(S^1)\subset B$, we can find a $C^1$ family of parameterized 
holomorphic disks near $\digamma$ with boundary values in $B^\prime$ and
nonetheless $C^1$ close to the boundary values of a neighborhood
of $\digamma$ in our original $7$-parameter family. Provided the norm of this section is
small,  each of the new disks will remain embedded, and will meet the hypersurface
$$\zz_1+\zz_2^2+\zz_3^2=0$$
that represents the quadric $\mathcal Q$ in our affine chart.

Let us now give this assertion  a more concrete geometrical interpretation.
Suppose that 
$P\subset \CP_3$ be the image of a general  $C^\infty$ embedding of $\RP^3$ into $\CP_3$
which satisfies the sole constraint that, with respect to the $C^3$ topology on the space of maps,
 it lies  in a sufficiently small  neighborhood ${\mathcal X}$ of
 the standard embedding. By shrinking ${\mathcal X}$ if necessary, we may assume that every
 such $P$ is totally real and does not meet the quadric ${\mathcal Q}\subset \CP_3$. 
 Since the complement of $Q$ is a tubular neighborhood of $\RP^3\subset \CP_3$, 
any such  
  $P$ may be represented by   a smooth section of the normal bundle 
  of $\RP^3$; and since the complex structure $J$ provides an isomorphism
  between the tangent and normal bundles of $\RP^3$, the freedom in choosing 
  $P$ amounts to that of choosing a vector field on $\RP^3$ of small $C^3$ norm. 
  Proposition \ref{globo}, in conjunction with Lemma \ref{regularity},  now tells us the following: 
    
 \begin{prop}\label{cadre}
Suppose that  $P\subset \CP_3$ is the image of a smooth embedding
$\RP^3\hookrightarrow  \CP_3$ which is sufficiently close to  the standard one in the 
$C^{3}$ topology. Then $P$ contains a uniquely determined  smooth family 
of embedded oriented circles $\ell_x \subset P$, $x\in S^2\times S^2$, each of  which 
bounds an  embedded holomorphic 
disk $D^2\subset \CP_3$ whose relative homology class 
 generates $H_2(\CP_3, P; \ZZ) \cong \ZZ$.
The corresponding family of holomorphic disks is smooth,
and the interiors of these disks smoothly foliate $\CP_3-P$. 
 \end{prop}

In fact, 
the existence of a $C^1$ family of such holomorphic disks  simply follows
from  elementary Fourier analysis and the inverse function theorem, 
and so may be rederived by essentially repeating the self-contained arguments
  given in \cite{lmzoll}. Once this is known, one can then use 
  Lemma \ref{regularity} to conclude that each of the constructed disks is actually smooth,
and  the smoothness of the constructed family then follows from Proposition 
\ref{globo} by showing that it locally
coincides
with the  
families of disks obtained by perturbing any given smooth disk through disks of
increasing regularity.  A less elementary, but  distinctly compelling,   
road to the same conclusion would be to appeal to      
the  non-linear elliptic methods that are now  standard in  the 
 theory of $J$-holomorphic curves \cite{mcsalt}.
 
\section{Constructing Self-Dual Metrics}

So far, we  have  associated a $4$-dimensional space of embedded holomorphic
disks  with each small perturbation of $\RP^3\subset \CP_3$. 
To finish our construction, we need to show that this $4$-dimensional
parameter space  carries a natural self-dual 
split-signature conformal structure. This will be obtained via  the 
following  general mechanism:

\begin{prop}\label{machine}
 Let $M$ be a smooth connected 
$4$-manifold,  and let $\varpi: {\mathcal X}\to M$ be a smooth 
$\CP_1$-bundle. Let   $\varrho  : {\mathcal X}\to {\mathcal X}$  be an involution
which 
commutes with  the projection $\varpi$, and has  as fixed-point  set 
${\mathcal X}_\varrho$  an $S^1$-bundle over $M$ which disconnects
${\mathcal X}$ into two closed $2$-disk bundles ${\mathcal X}_\pm$
with common boundary ${\mathcal X}_\varrho$. 
Suppose that $\Dye \subset T_\CC {\mathcal X}$ is
a    distribution of complex $3$-planes on ${\mathcal X}$
such that 
\begin{itemize}
\item 
$\varrho _* \Dye = \overline{\Dye}$; 
\item   the restriction of $\Dye$ to ${\mathcal X}_+$
is smooth 
and  involutive; 
\item $\Dye \cap \overline{\Dye} = 0$ on   ${\mathcal X}-{\mathcal X}_\varrho$; 
\item $\Dye \cap \ker \varpi_*$  is the $(0,1)$ tangent space of the $\CP_1$ fibers
of $\varpi$;  and 
\item  the restriction of $\Dye$ to a fiber of ${\mathcal X}$  has 
$c_1= -4$ with respect to the complex orientation.
\end{itemize}
Then $E=\Dye\cap T{\mathcal X}_\varrho$ is an integrable distribution of real   $2$-planes on 
${\mathcal X}_\varrho$, and  $M$ admits a 
  unique smooth split-signature  self-dual conformal structure $[g]$
  for which the $\beta$-surfaces are the projections
  via $\varpi$  of the integral manifolds of $E$.   
  \end{prop}
\begin{proof}
Let us begin by noticing that, since   $\Dye= \varrho^*\overline{\Dye}$ is continuous on the
closed sets ${\mathcal X}_+$ and 
${\mathcal X}_-$, it is continuous on all of $\mathcal X$. 

Now let $V^{0,1}$ be the $(0,1)$ tangent space of the fibers. By hypothesis, 
$V^{0,1}\subset \Dye$, so that  ${\wot}= \Dye /V^{0,1}$ is a well
defined rank-$2$ complex vector bundle. Also notice that, since $\Dye \cap \ker \varpi_* = V^{0,1}$,
 the fibers of 
${\wot}$ are carried injectively into $T_\CC M$ by $\varpi_*$. 
We may therefore define a  continuous map 
\begin{eqnarray*}
\psi :  {\mathcal X} & \to&  Gr_2(T_\CC M) \\
z&\mapsto & \varpi_* ({\wot}|_z)=  \varpi_* (\Dye |_z)
\end{eqnarray*}
which makes the diagrams 
\setlength{\unitlength}{1ex}
\begin{center}\begin{picture}(80,17)(0,3)
\put(10,17){\makebox(0,0){${\mathcal X}$}}
\put(18,19){\makebox(0,0){$\psi$}}
\put(18,5){\makebox(0,0){$M$}}
\put(28,17){\makebox(0,0){$Gr_2(T_\CC M)$}}
\put(11,15.5){\vector(2,-3){6}}
\put(25,15.5){\vector(-2,-3){6}}
\put(12,17){\vector(1,0){10}}
\put(36,10){{and}}
\put(52,17){\makebox(0,0){$\mathcal X$}}
\put(71,17){\makebox(0,0){$Gr_2(T_\CC M)$}}
\put(52,5){\makebox(0,0){$\mathcal X$}}
\put(71,5){\makebox(0,0){$Gr_2(T_\CC M)$}}
\put(50,11){\makebox(0,0){$\varrho$}}
\put(73,11){\makebox(0,0){$c$}}
\put(60,6.5){\makebox(0,0){$\psi$}}
\put(60,19){\makebox(0,0){${\psi}$}}
\put(71,15.5){\vector(0,-1){9}}
\put(52,15.5){\vector(0,-1){9}}
\put(53.5,17){\vector(1,0){11}}
\put(53.5,5){\vector(1,0){11}}
\end{picture}\end{center}
commute, where $c$ denotes the map induced by complex conjugation 
$T_\CC M\to T_\CC M$.

Now let $\zeta$ be a smooth, fiber-wise holomorphic
coordinate on $\mathcal X$, and notice that the corresponding 
vertical vector field $\partial / \partial \overline{\zeta}$
is  both smooth and  a section of $\Dye$.  Next, near any point of the interior of ${\mathcal X}_+$,
let ${\mathfrak w}_1$ and ${\mathfrak w}_2$ be any two local sections of $\Dye$ which are linearly independent
from  $\partial / \partial \overline{\zeta}$ and from each other. Then 
 the involutivity hypothesis
$[C^\infty (\Dye) , C^\infty (\Dye) ]\subset
C^\infty (\Dye)$
tells us that 
$$\frac{\partial}{\partial \overline{\zeta}} \left( \varpi_* ({\mathfrak w}_j) \right)
=  \varpi_* ({\mathcal L}_{\frac{\partial}{\partial \overline{\zeta}}}{\mathfrak w}_j)
=  \varpi_* \left(\left[  \frac{\partial}{\partial \overline{\zeta}} ,{\mathfrak w}_j \right]\right) 
\equiv 0 \bmod  \Big\langle  \varpi_* ({\mathfrak w}_1), \varpi_* ({\mathfrak w}_2)\Big\rangle ,$$ 
 and it follows that $\psi$ is a  fiber-wise holomorphic on the interior ${{\mathcal X}}_+$. 
Since $\psi = c\circ \psi \circ \varrho$, it thus follows that $\psi$ is also 
fiber-wise holomorphic on the interior ${{\mathcal X}}_-$.
However,  $\psi$ is also continuous across ${\mathcal X}_\varrho = {\mathcal X}_+
\cap {\mathcal X}_-$, so this implies   that $\psi$ is actually fiber-wise holomorphic on
all of ${\mathcal X}$.

By construction,  the restriction of ${\wot}$ to $\varpi^{-1}(x)$
is the  pull-back, via $\psi$, of the universal bundle ${\mathbb U}$ over 
$Gr_2 (T_{\CC} M|_x)\cong Gr_2 (\CC^4)$. Now consider the 
Pl\"ucker embedding
\begin{eqnarray*}
{\mathfrak P}  : Gr_2 (T_{\CC} M) & \hookrightarrow &  \PP (\wedge^2 T_{\CC} M)\\
\mbox{span} (w_1, w_2)   & \mapsto & [ w_1\wedge w_2] 
\end{eqnarray*}
and the  induced map $\hat{\psi} = {\mathfrak P}  \circ \psi : {\mathcal X} \to \PP (\wedge^2 T_{\CC} M)$.
Since  ${\mathfrak P}^*\O (-1)=\wedge^2 {\mathbb U}$,
we must have $\hat{\psi}^*\O (-1)= \wedge^2 {\wot}$. 
But  $V^{0,1}$ is the $(0,1)$ tangent space of $\varpi^{-1}(x)$, and 
hence $c_1(V^{0,1})=-2$ on any fiber of $\varpi$. 
  On the other hand, 
$c_1(\Dye)=-4$ on $\varpi^{-1}(x)$, by hypothesis. Adjunction therefore
tells us that $c_1({\wot}) = -2$ on any fiber. Thus the restriction  of $\hat{\psi}$ to
any fiber is a holomorphic map of degree   $2$ from $\CP_1$ to
the $4$-quadric $Q_4\subset \CP_5$. There are only two possibilities for 
this map:  either it is  the inclusion of a non-degenerate plane conic
$Q_1$ into $Q_4$, or else it is a ramified double cover of a projective
line $\CP_1\subset Q_4$ branched at two points. 

The latter possibility, however, is excluded by our hypotheses. Indeed, 
any line $\CP_1\subset Q_4\subset \CP_5$ corresponds to
the curve
in $Gr_2 (\CC^4)$ given by the pencil of all $2$-planes contained
in a $3$-dimensional subspace of $\CC^4$ and containing some
fixed line. If the image of $\varpi^{-1}(x)$ under $\hat{\psi}$ were a
line, we would thus have  
$$
\varpi_* (\Dye_z)+ \varpi_*(\Dye_{z^\prime})=
\varpi_* ({\wot}|_z)+ \varpi_*({\wot}|_{z^\prime})\subsetneq T_{\CC}M|_x
$$
for all $z, z^\prime \in \varpi^{-1}(x)$. However, since 
$\Dye + \overline{\Dye}=T_\CC {\mathcal X}$ away from ${\mathcal X}_\varrho$, and because
$\varrho^*\Dye = \overline{\Dye}$, we actually have 
 $$\varpi_* (\Dye|_z)+ \varpi_*(\Dye|_{\varrho (z)})= \varpi_* (\Dye|_z+ \overline{\Dye}|_z)= T_{\CC}M|_x$$
 for all $z\in \varpi^{-1}(x)$ with $\varrho (z)\neq z$. This contradiction shows
 that $\hat{\psi}[\varpi^{-1}(x)]$ cannot be a line. 
 
 Thus $\hat{\psi}$ holomorphically includes each fiber of $\varpi$ into $\PP (\wedge^2T_\CC M)$ 
 as  a non-degenerate conic curve. For each $x$, this conic is cut out by a unique 
 $3$-plane  $\Lambda_{-}^{\CC}|_x\subset \wedge^2T_{\CC} M|_x$. The restriction of the wedge
 product
  \begin{eqnarray*}
 \wedge^2T_{\CC} M \times  \wedge^2T_{\CC} M& \to &  \wedge^4T_{\CC} M \\
(\varphi , \omega ) & \mapsto  & \varphi \wedge \omega 
\end{eqnarray*}
 to $\Lambda_{-}^{\CC}|_x$ is, moreover, always a non-degenerate bilinear form, since, by construction, 
 $\PP (\Lambda_{-}^{\CC})$ always meets the quadric $\omega \wedge \omega =0$
 in the  non-degenerate  conic  $\hat{\psi}[\varpi^{-1}(x)]$.
  
  Now $\hat{\psi}$ is at least smooth
 on the interior of ${\mathcal X}_+$. By taking the images under $\hat{\psi}$ of three generic 
 smooth local sections of $\varpi$ which avoid ${\mathcal X}_\varrho$, we can thus locally span
 $\Lambda_-^\CC$ by three smooth local sections of $\wedge^2 T_\CC M$. Thus 
 $\Lambda_{-}^{\CC}\subset \wedge^2T_{\CC} M$ is a smooth sub-bundle. Moreover, 
 essentially the same argument shows that $\hat{\psi}$ is smooth on all of $\mathcal X$. 
 
 Since $\psi \circ \varrho
 = c\circ \psi$, we must therefore have $\Lambda_{-}^{\CC}= \CC\otimes \Lambda_-$ for a unique,
 smooth real vector sub-bundle $\Lambda_-\subset \wedge^2TM$ on which the 
 wedge product is non-degenerate. However, the wedge product must be {\em  indefinite}
 on every fiber of $\Lambda_{-}$, since ${\mathcal X}_\varrho$ meets every fiber of
 ${\mathcal X}$. Since  
 $$O(3,3)/[O(2,1)\times O(1,2)]= SL(4,\RR) /SO(2,2)$$
 it follows that there is a unique smooth split-signature conformal metric $[g]$
 on $M$ for which $\Lambda_-$ is the bundle of anti-self-dual bi-vectors 
 for an appropriate orientation of   $M$. For each metric $g$
 in this conformal class, $\Lambda_-$ then corresponds via 
  index-lowering to the bundle $\Lambda^-\subset \wedge^2$ of 
  real anti-self-dual $2$-forms.
   
 Now consider the subset of the complex tangent bundle of  ${\mathcal X}_\varrho$ defined by 
 $$E_\CC= \Dye \cap T_\CC {\mathcal X}_\varrho.$$
 Since each fiber of $T_\CC {\mathcal X}_\varrho$ has codimension $1$ in $T_\CC{\mathcal X}$,
 and since $T_\CC {\mathcal X}_\varrho$ does not contain the $1$-dimensional subspace 
 $V^{0,1}\subset \Dye$, the subspace $\Dye$ is always in general position 
 relative to  $T_\CC X_\varrho$. Hence $E_\CC$ is a smooth distribution of complex
 $2$-planes on ${\mathcal X}_\varrho$. However, $\varrho$ acts trivially on ${\mathcal X}_\varrho$, and hence
 $\varrho_*$ acts  
 on $T_\CC {\mathcal X}_\varrho$ via the identity. 
 The assumption  that $\varrho_*\Dye = \overline{\Dye}$ therefore 
 implies that $\overline{E_\CC}=E_\CC$. Hence $E_\CC$  is the  
 complexification of a smooth distribution of real $2$-planes
 $$E= \Dye \cap T{\mathcal X}_\varrho$$
 on ${\mathcal X}_\varrho$. Since $T{\mathcal X}_\varrho$ and
 $\Dye$ are both closed under Lie brackets, it follows that 
   $E$ is  Frobenius integrable. 
Thus ${\mathcal X}_\varrho$ is foliated by $2$-manifolds tangent to $E$. 
But the inclusion $E_\CC \hookrightarrow \Dye$ induces
an canonical isomorphism $E\otimes \CC\to  {\wot}|_{{\mathcal X}_\varrho}$, whereas $\psi$
identifies ${\mathcal X}_\varrho$ with the bundle of real $\beta$-planes for
$(M,[g])$. Thus each integral manifold of $E$ is sent via $\varpi$ 
to a $\beta$-surface of $(M,[g])$, and  $[g]$ is therefore self-dual by Proposition \ref{roger}.
Moreover, $[g]$ is uniquely determined by this last prescription, since,  at each point of $M$,  
 the union of the tangent spaces of these $\beta$-surfaces
 is  precisely the null cone of $[g]$, and
  the conformal class of any indefinite metric is 
completely determined by its null cone.
  \end{proof}

\begin{thm}\label{bild}
Let $P\subset \CP_3$ be a smooth, totally real submanifold 
 which, in   the $C^3$ topology,  is 
close to the standard `linear'  $\RP^3\subset \CP_3$; and, for clarity,  
 fix a quadric $Q\subset \CP_3$ which is disjoint from $P$.
For each $x\in Q \approx S^2 \times S^2$, let $D_x\subset \CP_3$
be the unique holomorphic disk of the family constructed in Proposition \ref{cadre} which
passes through $x$. For each $y\in P$, set 
$$S_y= \{ x \in S^2 \times S^2~|~ y \in D_x\}.$$
Then there is a unique, smooth  Zollfrei self-dual split-signature
conformal structure $[g]$ on $Q\approx S^2\times S^2$ whose   $\beta$-surfaces
are exactly the  $S_y$, $y\in P$. 
\end{thm}
\begin{proof}
Let $M=Q\approx S^2 \times S^2$, and let ${\mathcal X}_+\to M$
be the  $2$-disk bundle whose fiber over $x\in M$ is the holomorphic disk 
$D^2\subset \CP_3$ of the family passing through $x$. Thus there is a 
tautological smooth map ${\mathzap F}: {\mathcal X}_+\to \CP_3$ which sends the
interior of ${\mathcal X}_+$ diffeomorphically onto $\CP_3-P$, and
which sends $\partial {\mathcal X}_+\to P$. Recalling 
that ${\mathzap F}_* : T_\CC{\mathcal X}_+\to T_\CC\CP_3$
denotes the derivative of this map, let ${\mathzap F}_*^{1,0}: T_\CC  
{\mathcal X}_+\to T^{1,0}\CP_3$ denote be the $(1,0)$-component of this derivative, and
let
$$\Dye = \ker {\mathzap F}_*^{1,0}\subset T_\CC{\mathcal X}_+$$
denote the  kernel of this component. Since ${\mathzap F}$ is $C^1$ close to the corresponding map 
for the flat model, we may assume that ${\mathzap F}_*^{0,1}$ is everywhere
of maximal rank, as in the flat case. Thus $\Dye$ is a smooth complex bundle of rank 
$3$ on all of ${\mathcal X}_+$. Now if $V^{0,1}$  is the $(0,1)$-tangent space of
the fibers of the $D^2$-bundle ${\mathcal X}_+\to M$, then $V^{0,1}\subset 
\Dye$ because ${\mathzap F}$ is fiber-wise holomorphic. 
But because the $5$-manifold $\partial {\mathcal X}_+$ is 
sent to the $3$-manifold $P$ by ${\mathzap F}$, each fiber 
$$E = \ker {\mathzap F}_*|_{\partial {\mathcal X}_+}$$
has dimension $\geq 2$, and since $(E\otimes \CC) \oplus V^{0,1}\subset \Dye$,
we conclude that $E$ is in fact a smooth distribution of real $2$-planes
on $\partial {\mathcal X}_+$.

Now let ${\mathcal X}_-$ be a second copy of ${\mathcal X}_+$, 
and define  $\Dye\to {\mathcal X}_-$ to be the push-forward of the distribution
of complex $3$-planes $\overline{\Dye}\to {\mathcal X}_+$  
via the tautological diffeomorphism ${\mathcal X}_+\to {\mathcal X}_-$.
Similarly, let $V^{0,1}\to {\mathcal X}_-$ be the distribution of 
complex lines obtained from $\overline{V^{0,1}}\to {\mathcal X}_+$.
Let 
$${\mathcal X} = {\mathcal X}_+\cup_{\partial {\mathcal X}_+} {\mathcal X}_-$$
be the double of ${\mathcal X}_+$.
Then we have a canonical projection $\varpi : {\mathcal X}\to M$ which 
makes ${\mathcal X}$ into a $\CP_1$-bundle with vertical 
$(0,1)$ tangent space $V^{0,1}$. Moreover, 
 our two definitions of $\Dye$ agree along the hypersurface 
$\partial {\mathcal X}_+ =    \partial {\mathcal X}_-$, because 
both coincide with $V^{0,1}\oplus (E\otimes C )$ along this 
locus. Moreover,  $V^{0,1}= \Dye \cap \ker \varpi_*$
on all of $\mathcal X$. 

Let $\varrho : {\mathcal X}\to {\mathcal X}$ be the map which interchanges 
${\mathcal X}_\pm$ via the tautological diffeomorphism.  This is
an involution of ${\mathcal X}$ which commutes with 
$\varpi$, and its fixed-point set ${\mathcal X}_\varrho= \partial {\mathcal X}_+$
divides ${\mathcal X}$ into two disk bundles over $M$. 
By construction, we have 
$\varrho_*\Dye = \overline{\Dye}$. Moreover, 
$\Dye$ is smooth, involutive, and satisfies $\Dye \cap \overline{\Dye}=0$
on $\mbox{Int } {\mathcal X}_+$, since 
the diffeomorphism ${\mathzap F}$ from $\mbox{Int } {\mathcal X}_+$ to $\CP_3-P$ sends
$\Dye$ to $T^{0,1}\CP_3$. 

Finally,  observe that any holomorphic disk $(D^2, S^1)\to (\CP_3, P)$ obtained  by 
restricting ${\mathzap F}: {\mathcal X}_+\to \CP_3$ to a fiber of ${\mathcal X}_+\to M$
must have Maslov index $\kappa= 4$, since each such disk is obtained by deforming
a disk with $\kappa = 4$ from our flat model, and the Maslov index invariant
under  deformations \cite{mcsalt}. This index is by definition the winding number of 
$\wedge^3TP$ in $\wedge^3 T^{1,0}T\CP_3$ along $S^1= \partial D^2$, relative to 
any trivialization of $T^{1,0}\CP_3$  over $D^2$, remembering that the 
space of real lines in $\CC$ is exactly $\RP^1\approx S^1$. But recall that  ${\mathzap F}_*^{1,0}$
is surjective, so that we can identify ${\mathzap F}^*T^{1,0}\CP_3$ with the quotient 
$T_\CC{\mathcal X}_+/\Dye$. Since $\wedge^6 T_\CC {\mathcal X}_+$ is 
 the complexification of the trivial real bundle $\wedge^6 T{\mathcal  X}_+$,
it follows by adjunction that 
this Maslov index must be {\em minus} the Maslov index of $T{\mathcal X}_\varrho/E
\subset \Dye$. However, the latter winding number is also exactly the degree of 
$\wedge^3 \Dye$ on a $\CP_1$ fiber, since, by construction, $\Dye$ is defined
on the double of the disk precisely by gluing $\Dye|_{D^2}$ to $\overline{\Dye}|_{\overline{D^2}}$ 
so as to send  $T{\mathcal X}_\varrho/E$ to itself. Thus the evaluation of 
$c_1(\Dye)$ on a fiber of $\varpi$ is exactly $-\kappa = -4$.

The above arguments show that all the hypotheses of Proposition \ref{machine} are 
satisfied. Thus $M=S^2\times S^2$ admits a unique
self-dual split-signature metric $[g]$ for which the $\beta$-surfaces are the 
projections to $M$ of the integral manifolds of $E\to {\mathcal X}_\varrho$. 
By construction, however, $E$ is precisely the vertical tangent bundle of
the smooth submersion 
$${\mathzap F}|_{\partial {\mathcal X}_+}: {\partial {\mathcal X}_+}\to P,$$
so these $\beta$-surfaces are exactly of the form
$$S_y = \varpi [{\mathzap F}^{-1}(y)]$$
for $y\in P\subset \CP_3$, which is to say that
$$S_y= \{ x\in M ~|~ y \in D_y\},$$
as promised. 

Now  ${\mathcal X}_+$ is diffeomorphic to 
the Chern class $(2,2)$ disk bundle over 
$S^2\times S^2$. Thus  ${\mathcal X}_\varrho= \partial {\mathcal X}_+$
is the  Chern class $(2,2)$  circle bundle over $S^2\times S^2$, and
so is 
diffeomorphic to $\RP_3 \times S^2$. 
Moreover, since every disk of the family represents
the generator of $H_2 (\CP_3,P)$, the  long exact homotopy sequence 
$$\cdots \to H_2 (\CP_3, P;\ZZ ) \to H_1 (P; \ZZ ) \to H_1 (\CP_3; \ZZ)\to \cdots $$
of the pair $(\CP_3 ,P)$ tells us that the boundary of each 
disk generates $\pi_1 (P) = H_1 (P)\cong \ZZ_2$. It
follows that ${\mathzap F}|_{{\mathcal X}_\varrho}$ induces a surjection 
$\pi_1 ({\mathcal X}_\varrho)\to \pi_1 (P)$. But 
${\mathzap F}|_{{\mathcal X}_\varrho}$ is a proper submersion, and therefore
a smooth fibration, so we  have the long exact homotopy
sequence 
$$\cdots \to \pi_2 (P) \to \pi_1 (S_y)\to \pi_1 ({\mathcal X}_\varrho)\to \pi_1 (P)\to 
\pi_0 (S_y) \to \cdots $$
and it follows that the compact surface $S_y$ is connected and  simply connected. 
Thus every $\beta$-surface of $(M,[g])$ is a $2$-sphere. 
Hence $[g]$ is Zollfrei by Theorem \ref{tasty}, and we are done. 
\end{proof}

{\bf Theorem \ref{zfsd}} now follows  from Theorems \ref{zorro} and  \ref{bild}.

\section{The K\"ahler Case}
\label{kahler}

The protypical example which motivated our entire investiagtion 
of a Zollfrei self-dual manifold was the 
the indefinite product metric 
$g_0= \pi_1^*h - \pi_2^*h$ on $S^2\times S^2$.
Notice, however, that this metric
may be considered as 
an indefinite K\"ahler metric on $\CP_1\times \CP_1$, 
with K\"ahler form $\omega =\pi_1^*\mu - \pi_2^*\mu$,
where $\mu$ denotes the area form of $(S^2,h)$.
Notice that, since $h$ has constant Gauss curvature $1$, the scalar curvature 
of $g_0$ is $s=s_{g_0}=\pi_1^*s_h-\pi_2^*s_h=2-2= 0$. 
A pseudo-Riemannian metric with this last property is said to be 
{\em scalar-flat}.

Now, more generally, suppose that we have a scalar-flat indefinite 
K\"ahler metric
$g$ on a compact complex surface $(M^4,{\mathfrak J})$. From the outset, we 
choose to give $M$ the usual {\em complex} orientation,  but 
we will also need to
 systematically consider the reverse-oriented version  
 $\overline{M}$ of our  manifold. To see why, observe that  the 
K\"ahler form $\omega$   of $(M,g,{\mathfrak J})$ is a closed non-degenerate
$2$-form on $M$, and so may be considered as a {\em symplectic form}.
However, such a form determines an  orientation, and in the present
case this orientation is the {\em opposite} of the complex-manifold
orientation; thus,  
$(\overline{M}, \omega)$ becomes a symplectic $4$-manifold,
oriented according to standard  symplectic conventions.

Notice that while $\omega$ is a self-dual $2$-form on $(\overline{M}, g)$,
it  is instead {\em anti-self-dual} on $(M,g)$. With this potential source of confusion  
kept clearly in focus, 
standard Riemannian folklore \cite{bes4} immediate tells us the following: 

\begin{lem} \label{gau} 
Let $(M^4,{\mathfrak J},g)$ be a complex surface with an {\em indefinite} K\"ahler metric. Then 
$(M,g)$ is self-dual iff $g$ is scalar-flat. 
\end{lem}
\begin{proof}
The curvature of any  K\"ahler manifold, indefinite or not, is necessarily
of type $(1,1)$, so that the corresponding curvature operator  ${\mathcal R}$
kills $\Lambda^{2,0}\oplus \Lambda^{0,2}$, and  amounts to a 
linear map ${\mathcal R} : \Lambda^{1,1}\to \Lambda^{1,1}$. 
Now observe that 
if $(M^4,{\mathfrak J},g)$ is an {\em indefinite} K\"ahler manifold, equipped  with the 
complex orientation, 
we have 
\begin{eqnarray*}
 \Lambda^{1,1} & = &  \CC \omega \oplus \Lambda^+\\
\Lambda^-_\CC & = &   \CC \omega\oplus \Lambda^{2,0}\oplus \Lambda^{0,2}
\end{eqnarray*}
so that  the 
$$\left[W_-+\frac{s}{12}\right] : \Lambda^-\to \Lambda^-$$ 
block of the curvature operator 
kills  $\Lambda^{2,0}\oplus \Lambda^{0,2}$, 
and sends $\CC \omega$ to itself.
Since $s/4$ is the trace of this block, and since $W_-$ is its trace-free part,
we therefore have
$$W_-= \left(\begin{array}{ccc}-\frac{s}{12} &  &  \\ & -\frac{s}{12} &  \\ &  & \frac{s}{6}\end{array}\right)$$
in an appropriate basis. Hence $W_-=0$ iff $s=0$, as claimed. 
\end{proof} 

Now the structure group of an indefinite K\"ahler surface is 
$U(1,1)$, which is a connected Lie group.  
Every indefinite K\"ahler surface therefore 
carries a canonical space-time orientation. 
As a consequence, 
Theorem \ref{s2xs2}  tells us that any Zollfrei scalar-flat
indefinite K\"ahler surface is homeomorphic to $S^2\times S^2$. 
However,  a much stronger assertion is actually true:

\begin{thm} \label{bihol} 
Suppose that $(M^4,{\mathfrak J},g)$ be a complex surface with scalar-flat {\em Zollfrei} 
indefinite K\"ahler metric. Then 
$(M,{\mathfrak J})$ is biholomorphic to $\CP_1\times \CP_1$. 
\end{thm}
\begin{proof}
Let $S$ be any  $\beta$-surface in $(M,g)$. 
At each point of $S$ the image of 
 $\Lambda^2TS$ in $\Lambda^2$ is then 
  the span of a  real, non-zero simple element of $\Lambda^-$. But, up to a positive
  constant,  the general  such simple $2$-form 
can be written uniquely as $\omega +  \phi + \overline{\phi}$, where 
$\phi\in \Lambda^{2,0}$ is any element of unit norm. It follows that 
the restriction of $\omega$ to $S$ is non-zero at every point.
Hence  $S$ 
 is a symplectic submanifold of  $(\overline{M}, \omega )$.
 Moreover,  since $S$ is orientable and Lemma \ref{beta} asserts that 
 $S$ is either a $2$-sphere or a projective plane,  we have 
 $S\approx S^2$. 
 Proposition \ref{structure}  thus tells us  that $[S]\cdot [S] =-2$
 in $M$, and hence that $[S]\cdot [S] =+2$ in the reverse-oriented
 manifold $\overline{M}$. 
  Hence $(\overline{M}, \omega )$ contains a symplectic $2$-sphere
  $S$ of positive self-intersection, and  
  a fundamental   result of McDuff \cite{mcrules} therefore tells us that   $(\overline{M}, \omega )$
must be   diffeomorphic 
 to either $S^2 \times S^2$ or to $\CP_2 \# k \overline{\CP}_2$, 
 $k\geq 0$. 
 Since $M$ is spin by Proposition \ref{s2xs2},  it therefore follows that $M$ is therefore
 diffeomorphic to $S^2\times S^2$. 
 
 In particular, this shows that $M$ is a minimal complex surfae
 which  admits  a Riemannian metric of positive 
 scalar curvature, and   Seiberg-Witten theory \cite{FM,spccs} therefore tells us 
 $(M,{\mathfrak J})$ must have  Kodaira dimension $-\infty$.
 By the Kodaira-Enriques classification \cite{bpv,GH}, our simply connected 
 complex surface $(M,{\mathfrak J})$ is therefore {\em rational}, in the sense of being  obtained   
  from $\CP_2$ by blowing up and down, and since 
   $M$ is also spin, so it  follows that 
  $(M,{\mathfrak J})$ is an even Hirzebruch surface $\PP [\O  \oplus \O (2m)]\to \CP_1$.
  However, Kamada \cite{kamada} has shown that the existence of scalar-flat  
  indefinite K\"ahler metrics is 
  obstructed for $m\neq 0$  because a generalized form of the   Futaki invariant
 \cite{fuma} is non-zero in all  these cases. 
  Hence $m=0$, and $(M,{\mathfrak J})$ must be  
  biholomorphic to $\CP_1\times\CP_1$. 
 \end{proof}
 
 \begin{remark} The above result  would certainly become false if the 
 Zollfrei hypothesis were  dropped. For example, 
 one can easily construct  scalar-flat indefinite 
 K\"ahler metrics  on the product $\Sigma \times \Sigma$ 
 of any Riemann surface $\Sigma$ with itself, just by setting $g= \pi_1^*h -\pi_2^*h$,
 where $h$ is a  metric  on $\Sigma$ of constant sectional curvature. 
 \end{remark}
 
 \begin{thm}\label{twisthol}
 Let $(M,g,{\mathfrak J})$ be a scalar-flat Zollfrei indefinite K\"ahler metric. Then its
 twistor space $(Z,J)$, in the  sense  of definition \ref{chubby}, is biholomorphic
 to $\CP_3$. Moreover, this biholomorphism determines a preferred 
 non-singular quadric $Q\subset \CP_3$ obtained by thinking of  
 ${\mathfrak J}$ as section of ${\mathcal Z}_+\to M$.  \end{thm}
 \begin{proof}
 The complex structure ${\mathfrak J}$ defines a section of $\mbox{Int }{\mathcal Z}_+$, and because
 ${\mathfrak J}$ is parallel, the image of this section is tangent to ${\zap E}\subset {\zap D}$. 
 The composition of $\Psi$ with this section 
 is therefore a   holomorphically embedding of $(M, {\mathfrak J})$ 
 into $(Z,J)$. Moreover, as we saw in the proof of 
 Theorem \ref{difftwist}, 
 the normal bundle $\nu$ of $(M, {\mathfrak J})$ has
 Chern class $c_1(\nu ) = c_1 (M, {\mathfrak J})$. Since 
  Theorem \ref{bihol} tells us that $(M, {\mathfrak J})\cong \CP_1\times \CP_1$,
  this therefore gives us a hypersurface $Q\subset Z$ 
  cut out by a section of the corresponding  divisor line bundle  $L\to Z$, 
  such that $Q\cong \CP_1\times \CP_1$ and such that 
   $L|_Q\cong {\mathcal O}(2,2)$. 
   For each integer $m$, we therefore have an  exact sequence
   $$0\to {\mathcal O}(L^{m-1}) \to {\mathcal O}(L^{m}) \to {\mathcal O}_Q(2m,2m)\to 0$$
   of sheaves on $Z$. Since 
   $$H^q (\CP_1\times \CP_1 , {\mathcal O}(2m, 2m )) = 
    0~~\forall q, m > 0, $$
   it follows that, as $m\to \infty$, 
   $h^1(Z, {\mathcal O}(L^{m}))$ is decreasing, while  $h^2(Z, {\mathcal O}(L^{m}))$
   and $h^3(Z, {\mathcal O}(L^{m}))$
   remain constant. Hence \begin{eqnarray*}
\chi (Z, {\mathcal O}(L^{m})) & = & h^0 ({\mathcal O}(L^{m})) -h^1 ({\mathcal O}(L^{m})) +
h^2 ({\mathcal O}(L^{m})) -h^3 ({\mathcal O}(L^{m}))  \\
 & = & h^0 (Z, {\mathcal O}(L^{m})) + \mbox{const  } ~~~\forall m \gg 0.
\end{eqnarray*}
   However, Theorem \ref{difftwist} tells us that 
   $Z$ is diffeomorphic to $\CP_3$ in  a manner sending the Chern classes of $Z$ to 
    the Chern classes of $\CP_3$. Since
    $c_1(L)= \frac{1}{2}c_1 (Z)$, the Hirzebruch-Riemman-Roch  theorem therefore tells us that 
    $$\chi (Z, {\mathcal O}(L^{m}))= \chi (\CP_3 , {\mathcal O}(2m))= \frac{(2m+1)(2m+2)(2m+3)}{6}. 
    $$
    Hence $h^0 (Z, {\mathcal O}(L^{m}))$ grows cubically in $m$. The 
     complex $3$-fold $(Z,J)$ is therefore  Moishezon. Since
     $Z$ is also   diffeomorphic to 
    $\CP_3$,   Theorem  \ref{nakol} therefore tells us that $(Z,J)$ is biholomorphic to 
    $\CP_3$. Moreover, $Q\subset Z$ is carried by this biholomorphism to a 
    non-singular hypersurface of degree $2$. 
    \end{proof}
 
Now,  which totally
 real submanifolds
 of $\CP_3$ correspond to scalar-flat self-dual metrics? 
 The following 
 result provides the key to the answer.

 \begin{prop} \label{constr}
 Let $(M, g, {\mathfrak J})$ be a Zollfrei indefinite scalar-flat K\"ahler manifold,
 let $Q\subset Z\cong \CP_3$ be the quadric constructed in Theorem \ref{twisthol}, 
 and let $P=\Psi (F)$ be the space of $\beta$-surfaces in $M$. Then
 there is a meromorphic $3$-form $\Omega$ on $Z$
 which is holomorphic and non-zero on $Z-Q$ and which has the property  that 
 its pull-back  to $P$ is a real $3$-form. 
 \end{prop}
\begin{proof}
Consider a pseudo-orthonormal frame $e_1 , \ldots, e_4$ 
on some region of $M=\CP_1\times \CP_1$ in which 
$$\omega = \sqrt{2}\varphi_1 = e^1\wedge e^2 -e^3\wedge e^4.$$
Since $\nabla\omega=0$, we have
$$\theta^2_1= \theta^1_2=\theta^3_1=\theta^1_3=0,$$
so the connection on $\Lambda^-$ is determined by  a single  $1$-form 
$$\theta = \theta_3^2=-\theta_2^3.$$
The distribution ${\zap D}$  
is thus spanned by 
\begin{eqnarray*}
{\mathfrak w}_1 & = & (\zeta^2+1)e_1 -2\zeta e_3 + (\zeta^2-1)e_4 
\\ & & 
-\frac{1+\zeta^2}{2}\Big[(\zeta^2+1)\theta_1 - 2\zeta\theta_3 + (\zeta^2-1)\theta_4
\Big]  \frac{\partial}{\partial \zeta}\\
{\mathfrak w}_2  & = &  (\zeta^2+1)e_2+ (\zeta^2-1)e_3 + 2\zeta  e_4
\\ & & 
-\frac{1+\zeta^2}{2}\Big[(\zeta^2+1)\theta_2 + (\zeta^2-1) \theta_3 + 2\zeta  \theta_4
\Big]  \frac{\partial}{\partial \zeta}
\end{eqnarray*}
and $\partial/\partial \overline{\zeta}$, where $\theta_j=\theta (e_j)$. 
However,  $\varphi_2+ i \varphi_3$ is a unit section of the canonical line bundle
of $(M,g, {\mathfrak J})$, and 
\begin{eqnarray}
d \varphi_1 & = & 0 \nonumber \\
d \varphi_ 2 & = & -\theta \wedge \varphi_3 \label{oye} \\
d \varphi_3  & = & \theta \wedge \varphi_2 \nonumber
\end{eqnarray}
Hence $d\theta$ is just the Ricci form $\rho$ of $(M,g, {\mathfrak J})$.
But the Ricci form of any K\"ahler manifold is of type $(1,1)$, and 
in our case $\rho \wedge \omega =0$, since $(M,g, {\mathfrak J})$
is assumed to be scalar-flat. We thus conclude that 
\begin{equation}
\label{curva}
d\theta \wedge \varphi_{\zap j}=0, ~~~ {\zap j} = 1,2, 3.
\end{equation}
We remark in passing that this is simply a special case  of a more general fact: 
 namely,  $\Lambda^-$ has self-dual curvature
on any scalar-flat self-dual $4$-manifold.

Let us now set 
\begin{eqnarray*}
\Omega&=& -\frac{[(1+\zeta^2) \varphi_1 + (1-\zeta^2) \varphi_2 -2\zeta \varphi_3]\wedge
[2d\zeta + (1+\zeta^2)\theta]}{(1+\zeta^2)^2}\\
&=&  (\varphi_1+ \cos t ~\varphi_2 + \sin t ~\varphi_3 ) \wedge (dt- \theta )
\end{eqnarray*}
where $$t= -2\tan^{-1}\zeta =i \log (1-i\zeta ) - i \log (1+i\zeta ).$$ 
The  restriction of this form to $F=\partial {\mathcal Z}_+$  is a real,  geometrically meaningful, and 
globally defined $3$-form.   Indeed,  $\varphi_1+ \cos t ~\varphi_2 + \sin t ~\varphi_3$ is the 
tautological $2$-form on $F$, thought of
as the space of those real null anti-self-dual $2$-forms for which the inner product with the 
K\"ahler form 
$\omega$ is $\sqrt{2}$; and $dt- \theta $  the principle-connection $1$-form of the 
unit canonical 
bundle of $(M,g, {\mathfrak J})$. Since $\Omega$ is the unique
analytic continuaton of $\Omega|_F$ up the fiber disks of ${\mathcal Z}_+\to M$,
this shows that $\Omega$ is   globally defined on the  ${\mathcal Z}_+-Q$,
where $Q$ is the image of the section ${\mathfrak J}$, which is represented by 
$\zeta =i$. 

Next,   notice that $\Omega$ annihilates ${\mathfrak w}_1$, 
${\mathfrak w}_2$,  and $\partial/\partial \overline{\zeta}$. 
Thus $\Omega$ is a $(3,0)$-form   on $Z-(P\cup Q)$. 
Moreover, equations (\ref{oye}) and (\ref{curva}) tell us that 
\begin{eqnarray*}
d\Omega  & = &
 d(\varphi_1+ \cos t ~\varphi_2 + \sin t ~\varphi_3 ) \wedge (dt- \theta )\\
&&\hspace{1cm}
 +  (\varphi_1+ \cos t ~\varphi_2 + \sin t ~\varphi_3 ) \wedge d(dt- \theta )
\\&=&  ( -\sin t ~dt\wedge  \varphi_2 + \cos t ~d \varphi_2 
+\cos t ~ dt\wedge \varphi_3 +\sin t ~ d \varphi_3  )\wedge
(   dt - \theta )\\
&&\hspace{1cm} -   (\varphi_1+ \cos t ~\varphi_2 - \sin t ~\varphi_3 ) \wedge d \theta 
 \\ &=&
  ( -\sin t ~dt\wedge  \varphi_2 - \cos t ~\theta \wedge \varphi_3 
+\cos t ~ dt\wedge \varphi_3 +\sin t ~ \theta \wedge \varphi_2  )\wedge
(   dt - \theta )
 \\ &=& \sin t ~dt\wedge  \varphi_2 \wedge \theta 
- \cos t ~\theta \wedge \varphi_3 \wedge dt 
-\cos t ~ dt\wedge \varphi_3 \wedge \theta 
+\sin t ~ \theta \wedge \varphi_2 \wedge dt
 \\ &=& 0
\end{eqnarray*}
so the $(3,0)$-form $\Omega$ is actually $\overline{\partial}$-closed  on 
$Z-P$, where it is therefore a meromorphic $3$-form with only  a pole of order $2$ along $Q$. 
Moreover, the restriction of $\Omega$ to $\partial {\mathcal Z}_+ =F$
is a real closed $3$-form which kills the tangent space of the foliation
$\mathscr F$, since it annihilates  ${\mathfrak w}_1$ and ${\mathfrak w}_2$;
thus  $\Omega|_F$ is actually the pull-back of a  real $3$-form on $P$. 
This shows that $\Omega$ descends to a continuous $3$-form on
$Z-Q$ which is holomorphic on the complement of $P$. 
It is therefore holomorphic even across  $P$, 
by an iterated application of the Weierstrass removable singularities theorem. 
Identifying $Z$ with $\CP_3$ as in Theorem \ref{twisthol},
 $\Omega$ thus becomes a meromorphic $3$-form on $\CP_3$
 with a double pole at a quadric $Q$, and its pull-back to
 the totally real submanifold $P\subset \CP_3-Q$ is 
 real, as promised. 
\end{proof}

  Analogy with Pontecorvo's 
  characterization \cite{mano} of the twistor spaces of positive-definite scalar-flat  K \"ahler
  metrics would lead one to expect that the converse is also true. 
  Fortunately, this is indeed the case:
    
  \begin{prop} \label{transf}
 Let $(M, [g])$ be a space-time-oriented  Zollfrei self-dual  $4$-manifold
 whose twistor space $(Z,J)$ is biholomorphic to $\CP_3$.
 Suppose that there is a quadric 
 $Q\subset Z\cong \CP_3$ such that $P\cap Q= \emptyset$,
 and that there is a  meromorphic $3$-form $\Omega$ on $Z$
 which is holomorphic and non-zero on $Z-Q$ and which has the property  that 
 its pull-back  to $P$ is a real $3$-form. Then $Q$ determines 
 an integrable  complex structure ${\mathfrak J}$ on $M$ such that 
 $(M,{\mathfrak J})\cong \CP_1\times \CP_1$, and the conformal class
 $[g]$ contains a scalar-flat metric $g$ which is indefinite K\"ahler with 
 respect to the complex structure $\mathfrak J$. Moreover, this metric
 is uniquely determined up to an overall  multiplicative constant. 
 \end{prop}
\begin{proof} 
Since $Q$ represents double the generator of $H^2 (\CP_3 , \ZZ)$,
it generates $H^2 (\CP_3 - P, \ZZ)= H^2 ({\mathcal Z}_+ , \ZZ)$,
and so has intersection number $1$ with any fiber disk in ${\mathcal Z}_+$.
Thus $Q$ represents a section of $\mbox{Int }{\mathcal Z}_+$,
and may be interpreted as an almost-complex structure $\mathfrak J$.
Moreover, the induced  projection $Q\to M$ is a diffeomorphism, and 
the pull-back of $\mathfrak J$ to $Q$ is exactly the given complex structure
on $Q\cong \CP_1\times \CP_1$, so $\mathfrak J$ is, in particular integrable.

Near an arbitrary point of $M$,
 choose a local pseudo-orthonormal frame so that $e_2= {\mathfrak J} e_1$
 and $e_4= {\mathfrak J}e_3$. Then $Q$ is represented in the corresponding
 local coordinates on ${\mathcal Z}_+$ by $\zeta =i$. 
Now pull
$\Omega$ back to ${\cal Z}_+$, and observe that 
we must then have 
$$\Omega= -   \frac{{\zap f}}{(1+\zeta^2)^2}
[(1+\zeta^2) \varphi_1 + (1-\zeta^2) \varphi_2 -2\zeta \varphi_3]\wedge
[2d\zeta - (1-\zeta^2) \theta_1^3 - 2\zeta \theta_1^2 + (1+\zeta^2) \theta_3^2]$$
for some function $\zap f$
on ${\mathcal Z}_+$, since $\Omega$  annihilates ${\mathfrak w}_1$, 
${\mathfrak w}_2$, and $\partial /\partial \overline{\zeta}$. 
Moreover, this function ${\zap f}$ is holomorphic in $\zeta \neq i$, 
bounded on the entire half-plane, and real when $\zeta$ is real. 
Hence $\zap f$ is independent of $\zeta$,  by Liouville's Theorem and the 
reflection principle.  In particular, the $\zeta$-derivative of $\zap f$ vanishes
at $\zeta = i$, so the residue 
 $\omega$ of $\Omega$ at $\zeta = i$  is a multiple of 
$\varphi_1$. However, this residue is also   a closed nowhere-zero
$2$-form on $M$, as, up to an overall constant,  it may instead be obtained by restricting $\Omega$
to $F= \partial {\mathcal Z}_+ = \Psi^{-1}(P)$ and integrating along the 
fibers of ${\zap p}: F\to M$. 
But this means that $\omega$ is the K\"ahler form with respect to ${\mathfrak J}$ 
of
 an indefinite  K\"ahler metric $g$ in the self-dual conformal class 
$[g]$. Since such a metric must also be scalar-flat by Lemma \ref{gau},
the claim therefore follows. 
\end{proof}

{\bf Theorem \ref{D}} now follows immediately from Propositions \ref{constr} and
\ref{transf}. since a projective transformation is all that is needed to arrange for 
the quadric $Q$ to be given by $z_1^2 + z_2^2+z_3^2+z_4^2=0$, and for 
the associated $3$-form to be some real constant times
$$
\Omega = \frac{\left(z_j\frac{\partial}{\partial z_j}\right) \hook (dz^1 \wedge dz^2 \wedge dz^3\wedge dz^4)}{[z_1^2 + z_2^2+z_3^2+z_4^2]^2}.
$$
Of course, requiring  that the pull-back of $\Omega$ to $P$ be real has been 
re-interpreted in the statement of 
{Theorem \ref{D}} as the condition that the pull-back of $\phi = \Im m ~\Omega$
should vanish.

It remains only to  ask whether there are many 
submanifolds $P$ near  $\RP^3\subset \CP_3$ on which 
$\phi= \Im m ~\Omega$ vanishes. In fact, the
condition in question is a weakening of the 
{\sl special Lagrangian} condition studied by 
McLean \cite{mclean}, and similar arguments 
will now   show that such 
submanifolds exist in considerable  profusion:

\begin{prop}
For any integer $k\geq 1$ and any $\alpha \in (0,1)$, the space
of compact $C^{k,\alpha}$ totally real submanifolds $P\subset \CP_3-Q$
near $\RP^3$ 
on which $\phi = \Im m ~ \Omega$ vanishes is a Banach manifold 
 whose tangent space at $P$ consists of real $C^{k,\alpha}$ vector
 fields $v$ on $P$ for which $\mbox{div }v =0$ with respect to the standard volume form 
on $\RP^3$. 
\end{prop}
\begin{proof}
The normal bundle of $\RP^3\subset \CP_3$ may be identified with $T\RP^3$
 via $J$, so some tubular neighborhood of $\RP^3$ must be diffeomorphic to
 $T\RP^3$. In fact, we can even take this tubular to be all of $\CP_3-Q$  by invoking the 
   real-analytic diffeomorphism 
 $$\gimel : T\RP^3\longrightarrow \CP^3-Q
 $$ 
 given  by $$\pm (\vec{x},\vec{y})\mapsto \Big[\vec{x}+ i\frac{\vec{y}}{\sqrt{1+|\vec{y}|^2}}\Big]$$
 for $\vec{x}, \vec{y}\in \RR^4$ with $|\vec{x}|^2=1$ and $\vec{x}\cdot \vec{y}=0$. 
Thus, for any 
 integer $k\geq 1$ and any $\alpha \in (0,1)$, each real-valued 
 $C^{k,\alpha}$ vector field $v$ on $\RP^3$ defines a new
 embedding 
 \begin{eqnarray*}
\saturn_v: \RP^3  &\longrightarrow & \CP_3- Q \\
y & \longmapsto & \gimel  (Jv_y)
\end{eqnarray*}
and every other compact submanifold of  $\CP_3$ which is
$C^{k,\alpha}$ close to $\RP^3$ can be so parameterized 
in a unique manner. Now let ${\mathfrak B}^{k,\alpha}$
be the  Banach space of 
$C^{k,\alpha}$ vector fields on $\RP^3$,  and  let ${\mathfrak C}^{k,\alpha}$
be the Banach space of $C^{k,\alpha}$ real-valued $3$-forms on $\RP^3$
with integral $0$ on $\RP^3$. Let $\mu$ be the standard volume form on 
$\RP^3$.  We may then define a smooth map of 
Banach manifolds 
 \begin{eqnarray*}
\leo : {\mathfrak B}^{k,\alpha}  \times {\mathfrak C}^{k,\alpha}&\longrightarrow &  {\mathfrak C}^{k-1,\alpha}\times 
{\mathfrak B}^{k-1,\alpha}\\
(v, f \mu ) & \longmapsto & (\saturn_v^*\phi , \mbox{curl }  v + \mbox{grad }f) 
\end{eqnarray*}
whose derivative at $0$  is 
$$(v,f)\mapsto ( \mbox{div } v, \mbox{curl } v+ \mbox{grad }f).$$
Now this is essentially just the elliptic operator $d+d^*: 
\Lambda^{\mbox{\tiny even}}\to \Lambda^{\mbox{\tiny odd}},$ 
and so has trivial kernel and cokernel because $H^2(\RP^3 , \RR)=0$. 
The interior 
Schauder estimates for elliptic equations therefore imply that $\leo_{*0}$
is a Banach-space isomorphism. Hence the  inverse function
theorem for Banach spaces implies  that $\leo$ becomes a diffeomorphism 
when restricted to
 some neighborhood ${\mathcal U}={\mathcal U}_1\times {\mathcal U}_2\subset 
 {\mathfrak B}^{k,\alpha}  \times {\mathfrak C}^{k,\alpha}$ of the origin. 
 Thus $(\leo|_{\mathcal U})^{-1} (  \{0\}\times {\mathfrak B}^{k-1,\alpha})$  
is  a Banach manifold. By inspection, however, this set is 
of the form ${\mathcal M}\times {\mathcal U}_2$, where $\mathcal M$
is its projection to ${\mathfrak B}^{k,\alpha}$. Hence $\mathcal M$
is a Banach manifold, and  represents the desired moduli space of 
solutions. Moreover,
$$T_0{\mathcal M} = \{ v \in {\mathfrak B}^{k,\alpha}\subset \Gamma (T\RP^3)~|~ \mbox{div } v=0\},$$
precisely as claimed, so we are done. 
 \end{proof}
 
 Thus, 
while  self-dual split-signature
conformal structures on $S^2\times S^2$ essentially depend on  a 
 vector field on $\RP^3$, 
scalar-flat K\"ahler metrics correspond to the case in which  the  vector field
is {\sl divergence free}.
So far, though, this is just a abstract existence statement.
 Nonetheless,  one can do much better in the real-analytic case. 
 Indeed, let $v$ be a divergence-free real-analytic vector field  on $\RP^3$;
 in other words, let  $v= \mbox{curl }w$ for 
 $w$ some real-analytic vector field on $\RP^3$. Then $Jv$ corresponds
 to the section $Jv+iv$ of  $(T^{1,0}\CP_3)|_{\RP^3}$.
 Because $v$ is locally represented by 
power series, $Jv+iv$ can then be uniquely extended to some neighborhood $U\subset \CP_3-Q$ 
of $\RP^3$ as a holomorphic  vector field ${\mathfrak v}$, and we then have
a real-analytic $1$-parameter family $\{ \psi_t~|~ t\in (-\varepsilon, \varepsilon )\}$ of biholomorphisms from  neighborhoods $U_t\subset U$ 
of
$\RP^3$ to $U$ obtained  by following the integral curves 
  of $\Re e~ {\mathfrak v}$. 
Notice, however, that ${\mathcal L}_{\mathfrak  v}\Omega=0$, since this expression 
is the analytic continuation of  $\mbox{div } iv$ from  $\RP^3$ to $U$. 
The constructed biholomorphisms therefore satisfy $\psi_t^*\Omega =\Omega$. 
Hence $P= \psi_t (\RP^3)$ is a submanifold on which $\phi = \Im m ~ \Omega$
vanishes.

\vfill

\begin{ack} 
The first author would like to  express his gratitude to  Franc Forstneri\v{c}, Bill Goldman, 
Denny Hill,    Matthias Kreck, Blaine Lawson, 
Yair Minsky, Yom-Tung Siu, and Dennis Sullivan 
for their friendly help in 
 drawing his attention to some key  references. 
He would also like to thank  to  Jeff Cheeger and the 
 Courant Institute of Mathematics  
 for their  hospitality during the initial phase  of the writing of this paper. 
\end{ack}

\hspace{1in}

\noindent 
{\sc 
Department of Mathematics, SUNY, Stony Brook, NY 11794-3651 USA\\
The Mathematical Institute, 24-29 St Giles,
Oxford OX1 3LB,  England}

\pagebreak 
%
%\bibliographystyle{siam}
%\bibliography{lebrun}

  \end{document}